\documentclass{nyjm}
\usepackage{amsfonts,amsthm,amsmath,amssymb,mathtools}
\usepackage[all,cmtip]{xy}
\usepackage[english]{babel}
\usepackage[utf8]{inputenc}
\usepackage[OT1]{fontenc}
\usepackage{mathrsfs}
\usepackage{xpunctuate}
\usepackage{enumitem}
\usepackage[backref]{hyperref}
\usepackage{geometry}
\numberwithin{equation}{section}

\theoremstyle{plain}
\newtheorem{theorem}{Theorem}[section]
\newtheorem*{theorem*}{Theorem}
\newtheorem*{conjecture}{Conjecture}
\newtheorem{proposition}[theorem]{Proposition}
\newtheorem*{defthm*}{Definition/Theorem}
\newtheorem{lemma}[theorem]{Lemma}
\newtheorem*{lemma*}{Lemma}
\newtheorem{corollary}[theorem]{Corollary}
\makeatletter
\newtheorem*{intr@thm}{\intr@thmname}
\newenvironment{introtheorem}[1]{%
	\def\intr@thmname{Theorem #1}
	\begin{intr@thm}}
	{\end{intr@thm}}
\newtheorem*{intr@crl}{\intr@crlname}
\newenvironment{introcorollary}[1]{%
	\def\intr@crlname{Corollary #1}
	\begin{intr@crl}}
	{\end{intr@crl}}
\newtheorem*{intrd@f}{\intrd@fname}
\newenvironment{introdefinition}[1]{%
	\def\intrd@fname{Definition #1}
	\begin{intrd@f}}
	{\end{intrd@f}}
\makeatother

\theoremstyle{remark}
\newtheorem{remark}[theorem]{Remark}
\newtheorem{example}[theorem]{Example}

\theoremstyle{definition}
\newtheorem{definition}[theorem]{Definition}
\newtheorem*{definition*}{Definition}

\newtheorem{notation}[theorem]{Notation}

\newcommand{\expp}[1]{p^{#1}}
\newcommand*{\Gal}[2]{\mathrm{Gal}(#1/#2)} 
\newcommand*{\splitG}[1][0]{D_{#1}} 
\newcommand*{\filindex}[3][,]{h_{#2#1#3}}
\newcommand*{\subsplit}[1][0]{\Sigma_{#1}} 
\newcommand*{\catMod}[2][{[\splitG]}]{\mathbf{Mod}^{\mathrm{#2}}_{\ZZ_p#1}} 
\newcommand*{\cofg}{co\textrm{-}f.g.} 
\newcommand*{\ZZ}{\mathbb{Z}} 
\newcommand*{\NN}{\mathbb{N}} 
\newcommand*{\QQ}{\mathbb{Q}} 
\DeclareMathOperator{\pr}{pr} 

\makeatletter 
\DeclareRobustCommand{\cev}[1]{%
  {\mathpalette\do@cev{#1}}%
}
\newcommand{\do@cev}[2]{%
  \vbox{\offinterlineskip
    \sbox\z@{$\m@th#1 x$}%
    \ialign{##\cr
      \hidewidth\reflectbox{$\m@th#1\vec{}\mkern4mu$}\hidewidth\cr
      \noalign{\kern-\ht\z@}
      $\m@th#1#2$\cr
    }%
  }%
}
\makeatother

\newcommand{\rint}[1]{\mathcal{O}_{#1}} 
\newcommand{\unit}[1]{\rint{#1}^\times} 
\newcommand*{\realemb}[1]{r_1(#1)} 
\newcommand*{\cpxemb}[1]{r_2(#1)} 
\newcommand{\clgp}[1]{Cl_{#1}} 
\newcommand*{\clsyl}[1]{A_{#1}} 
\newcommand*{\hsyl}[1]{{\varepsilon_{#1}}} 
\newcommand*{\algnm}[1]{\mathrm{Nm}_{#1}}
\newcommand*{\arnm}[2]{\mathrm{N}_{#1/#2}}
\newcommand*{\arextsymb}{\iota} 
\newcommand*{\arext}[2]{\arextsymb_{#1/#2}}
\newcommand*{\princ}[1]{\mathrm{Pr}_{#1}} 
\newcommand*{\ideal}[1]{\mathrm{Id}_{#1}} 
\newcommand*{\aboveT}[1]{\Pi_{#1}} 
\newcommand*{\Ram}[2]{\operatorname{T}(#1/#2)} 
\newcommand*{\oram}[2]{\tau(#1/#2)} 
\newcommand*{\idele}[1]{\mathbb{A}^\times_{#1}} 
\newcommand*{\idunit}[1]{\mathcal{U}_{#1}} 
\newcommand*{\idclg}[1]{\mathcal{C}_{#1}} 
\newcommand*{\q}[1]{Q_{#1}} 
\def\primen@me{p}
\def\conjprimen@me{q}
\def\oth@rprime{l}
\newcommand*{\primeid}{\mathfrak{\primen@me}}
\newcommand*{\conjprimeid}{\mathfrak{\conjprimen@me}}
\newcommand*{\otherprime}{\mathfrak{\oth@rprime}}
\newcommand*{\Primeid}[1]{\mathfrak{\MakeUppercase{\primen@me}}_{#1}}
\newcommand*{\Conjprimeid}[1]{\mathfrak{\MakeUppercase{\conjprimen@me}}_{#1}}
\newcommand*{\Otherprime}[1][\empty]{\mathfrak{\MakeUppercase{\oth@rprime}}_{#1}}
\newcommand*{\semiloc}[2][\empty]{\mathcal{V}^{#1}_{#2}} 
\newcommand*{\ramdeg}[2]{e_{#1/#2}} 
\newcommand*{\indeg}[3][/]{f_{#2#1#3}} 
\newcommand*{\logindeg}[1][\empty]{\ell_\infty^{#1}} 
\newcommand*{\locdeg}[2]{d_{#1/#2}} 
\newcommand*{\nprimes}[1][\primeid]{g_{#1}} 
\newcommand*{\fundclass}[2]{u_{#1/#2}} 
\newcommand*{\loclev}[1]{E_{#1}} 
\newcommand*{\invJP}[1][\loclev{\empty}]{\operatorname{inv}_{#1}} 

\newcommand{\cch}[1]{\delta^{(#1)}} 
\newcommand{\cchT}[1]{\widehat{\delta}^{(#1)}} 
\newcommand*{\tateiso}[3][\empty]{(\cup\varkappa_{#2}^{#1})^{#3}} 
\DeclarePairedDelimiterXPP{\hh}[3]{\widehat{h}^{#1}}{(}{)}{}{#2,#3} 
\DeclarePairedDelimiterXPP{\hhN}[3]{h^{#1}}{(}{)}{}{#2,#3} 
\DeclarePairedDelimiterXPP{\HH}[3]{\widehat{H}^{#1}}{(}{)}{}{#2,#3} 
\DeclarePairedDelimiterXPP{\HHm}[2]{\widehat{H}^{#1}}{(}{)}{}{#2} 
\DeclarePairedDelimiterXPP{\HHN}[3]{H^{#1}}{(}{)}{}{#2,#3} 
\DeclareMathOperator{\res}{Res} 
\DeclareMathOperator{\defl}{Defl} 
\DeclareMathOperator{\infl}{Inf} 

\newcommand*{\pallino}[1][{\field}]{{#1{\bullet}}} 
\newcommand*{\nsu}{j} 
\newcommand*{\nsd}{k} 
\newcommand*{\normsys}{\mathbf{NS}_{\DD}} 
\newcommand*{\hatnsu}[2][\jmath]{\widehat{#1}_{#2}} 
\makeatletter
\DeclarePairedDelimiterXPP{\s@sHH}[2]{\p@rH^{#1}}{(}{)}{}{#2} 
\newcommand{\sysHH}[3][\empty]{%
\def\p@rH{\widehat{\calH}_{#1}}
\s@sHH{#2}{#3}%
}
\makeatother

\newcommand*{\sysnorm}[1]{\overline{\mathcal{N}#1}} 

\makeatletter
\def\r@otsl#1#2#3#4{\mu_{#4}(#2)}
\def\r@ots#1{\mu_{#1}}
\def\roots{\futurelet\radici\r@otschose}
\def\r@otschose{%
	\ifx\radici[
		\expandafter\r@otsl
	\else
		\expandafter\r@ots
\fi}
\makeatother 
\newcommand{\defect}[2]{\beta_{#1}} 

\newcommand*{\iwLambda}[1][]{\Lambda#1} 
\newcommand*{\topgen}{\gamma_0} 
\newcommand*{\iwmu}[2][]{\mu_{#2}^{#1}} 
\newcommand*{\iwlambda}[2][]{\lambda_{#2}^{#1}} 
\newcommand*{\lambdacohom}[2][\unit{\pallino}]{\lambda(#2,#1)}
\newcommand*{\iwnu}[2][]{\nu_{#2#1}} 
\makeatletter
\def\th@X{X} 
\newcommand*{\iwXL}{\th@X_{\iw}} 
\newcommand*{\iwXK}{\th@X_{\fk}} 
\makeatother
\newcommand*{\indextr}[1][\fromkprimeid]{{n_\mathrm{tr}^{#1}}} 
\newcommand*{\Lcyc}{L_\mathrm{cyc}} 
\newcommand*{\Ldue}{\widetilde{L}} 
\newcommand*{\DD}{\mathcal{D}} 
\newcommand*{\GG}{\mathcal{G}} 
\newcommand*{\Minfty}{M_\infty} 
\newcommand*{\Frob}[2]{\operatorname{Frob}(#1,#2)} 

\newcommand*{\nuceta}[1]{\eta_{#1}} 

\DeclareMathOperator{\rank}{rk} 
\DeclareMathOperator{\corank}{cork} 
\DeclarePairedDelimiterX\gindex[2]{(}{)}{#1:#2} 
\DeclarePairedDelimiter\gorder{\lvert}{\rvert} 
\newcommand*{\kernel}[1]{\operatorname{Ker} #1} 
\newcommand*{\cokernel}[1]{\operatorname{Coker} #1} 
\newcommand*{\im}[1]{\operatorname{Im} #1}
\newcommand*{\tors}[2][p]{{{#2}[#1]}}
\newcommand*{\id}[1]{\mathrm{id}_{#1}} 
\DeclareMathOperator{\Hom}{Hom}
\newcommand*{\sch}{\delta} 
\makeatletter
\newcommand{\Ind}[2]{%
\def\gr@up{#1}%
\def\subgr@up{#2}%
\operatorname{Ind}_{\gr@up}^{\subgr@up} 
					} 
\makeatother

\makeatletter
\newcommand*{\dblsys}[2][\empty]{\mathbf{DS}^{\mathrm{#1}}_{#2}}
\newcommand*{\dsu}{\varphi} 
\newcommand*{\dsd}{\pi}  
\newcommand*{\capzero}[2][\scriptscriptstyle\vartriangle]{{^{#1}{#2}}} 
\newcommand*{\capbar}[2][\scriptscriptstyle\blacktriangle]{{_{#1}{#2}}} 
\newcommand*{\nuczero}[2][\dagger]{{^{#1}{#2}}} 
\newcommand*{\nucbar}[2][\ddagger]{{_{#1}{#2}}} 
\newcommand*{\bddsys}[2][\empty]{\mathbf{B}^{\mathrm{#1}}_{#2}} 
\newcommand*{\quotsys}[2][\empty]{{\dblsys[#1]{#2}}/{\bddsys[#1]{#2}}} 
\newcommand*{\natr}{\Phi} 
\newcommand*{\bddcong}[1][\bddsys{\DD}]{\cong_{#1}} 
\newcommand*{\bddeq}{\propto} 
\newcommand*{\longdashrightarrow}{\dashrightarrow} 
\makeatother
\newcommand*{\stepfunctor}{\mathscr{S}}
\newcommand*{\limemb}{\mathscr{L}} 
\newcommand*{\calH}{\mathcal{H}}
\newcommand*{\calX}{\mathcal{X}}
\newcommand*{\calY}{\mathcal{Y}}
\newcommand*{\calZ}{\mathcal{Z}}
\newcommand*{\calB}{\mathcal{B}}
\newcommand*{\calC}{\mathcal{C}}

\newcommand*{\longhookrightarrow}{\lhook\joinrel\relbar\joinrel\rightarrow} 
\newcommand*{\longhookdashedarrow}{\lhook\dashrightarrow} 
\newcommand*{\longtwoheaddashedarrow}{\dashrightarrow\mathrel{\mkern-19mu}\dashrightarrow}
\newcommand{\ie}{\textit{i.~e.~}} 
\newcommand*{\rsp}{{resp.~}} 
\newcommand{\loccit}{\emph{loc.\kern3pt cit}\xperiod} 
\newcommand{\ibid}{\emph{ibid}\xperiod} 
\newcommand{\fk}{\mathrm{fake}} 
\newcommand{\cyc}{\mathrm{cyc}} 
\newcommand{\iw}{\mathrm{Iw}} 
\title{Cohomology of normic systems and fake $\ZZ_p$-extensions}
\author[L.~Caputo]{Luca~Caputo}
\email{luca.caputo@gmx.com}
\author[F.~A.~E.~Nuccio]{Filippo A.~E.~Nuccio Mortarino~Majno~di~Capriglio}
\email{filippo.nuccio@univ-st-etienne.fr}
\address{Université Jean Monnet Saint-Étienne, CNRS UMR 5208, Institut Camille Jordan, F-42023 Saint-Étienne, France}

\date{November 2$^\text{nd}$, 2023}
\subjclass{Primary 11R23, 11R20; Secondary 11R29, 11R34}
\keywords{Iwasawa theory, dihedral Galois extension, normic systems, units in profinite extensions}
\begin{document}
\begin{abstract}

We set up a general framework to study Tate cohomology groups of Galois modules along $\ZZ_p$-extensions of number fields. Under suitable assumptions on the Galois modules, we establish the existence of a five-term exact sequence in a certain quotient category whose objects are simultaneously direct and inverse systems, subject to some compatibility. The exact sequence allows one, in particular, to control the behaviour of the Tate cohomology groups of the units along $\ZZ_p$-extensions.

As an application, we study the growth of class numbers along what we call ``fake $\ZZ_p$-extensions of dihedral type''. This study relies on a previous work, where we established a class number formula for dihedral extensions in terms of the cohomology groups of the units.

\end{abstract}
\maketitle
\section{Introduction}
The main goal of this work is to set up a convenient algebraic framework to study Tate cohomology groups along $\ZZ_p$-extensions of number fields. Before giving more details, let us fix some notation. Let $p$ be a prime number, let $F$ be number field and fix a $\ZZ_p$-extension $L_\infty/F$, by which we mean a Galois extension such that $\Gamma=\Gal{L_\infty}{F}$ is isomorphic to $\ZZ_p$: in particular,
\[
L_\infty=\bigcup_{n\geq 0} L_n
\]
where we set $L_0=F$ and where each $L_n/F$ is a cyclic Galois extension of degree $p^n$. Write $\Gamma_n$ for the open subgroup $\Gamma_n=\Gamma^{p^n}=\Gal{L_\infty}{L_n}$ and set $G_n=\Gamma/\Gamma_n=\Gal{L_n}{F}$; more generally, for all $m\geq n\geq 0$, set  $G_{m,n}=\Gal{L_m}{L_n}$. In this setting, one can attach to each field $L_n$, often regarded as a ``layer'', several interesting arithmetic objects: the unit group $\unit{L_n}$, the ideal class group $\clgp{L_n}$ or its $p$-Sylow subgroup $\clsyl{L_n}$, the group $\idunit{L_n}\subseteq \idele{L_n}$ of idelic units, the idèle class group $\idclg{L_n}$, and so forth. Let $\{B_n\}_{n\geq 0}$ denote any of the above collections. Since all the $B_n$'s are $G_n$-modules, the Tate cohomology groups $\HH{i}{G_n}{B_n}$ are defined, for every $i\in\ZZ$. Moreover, for all $m\geq n\geq 0$, the inclusion $L_n^\times\hookrightarrow L_m^\times $ and the norm $\arnm{L_m}{L_n}\colon L_m^\times\to L_n^\times$ induce $G_{m}$-morphisms $\arext{L_m}{L_n}\colon B_n\to B_m$ and $\arnm{L_m}{L_n}\colon B_m\to B_n$. Almost by definition, these maps satisfy 
\begin{equation}\label{eq:norm_relations}
\arnm{L_m}{L_n}\circ\arext{L_m}{L_n}=p^{m-n}\qquad\text{ and }\qquad\arext{L_m}{L_n}\circ \arnm{L_m}{L_n}=\algnm{G_{m,n}}
\end{equation}
where, for an arbitrary group $H$, we denote by $\algnm{H}$ the norm element in $\ZZ[H]$. These morphisms are fundamental in the study of the behaviour of the groups $B_n$ along the $\ZZ_p$-extension but they do not systematically induce maps $\dsu_{n,m}\colon\HH{i}{G_n}{B_n}\to\HH{i}{G_m}{B_m}$ or $\dsd_{m,n}\colon \HH{i}{G_m}{B_m}\to\HH{i}{G_n}{B_n}$. Indeed, although the maps
\[
\arext{L_m}{L_n}^*\colon\HH{i}{G_n}{B_n}\longrightarrow\HH{i}{G_n}{B_m^{G_{m,n}}}\qquad\text{ and }\qquad
\arnm{L_m}{L_n}^*\colon\HH{i}{G_n}{B_m^{G_{m,n}}}\longrightarrow\HH{i}{G_n}{B_n}
\]
are always defined, they do not have the expected domain or codomain. When $i\geq 1$, one could take $\dsu_{m,n}=\infl\circ\arext{L_m}{L_n}^*$, where $\infl$ denotes the inflation; and, when $i\leq -1$, one could consider the deflation $\defl$ (see \cite{Wei59}), setting $\dsd_{m,n}=\arnm{L_m}{L_n}^*\circ\defl$. But, in general, for a given $i\in\ZZ$, only one of these would be defined. In particular, a relation like~\eqref{eq:norm_relations} involving $\dsu_{n,m}$ and $\dsd_{n,m}$ could not even be stated. Moreover, taking $G_n$-cohomology kills the $G_n$-action, so all abelian groups $\HH{i}{G_n}{B_n}$ are trivial $G_n$-modules. 

On the other hand, these Tate cohomology groups are certainly interesting arithmetic objects: they have the advantage of always being \emph{finite} (at least for all $B_n$ as above) and they occur naturally in class field theory. For example,
\[
\kernel\bigl(\HH{1}{G_n}{\unit{L_n}}\longrightarrow\HH{1}{G_n}{\idunit{L_n}}\bigr)\cong \kernel\bigl(\arext{L_n}{F}\colon \clgp{F}\to \clgp{L_n}\bigr)\qquad\text{ for all } n\geq 0
\]
(see, for instance,~\cite[Proposition~2.2]{Nuc10}) or, even more fundamentally,
\[
\HH{2}{\Gal{H_n}{L_n}}{\idclg{H_n}}\cong \clgp{L_n}\qquad\text{ for all }n\geq 0
\]
where $H_n$ denotes the Hilbert class field of $L_n$. To give one more example, one that has been at the origin of our investigation, assume that $p$ is odd and suppose that there exists a subfield $k\subseteq F$ such that, for every $n$, $L_n/k$ is Galois with dihedral Galois group $\splitG[n]$ of order $2p^n$. Denote by $K_n$ a subfield of index $2$ in $L_n/k$: then the equality
\begin{equation}\label{eq:from_glasgow}
\frac{\gorder{\clgp{L_n}}\cdot \gorder{\clgp{k}}^2}{\gorder{\clgp{F}}\cdot \gorder{\clgp{K_n}}^2}=\frac{\gorder{\HH{0}{\splitG[n]}{\unit{L_n}}}}{\gorder{\HH{-1}{\splitG[n]}{\unit{L_n}}}}
\end{equation}
holds up to a power of $2$ (see~\cite[Theorem~3.14]{CapNuc20}).

Hence, it looks compelling to regard the assignment $n\mapsto \HH{i}{G_n}{B_n}$ as an analogue of $n\mapsto B_n$, or perhaps of $n\mapsto \clsyl{L_n}$, since the $\clsyl{L_n}$ are also finite groups. Yet, this seems to immediately break down, at least insofar techniques from Iwasawa theory are involved: to explain why, it is probably necessary to clarify what we mean by ``studying'' or ``analysing'' the above assignment. As mentioned, the groups $\HH{i}{G_n}{B_n}$ are all finite, and the most basic question one could ask is to describe the behaviour of their orders $\gorder{\HH{i}{G_n}{B_n}}$ as $n\to +\infty$. After all, the whole subject of Iwasawa theory begun with the celebrated
\begin{theorem*}[Iwasawa] There exists three integers $\iwmu{L_\infty},\iwlambda{L_\infty},\iwnu{L_\infty}$ and an index $n_0$ such that
	\begin{equation}\label{eq:original_iwasawa}
	\gorder{\clsyl{L_n}}=p^{\iwmu{k_\infty}p^n+\iwlambda{k_\infty}n+\iwnu{k_\infty}}\qquad\text{ for all }n\geq n_0.
	\end{equation}
\end{theorem*}
Recall now the usual strategy of proof of Iwasawa's theorem: one first considers the inverse limit, with respect to the norm maps, $\iwXL=\varprojlim \clsyl{L_n}$ and regards it as a module over the completed group algebra $\iwLambda=\ZZ_p[\![\Gamma]\!]$, showing that as such it is finitely generated and torsion. Secondly, a fine analysis of the Galois action on $\iwXL$, coupled with global class field theory, shows that one can closely relate $\clsyl{L_n}$ with the co-invariants $(\iwXL)_{\Gamma_n}$; finally, a structure theorem for $\iwLambda$-torsion modules of finite type yields the formula. Inherent to this approach is the asymptotic flavour of Iwasawa's result: at several stages, some adjustment is required, which modifies the outcome by a ``finite, bounded error term''. Now, all this breaks down when replacing the finite groups $\clsyl{L_n}$ with any of the $\HH{i}{G_n}{B_n}$: already, they might not form an inverse system, when $i\geq 1$; and, even if they do, the triviality of the Galois action implies that the co-invariants of the inverse limit coincide with the whole limit and are, in general, infinite. Hence, they are of no use to recover the finite groups $\HH{i}{G_n}{B_n}$.

To describe our approach, it is important to consider the somewhat dual strategy of considering the direct limit
\[
A_{L_\infty}=\varinjlim_{\arext{L_m}{L_n}}\clsyl{L_n}.
\]
As discussed in~\cite[Introduction]{Gre73} the modules $\iwXL$ and the Pontryagin dual $\Hom_{\ZZ_p}(A_{L_\infty},\QQ_p/\ZZ_p)$ essentially carry the same information, and could be used interchangeably; Iwasawa himself occasionally works with the latter module instead of the former (see, for instance,~\cite{Iwa81}). When $B_n=\unit{L_n}$ rather than $\clsyl{L_n}$, the direct limit occurs, for instance, in~\cite[\S5]{Iwa83}: there, Iwasawa claims that the inflation maps
\begin{equation}\label{eq:iwasawa_limit}
\infl\colon\HH{i}{G_n}{\unit{L_n}}\longrightarrow \tors[p^n]{\bigl(\varinjlim\HH{i}{G_m}{\unit{L_m}}\bigr)}=\tors[p^n]{\HHN{i}{\Gamma}{\unit{L_\infty}}}\qquad\text{ for }i=1,2
\end{equation}
have kernel and cokernel which have bounded order as $n$ grows. Given that the structure of the direct limits $\varinjlim\HH{i}{G_n}{B_n}=\HH{i}{\Gamma}{B_\infty}$ is, in general, quite explicit (but ``up to finite groups'') we interpret this as a weak analogue of Iwasawa's descent theorem leading to~\eqref{eq:original_iwasawa}, replacing the operation of taking co-invariants by taking the $p^n$-torsion subgroup. Hence, we set out to find an algebraic setting where the boundedness of kernels and cokernels as in~\eqref{eq:iwasawa_limit} could be proven in general.
But we were still confronted with two problems: first, in some formul\ae, for instance in~\eqref{eq:from_glasgow}, negative cohomological degrees must be considered and in this case the Tate cohomology groups naturally form an inverse, rather than a direct, system; and secondly, that working ``up to finite groups'', would turn all inflation maps occurring in~\eqref{eq:iwasawa_limit} into the $0$ map.

At this point, we were inspired by Vauclair's approach in~\cite{Vau09}, where he defines \emph{normic systems} in quite a general context: in particular, all the examples $\{B_n\}_{n\geq 0}$ above are normic systems. Rather than taking limits, he works in the category whose objects are collections of $G_n$-modules simultaneously carrying the structure of a direct and an inverse system, subject to some compatibility. With this in mind, we define, at least under the assumptions~\ref{cond:inj} and~\ref{cond:surj} of Definition~\ref{def:conditions_Inj_Gal} on the normic system $\{B_n\}_{n\geq 0}$, ascending and descending morphisms
\[
\xymatrix@C=3em{
	\HH{i}{G_m}{B_m}\ar@<.35em>[0,2]^{\dsd_{m,n}}&&\HH{i}{G_n}{B_n}\ar@<.35em>[0,-2]^{\dsu_{n,m}}&\text{for all }m\geq n\geq 0}
\]
turning the collection of abelian groups $\HH{i}{G_n}{B_n}$ in what we call a ``double system''. Inside the category\footnote{In this introduction, we confine ourself to a special case of our main result, which holds for profinite groups $\DD$ which are more general than $\Gamma$. Although this generalised approach is crucial for the dihedral applications we have in mind, we prefer to stick to the case $\DD=\Gamma$ here, to lighten notation. In particular, all subscripts $\Gamma$ decorating the several categories that we are going to introduce, actually read as $\DD$ in the body of the paper.} $\dblsys{\Gamma}$ of double systems, we identify a a certain thick subcategory $\bddsys{\Gamma}$ of double systems of bounded orders (see Definition~\ref{def:bddsys}), such that the corresponding quotient category $\quotsys{\Gamma}$ turns out to be the framework we were looking for. Indeed, upon restricting to a subcategory $\quotsys[\cofg]{\Gamma}$ of $\quotsys{\Gamma}$ defined by some very natural co-finiteness condition of the direct limit, we define an endofunctor
\[
\limemb\colon\quotsys[\cofg]{\Gamma}\longrightarrow \quotsys[\cofg]{\Gamma},\qquad\calX=(X_n)_{n\geq 0}\longmapsto(X_\infty[p^n])_{n\geq 0}
\]
attaching to a double system $\calX=(X_n)_{n\geq 0}$ the double system $(X_\infty[p^n])_{n\geq 0}$ (endowed with suitable transition morphisms), where $X_\infty=\displaystyle{\varinjlim X_n}$ . Writing $a_n \bddeq b_n$ to mean that two sequences $\{a_n\},\{b_n\}$ of natural numbers are eventually proportional (see Definition~\ref{def:eventually_prop}), we obtain the following
\begin{introcorollary}{\ref{cor:ha_vinto_caputo}}
Let $\calX=(X_n)_{n\geq 0}$ be a double system in $\dblsys[\cofg]{\Gamma}$, and suppose that $\calX\bddcong[\bddsys{\Gamma}]\limemb(\calX)$. Then
\[
\gorder{X_n}\bddeq p^{\lambda_\calX n}
\]
where $\varinjlim\calX\cong (\QQ_p/\ZZ_p)^{\lambda_\calX}\oplus\text{\textup{(finite group)}}$.
\end{introcorollary}
In this language, we can reinterpret Iwasawa's claim about the boundedness of the kernels and cokernels in~\eqref{eq:iwasawa_limit} as the statement that the double systems $\sysHH{i}{\unit{\pallino}}$, whose components are the cohomology groups $\HH{i}{G_n}{\idunit{L_n}}$, satisfy $\sysHH{i}{\unit{\pallino}}\bddcong\limemb(\sysHH{i}{\unit{\pallino}})$ for $i=1,2$. This is proved in Theorem~\ref{thm:i_cinque_dell'apocalisse} (see also Remark~\ref{rmk:Iwasawa83}).

It is worth noting that the rigidity yielded by requiring that the objects in the $\dblsys{\Gamma}$ are simultaneously inverse and direct systems is forced upon us by an interesting phenomenon: during the proof of Theorem~\ref{thm:i_cinque_dell'apocalisse}, two morphisms in $\dblsys{\Gamma}$ need to be shown to become injective in $\quotsys{\Gamma}$, in order to perform some homological algebra. Now, the analogous definitions $\bddsys[\mathrm{inv}]{\Gamma}$ (\rsp~$\bddsys[^\mathrm{dir}]{\Gamma}$) of the thick subcategory $\bddsys{\Gamma}$ can be given for the category of inverse (\rsp~direct) systems. Yet, it turns out that the kernel of the first morphism (regarded simply as an inverse system) does not belong to $\bddsys[\mathrm{inv}]{\Gamma}$, and the kernel of the second one (regarded simply as a direct system) does not belong to $\bddsys[\mathrm{dir}]{\Gamma}$: so, working either with direct or with inverse system alone would break our homological argument. On the other hand, since both kernels belong to $\bddsys{\Gamma}$, both $\dblsys{\Gamma}$-morphisms become injective in $\quotsys{\Gamma}$, as wanted: for a more comprehensive analysis of this phenomenon, we refer to Remark~\ref{rmk:quel_vecchio_volpone_del_Nuccio}. We explicitly mention Yamashita's paper~\cite{Yam84}, which has been very helpful at this point: she regards Iwasawa's boundedness claim concerning~\eqref{eq:iwasawa_limit} in the setting of abelian groups endowed with two (not necessarily compatible) structures of direct and inverse systems, and this was the starting point for our definitions of $\dblsys{\Gamma}$ and $\bddsys{\Gamma}$.

With this formalism at our disposal, we are in shape to ``let $n$ go to $+\infty$'' in~\eqref{eq:from_glasgow}. To state our result, as well as to motivate our title, we record the
\begin{introdefinition}{\ref{def:fake}} Let $p$ be an odd prime. Suppose that there exists a subfield $k\subseteq F$ such that for every $n$, $L_n/k$ is Galois, with dihedral Galois group $\splitG[n]$ of order $2p^n$. Denote by $K_n$ a subfield of index $2$ in $L_n/k$, chosen so that $K_n\supseteq K_{n-1}$ for all $n\geq 1$, and put $K_\infty = \bigcup K_n$. The extension $K_\infty/k$ is said to be a fake $\ZZ_p$-extension of dihedral type, the extensions $L_\infty/k$ is the Galois closure of the fake $\ZZ_p$-extension and the field $F$ is said to be the normalizing quadratic extension.
\end{introdefinition}
With this definition, one of our main results is the following
\begin{introtheorem}{\ref{thm:iwasawa_formula}}
Let $K_\infty/k$ be a fake $\ZZ_p$-extension of dihedral type. There exist constants $\iwmu{\fk},\iwnu{\fk}\in\ZZ[\frac{1}{2}]$ and $\iwlambda{\fk}\in\ZZ$ such that
\begin{equation*}
\gorder{\clsyl{K_n}}=p^{\iwmu{\fk}p^n+\iwlambda{\fk}n+\iwnu{\fk}}\qquad\text{ for all }n\gg 0.
\end{equation*}
\end{introtheorem}
Observe that in a fake $\ZZ_p$-extension of dihedral type, none of the subextensions $K_n/k$ is normal, so the groups $\clsyl{K_n}$ are not Galois modules and Iwasawa's original formula do not apply. We refer to~\S\ref{subsec:formula} for a more precise version of the above statement, describing the invariants $\iwmu{\fk}$ and $\iwlambda{\fk}$ in terms of explicit arithmetic quantities. In \S\ref{subsec:special_cases} we derive, in some special cases, explicit bounds for the values taken by the invariant $\iwlambda{\fk}$: as an example, let us state the following result\footnote{In its actual formulation, the quoted corollary is more precise: in order to avoid too much notation, we content ourselves with a slightly weaker statement in this Introduction, referring to~\S\ref{subsec:special_cases} for the full statement.}
\begin{introcorollary}{\ref{cor:formula_quad}}
Let $K_\infty/\QQ$ be a fake $\ZZ_p$-extension of dihedral type: in particular, the normalizing extension $F/\QQ$ is imaginary quadratic and the Galois closure $L_\infty/\QQ$ is the anticyclotomic extension of $F$. Then,
\[
\iwlambda{\fk}=
\begin{cases}
\frac{\iwlambda{\iw}+1}{2}&\text{ if $p$ splits in $F/\QQ$}\\&\\
\frac{\iwlambda{\iw}}{2}&\text{ if $p$ does not split in $F/\QQ$}
\end{cases}
\]
and therefore $\iwlambda{\iw}$ is odd if $p$ splits in $F$ and it is even if $p$ does not split.
\end{introcorollary}
Corollaries~\ref{cor:special_case_inert} and~\ref{cor:Ass+CM}  extend the above one to more general CM fields beyond the imaginary quadratic case.

In the Galois setting of $\ZZ_p$-extensions, the invariant $\iwlambda{\iw}$ responsible for the linear growth in the exponent of~$\clsyl{L_n}$ can be interpreted as the $\QQ_p$-dimension of $\iwXL\otimes\QQ_p$, where $\iwXL=\varprojlim\clsyl{L_n}$. In the final Section~\S\ref{subsec:structure_limit} we extend this result to our non-Galois setting, by studying the structure of the projective limit
\[
\iwXK=\varprojlim\clsyl{K_n}
\]
with respect to norm maps. We prove the following
\begin{introtheorem}{\ref{thm:lambda_equal_dim}}
	Given a fake $\ZZ_p$-extension of dihedral type $K_\infty/k$, we have $\iwlambda{\fk}=\dim_{\QQ_p} \iwXK\otimes_{\ZZ_p} \QQ_p$.
\end{introtheorem}
We conclude this Introduction by observing that our results concerning pro-dihedral extensions are not entirely new. We refer to Remark~\ref{rmk:jaulent} for a comparison with Jaulent's work~\cite{Jau81}, to Remark~\ref{rmk:Gillard_Carrol_Kisilevsky} for a comparison with the works~\cite{Gil76} by Gillard and~\cite{CarKis81} by Carroll--Kisilevsky, as well as to Remark~\ref{rmk:Finis_Hida_Tilouine} for some results concerning the anticyclotomic $\iwmu{\iw}$ invariant.

This paper has a long history, which has been summarized in~\cite[Acknowledgment]{CapNuc20}. As mentioned in~\loccit, the original \texttt{arXiv} preprint has been split up in two articles, the first being~\cite{CapNuc20}. This work is the second one, focusing on fake $\ZZ_p$-extensions. We are grateful to the anonymous referee for a thorough reading of our manuscript and for several remarks that improved the clarity and the readability of the text.
	
\section{Algebraic preliminaries}\label{sec:algprel}
\subsection{Group cohomology: notation and generalities}
Given a group $G$ and a $G$-module $B$, we denote by~$B^G$ (\rsp~$B_G$) the maximal submodule (\rsp~the maximal quotient) of $B$ on which $G$ acts trivially.
Moreover, let $\algnm{G}=\sum_{g\in G}g\in\ZZ[G]$ be the norm, $B[N_G]$ the kernel of multiplication by $N_G$ and $I_{G}$ the augmentation ideal of $\ZZ[G]$ defined as
\[
I_G=\kernel{(\ZZ[G]\longrightarrow \ZZ)}\qquad\text{ where the morphism is induced by }g\mapsto 1\text{ for all }g\in G.
\]
For every $i\in \ZZ$, let $\HH{i}{G}{B}$ denote the $i$th Tate cohomology group of $G$ with values in $B$. Similarly, for $i\geq 0$, let $\HHN{i}{G}{B}$ be the $i$th cohomology group of $G$ with values in $B$. For standard properties of these groups (which will be implicitly used without specific mention) the reader is referred to \cite{NeuSchWin08}. If $B'$ is another $G$-module and $f\colon B\to B'$ is a homomorphism of abelian groups, we say that $f$ is $G$-equivariant (\rsp~$G$-antiequivariant) if $f(g b)=g f(b)$ (\rsp~$f(g b)=-g f(b)$) for every $g\in G$ and $b\in B$. 

\subsection{Normic systems and their Tate cohomology} \label{subsec:normic_sys} Fix a prime number $p$. Arithmetic objects attached to layers of $\ZZ_p$-extensions have a structure of normic system, in the sense of Vauclair (see~\cite{Vau09}). After recalling the definition of a normic system, we will study its Tate cohomology, which is our primary goal in this paper. 

We start with the definition of a normic system, in a slightly different context than the one of ~\cite[Définition~2.1]{Vau09}. Let $\Gamma$ be a profinite group isomorphic to $\ZZ_p$, endowed with a decreasing filtration by closed subgroups $\Gamma_0=\Gamma\supseteq \Gamma_1\cdots\supseteq \Gamma_n\cdots$ for which there exists $g=g(\Gamma)\in\NN$, such that for all $m\geq n\geq g$ the equality $\gindex{\Gamma_n}{\Gamma_m}=p^{m-n}$ holds. Assume also that we are given an exact sequence of topological groups
\begin{equation}\label{eq:seq_gal_groups}
1 \longrightarrow \Gamma \longrightarrow \DD \longrightarrow  \splitG\longrightarrow 1.
\end{equation}
Set $G_{n}=\Gamma/\Gamma_n$ and $G_{m,n}=\Gamma_n/\Gamma_m$ for all $m\geq n$. Since the subgroups $\Gamma_n$ are closed, they are characteristic in~$\Gamma$ and hence normal in $\DD$, so that the quotient groups $\splitG[n]=\DD/\Gamma_n$ are defined for all $n\geq 0$; in particular,~$G_n$ can be regarded as a subgroup of $\splitG[n]$.

\begin{remark} 
The most common choice for the filtration $\{\Gamma_n\}$ in the above setting is given by $\Gamma_n=\Gamma^{p^n}$ for all~$n\geq 0$, and the reader can have this in mind in most of what follows. In that case, $g=0$ and $\gindex{\Gamma_n}{\Gamma_m}=p^{m-n}$ for all $m\geq n\geq 0$. The reason for the slight generality considered above comes from the arithmetic setting considered in Section \ref{subsection:arithmetic_set-up}, where $\Gamma$ is the local decomposition group inside a global Galois group. When the corresponding prime ideal splits, a shift in the numbering occurs and a more general filtration than $\{\Gamma^{p^n}\}$ needs to be considered.
\end{remark}

\begin{definition}\label{def:normic_system}
Let $(\Gamma, \{\Gamma_n\}_{n\in\mathbb{N}}, \DD)$ be as above. A $(\Gamma, \{\Gamma_n\}_{n\in\mathbb{N}}, \DD)$-normic system $\calB=(B_n,\nsu_{n,m},\nsd_{m,n})_{n,m\in\mathbb{N}}$ (or simply a $\DD$-normic system or a normic system if the groups are understood) is a collection of $\ZZ_p[\DD]$-modules~$B_n$ together with homomorphisms of $\ZZ_p[\DD]$-modules
\[
\xymatrix@C=2.5em{
	B_m\ar@<.5ex>[0,2]^{\nsd_{m,n}}&&B_n\ar@<.5ex>[0,-2]^{\nsu_{n,m}}&&\text{for all }m\geq n\geq 0}
\]
satisfying
\begin{itemize}
\item $(B_n,\nsu_{n,m})$ (\rsp~$(B_n,\nsd_{m,n})$) is a direct (\rsp~inverse) system of $\ZZ_p[\DD]$-modules: in particular, for all~$n\geq 0$, the compatibilities $\nsu_{n+1,n+2}\circ\nsu_{n,n+1}=\nsu_{n,n+2}$ and $\nsd_{n+1,n}\circ\nsd_{n+2,n+1}=\nsd_{n+2,n}$ hold; 
\item $B_n$ is fixed by $\Gamma_n$ (in particular, it can be regarded as a $\ZZ_p[\splitG[n]]$-module);
\item for all $m\geq n\geq 0$, $\nsu_{n,m}\circ \nsd_{m,n} = \algnm{G_{m,n}}$ and $\nsd_{m,n}\circ \nsu_{n,m} = \gindex{\Gamma_n}{\Gamma_m}$.
\end{itemize}
If $\calB=(B_n,\nsu_{n,m},\nsd_{m,n})$ and $\calB'=(B'_n,\nsu'_{n,m},\nsd'_{m,n})$ are two normic systems, a collection $f=(f_n)$ of maps $f_n\colon B_n \to B_n'$ is a morphism of normic systems if it is both a morphism of direct systems $(B_n,\nsu_{n,m})\to (B_n',\nsu'_{n,m})$ and of inverse systems $(B_n,\nsd_{m,n})\to (B_n',\nsd'_{m,n})$. This defines the category $\normsys$ of normic systems.
\end{definition}

The following lemma shows that $\normsys$ is abelian, using the concept of category of diagrams. For a general reference about categories of diagrams see~\cite[\S1.6~and~\S1.7]{Gro57}: we adopt notation and definitions from \ibid.

\begin{lemma}\label{lemma:normsys_abelian} 
The category $\normsys$ is a category of diagrams with commuting relations valued in the abelian category of $\ZZ_p[\DD]$-modules, so it is abelian.
\end{lemma}
\begin{proof} 
Consider the diagram scheme $S= (\NN,\Psi_\mathrm{d}\times\Psi_\mathrm{i},d)$ where 
\[
\Psi_\mathrm{d}=\{(n,m)\in\NN^2\text{ such that }n\leq m\}\qquad\text{ and }\qquad\Psi_\mathrm{i}=\{(m,n)\in\NN^2\text{ such that }m\geq n\}
\]
and the map $d$ sends every $\vec{v}_{n,m}=(n,m)\in\Psi_\mathrm{d}$ to the pair $(n,m)$ and every $\cev{v}_{m,n}=(m,n)\in\Psi_\mathrm{i}$ to the pair $(m,n)$. Fix a topological generator $\topgen\in\Gamma$ and for every $(i,j)\in\NN^2$ consider the set $R_{i,j}$ of $\ZZ_p[\DD]$-relations
\begin{equation}\label{eq:def_R}
R_{i,j}=\begin{cases}

\Bigl\{e_i,\topgen^{\gindex{\Gamma}{\Gamma_i}}=e_i,\cev{v}_{k,i}\vec{v}_{i,k}=\gindex{\Gamma_i}{\Gamma_k}e_i,\displaystyle{\vec{v}_{\ell,i}\cev{v}_{i,\ell}=\sum_{a=1}^{\gindex{\Gamma_\ell}{\Gamma_i}}\topgen^{ap^\ell}}\Bigr\}_{k\geq i,\ell\leq i}&\text{ if }i=j\\
\Bigl\{\vec{v}_{k,j}\vec{v}_{i,k}=\vec{v}_{i,j}\Bigr\}_{i\leq k\leq j}&\text{ if }j>i\\
\Bigl\{\cev{v}_{k,j}\cev{v}_{i,k}=\cev{v}_{i,j}\Bigr\}_{i\geq k\geq j}&\text{ if }i>j\\
\end{cases}
\end{equation}
where the $e_i$ are auxiliary elements corresponding to the identity, as \ibid. Set $\Sigma=(S,R)$: then $\normsys$ is the category of commutative diagrams $D\colon\Sigma\to \mathbf{Mod}_{\ZZ_p[\DD]}$ (where $\mathbf{Mod}_{\ZZ_p[\DD]}$ is the category of $\ZZ_p[\DD]$-modules) and is therefore abelian thanks to~\cite[Proposition~1.6.1]{Gro57}.
\end{proof}

Given a normic system $\calB=(B_n,\nsu_{n,m},\nsd_{m,n})$ one can consider the Tate cohomology groups $\HH{i}{G_n}{B_n}$ for $i\in\ZZ$ as well as the usual cohomology groups $\HHN{i}{G_n}{B_n}$ for $i\geq 0$. These groups acquire the structure of $\ZZ_p[\splitG]$-modules via the conjugation action of $\splitG$.

\begin{definition}\label{def:conditions_Inj_Gal}
We say that the normic system $\calB=(B_n,\nsu_{n,m},\nsd_{m,n})$ satisfies condition~\ref{cond:inj} if 
\begin{equation}\label{cond:inj}\tag*{\textbf{\textup{(Inj)}}}
\begin{minipage}{.8\textwidth}
\baselineskip=1.25\baselineskip
$\nsu_{n,m}$ is injective for all $m\geq n\geq 0$
\end{minipage}
\end{equation}
and that it satisfies condition~\ref{cond:surj} if
\begin{equation}\label{cond:surj}\tag*{\textbf{\textup{(Gal)}}}
\begin{minipage}{.8\textwidth}
\baselineskip=1.25\baselineskip
$\nsu_{n,m}(B_n)=B_m^{G_{m,n}}$ for all $m\geq n\geq 0$.
\end{minipage}
\end{equation}
\end{definition}

Let $\calB=(B_n,\nsu_{n,m},\nsd_{m,n})$ be a normic system. Our next task is to construct, under suitable hypotheses, functorial morphisms of $\ZZ_p[\splitG]$-modules
\[
\xymatrix@C=3em{
	\HH{i}{G_m}{B_m}\ar@<.35em>[0,2]^{\HHm{i}{\nsd}=\HHm{i}{\nsd_{m,n}}}&&\HH{i}{G_n}{B_n}\ar@<.35em>[0,-2]^{\HHm{i}{\nsu}=\HHm{i}{\nsu_{n,m}}}&\text{for all }m\geq n\geq 0.}
\]
The ``ascending'' maps $\HHm{i}{\nsu}$ are well-defined independently of any assumption on the normic system, and we show in \S\ref{subsec:comparison_inf} that they are related to inflation in positive degrees. The definitions of the ``descending'' maps $\HHm{i}{\nsd}$ require condition~\ref{cond:inj} in odd degrees and both conditions~\ref{cond:inj} and \ref{cond:surj} in even degrees.
\subsubsection{Cohomology maps in degree $-1$} \label{subsubsec:minus1}
We first define $\HHm{-1}{\nsu}$ by 
\[
\HHm{-1}{\nsu}(y\mod{I_{G_n}}B_n)=\nsu_{n,m}(y)\mod{I_{G_m}B_m} \qquad \text{for all $y\in\HH{-1}{G_n}{B_n}$.}
\] 
Observe that $\HHm{-1}{\nsu}$ is well-defined, because
\begin{align*}
\nsu_{n,m}(B_n[\algnm{G_n}])\subseteq B_m[\algnm{G_m}]\qquad\text{ and }\qquad \nsu_{n,m}(I_{G_n}B_n)\subseteq I_{G_m}B_m,
\end{align*}
since $\nsu_{n,m}$ is $G_m$-equivariant. 

Assume now that $\calB$ satisfies~\ref{cond:inj}. Then we define $\HHm{-1}{\nsd}$ by 
\[
\HHm{-1}{\nsd}(x\mod{I_{G_m}B_m})=\nsd_{m,n}(x)\mod{I_{G_n}B_n} \qquad \text{for all $x\in\HH{-1}{G_m}{B_m}$.}
\] 
To see that $\HHm{-1}{\nsd}$ is also well-defined,
observe that
\begin{align*}
\nsu_{n,m}\circ\algnm{G_n}\circ \nsd_{m,n}&=\algnm{G_n}\circ \nsu_{n,m}\circ \nsd_{m,n}\\
&=\algnm{G_n}\circ\algnm{G_{m,n}}\\
&=\algnm{G_m}.
\end{align*}
Since $\nsu_{n,m}$ is assumed to be injective, this shows that $\nsd_{m,n}(B_m[\algnm{G_m}])\subseteq B_n[\algnm{G_n}]$. Since $\nsd_{m,n}$ is $G_m$-equivariant, it is immediate that $\nsd_{m,n}(I_{G_m}B_m)\subseteq I_{G_n}B_n$, so $\HHm{-1}{\nsd}$ is well-defined. 

Observe that we have
\[
\HHm{-1}{\nsd_{m,n}}\circ \HHm{-1}{\nsu_{n,m}} = \gindex{\Gamma_n}{\Gamma_m} \qquad \text{ and }\qquad\HHm{-1}{j_{n,m}}\circ \HHm{-1}{\nsd_{m,n}} = \gindex{\Gamma_n}{\Gamma_m}.
\]
Indeed, using the relation between $\nsd_{m,n}$ and $\nsu_{n,m}$ in the definition of a normic system, we get 
\begin{align*}
\HHm{-1}{\nsd_{m,n}}\circ \HHm{-1}{\nsu_{n,m}}(x \mod{I_{G_n}B_n}) & = \nsd_{m,n}\circ\nsu_{n,m}(x)\mod{I_{G_n}B_n}\\
&=\gindex{\Gamma_n}{\Gamma_m}x\mod{I_{G_n}B_n}\\
\intertext{and}
\HHm{-1}{\nsu_{n,m}}\circ \HHm{-1}{\nsd_{m,n}}(x \mod{I_{G_n}B_n)}& = \nsu_{n,m}\circ\nsd_{m,n}(x)\mod{I_{G_m}B_m}\\&=\algnm{G_{m,n}}x\mod{I_{G_m}B_m}\\& = \gindex{\Gamma_n}{\Gamma_m}x \mod{I_{G_m}B_m}
\end{align*}
where the last equality follows from the fact that the action of $G_{m,n}$ is trivial on $\HH{-1}{G_m}{B_m}$. Moreover, both $\HHm{-1}{\nsd_{m,n}}$ and $\HHm{-1}{\nsu_{n,m}}$ are morphisms of $\ZZ_p[\splitG]$-modules, because both $\nsd_{m,n}$ and $\nsu_{n,m}$ are $\splitG[m]$-equivariant.

\subsubsection{Cohomology maps in degree $0$}\label{subsubsub:HH0}
We first define $\HHm{0}{\nsu}$ by 
\[
\HHm{0}{\nsu}(y\mod{\algnm{G_n}B_n})=\gindex{\Gamma_n}{\Gamma_m}\nsu_{n,m}(y)\mod{\algnm{G_m}B_m} \qquad \text{for all $y\in\HH{0}{G_n}{B_n}$.}
\] 
Observe that, since $\nsu_{n,m}$ takes values in $B_m^{G_{m,n}}$, there are inclusions $\gindex{\Gamma_n}{\Gamma_m}\nsu_{n,m}(B_n)^{G_n} \subseteq \gindex{\Gamma_n}{\Gamma_m}B_m^{G_m}\subseteq B_m^{G_m}$ and 
\[
\gindex{\Gamma_n}{\Gamma_m}\nsu_{n,m}\circ \algnm{G_n} = \algnm{G_n}\circ \nsu_{n,m} \circ \gindex{\Gamma_n}{\Gamma_m} =  \algnm{G_n}\circ \nsu_{n,m}\circ \nsd_{m,n} \circ \nsu_{n,m} =  \algnm{G_m}\circ \nsu_{n,m}
\]
so that $\gindex{\Gamma_n}{\Gamma_m}\nsu_{n,m}(\algnm{G_n}B_n) \subseteq \algnm{G_m}B_m$, showing that $\HHm{0}{\nsu}$ is well-defined. 

Assume now that $\calB$ satisfies both~\ref{cond:inj} and~\ref{cond:surj}. Then we define $\HHm{0}{k}$ by
\[
\HHm{0}{\nsd}(x\mod{\algnm{G_m}B_m})=(\nsu_{n,m})^{-1}(x)\mod{\algnm{G_n}B_n} \qquad \text{for all $x\in\HH{0}{G_m}{B_m}$.}
\]  
To check that $\HHm{0}{\nsd}$ is well-defined, observe that~\ref{cond:inj} and~\ref{cond:surj} imply that $\nsu_{n,m}$ defines a $G_n$-isomorphism of~$B_n$ onto $B_m^{G_{m,n}}$, which restricts to an isomorphism of $B_n^{G_n}$ onto $B_m^{G_m}$. Therefore, given $x \in B_m^{G_m}$, the element $(\nsu_{n,m})^{-1}(x)\in B_n^{G_n}$ is uniquely defined. Moreover, if $x= \algnm{G_m}z$ for some $z\in B_m$, then $w=\algnm{G_{m,n}}z$ belongs to $B_m^{G_{m,n}}$, so that 
\[
x = \algnm{G_m}z = \algnm{G_n}(w)=  \algnm{G_n}\circ \nsu_{n,m}\bigl((\nsu_{n,m})^{-1}(w)\bigr) = \nsu_{n,m} \circ \algnm{G_n}\bigl((\nsu_{n,m})^{-1}(w)\bigr),
\]  
\ie $(\nsu_{n,m})^{-1}(x)\in \algnm{G_n}B_n$. This concludes the proof that $\HHm{0}{\nsd}$ is well-defined. 

As in the case $i=-1$, we have 
\[
\HHm{0}{\nsd_{m,n}}\circ \HHm{0}{\nsu_{n,m}} = \gindex{\Gamma_n}{\Gamma_m} \qquad\text{and}\qquad \HHm{0}{\nsu_{n,m}}\circ \HHm{0}{\nsd_{m,n}} = \gindex{\Gamma_n}{\Gamma_m}
\]
and both $\HHm{0}{\nsd_{m,n}}$ and $\HHm{0}{\nsu_{n,m}}$ are morphisms of $\ZZ_p[\splitG]$-modules, because $\nsu_{n,m}$ is $\splitG[m]$-equivariant.

\subsubsection{Cohomology maps in arbitrary degree}\label{subsubsec:arbitrarydegree}
Fix a topological generator $\topgen\in\Gamma$. For each $n\in \mathbb{N}$, $\topgen$ maps to a generator $g_n\in G_n$ and we let $\chi_n\in \Hom(G_n,\QQ/\ZZ)=\HH{1}{G_n}{\QQ/\ZZ}$ be the map sending $g_n$ to $\gindex{\Gamma}{\Gamma_n}^{-1}$. We claim that
\[
\infl(\chi_n)=\gindex{\Gamma_n}{\Gamma_m}\chi_m\qquad\text{ for }m\geq n
\]
where $\infl\colon \HH{1}{G_n}{\QQ/\ZZ}\to \HH{1}{G_m}{\QQ/\ZZ}$ is the inflation map. Indeed, since $g_m$ maps to $g_n$, $\infl(\chi_n)\in\Hom(G_m,\QQ/\ZZ)$ sends $g_m$ to $\chi_n(g_n)=\gindex{\Gamma}{\Gamma_n}^{-1}$, which in turn is equal to the ratio $\gindex{\Gamma_n}{\Gamma_m}\gindex{\Gamma}{\Gamma_m}^{-1}=\gindex{\Gamma_n}{\Gamma_m}\chi_m(g_m)$. Since the connecting homomorphism $\cch{1}\colon \HH{1}{G_n}{\QQ/\ZZ}\to\HH{2}{G_n}{\ZZ}$ corresponding to the exact sequence 
\[
0\longrightarrow \ZZ \longrightarrow \QQ \longrightarrow \QQ/\ZZ \longrightarrow 0
\]
induced by the inclusion $\ZZ\subset \QQ$ is an isomorphism, there exists a unique element $\varkappa_n\in\HH{2}{G_n}{\ZZ}$ such that 
\[
\cch{1}(\chi_n)=\varkappa_n.
\]

We now go back to the cohomology of the normic system $\calB$. Fix $i\in\ZZ$ and, for all $h\geq 0$, let 
\[
\tateiso{n}{h}\colon\HH{i}{G_n}{B_n} \longrightarrow \HH{i+2h}{G_n}{B_n}
\]
be the isomorphism obtained by taking the cup product with $\varkappa_n$ $h$ times. Extend the definition to $h<0$ by setting
\[
\tateiso{n}{h}=\text{ inverse of }\tateiso{n}{-h}.
\] 
If $i$ is odd (\rsp~even), write $i=2h-1$ and set $\varepsilon=-1$ (\rsp~$i=2h$ and $\varepsilon=0$).
Then we define, for $m\geq n\geq 0$,
\[
\HHm{i}{\nsu}=\HHm{i}{\nsu_{n,m}} = \tateiso{m}{h} \circ \HHm{\varepsilon}{\nsu_{n,m}} \circ \tateiso{n}{-h}\colon \HH{i}{G_n}{B_n}\longrightarrow \HH{i}{G_m}{B_m}.
\]
Similarly, if $i$ is odd (\rsp~even) and $\calB$ satisfies~\ref{cond:inj} (\rsp~both~\ref{cond:inj} and~\ref{cond:surj}), we define 
\[
\HHm{i}{\nsd}=\HHm{i}{\nsd_{m,n}} = \tateiso{n}{h} \circ \HHm{\varepsilon}{\nsd_{m,n}} \circ \tateiso{m}{-h}\colon \HH{i}{G_m}{B_m}\longrightarrow \HH{i}{G_n}{B_n} .
\]

\begin{remark} 
The definition of these homomorphisms is independent of the choice of $\topgen$. To see this, let $\topgen'\in\Gamma$ be another topological generator, that must be of the form $\topgen^u$ for some $u\in\ZZ_p^\times$. Let $g'_n=g_n^u\in G_n$ be the image of $\topgen'$ and let $\chi'_n\in\HH{1}{G_n}{\QQ/\ZZ}$ be the map sending $g_n'$ to $\gindex{\Gamma}{\Gamma_n}^{-1}$: it corresponds to a unique class $\varkappa_n'\in \HH{2}{G_n}{\ZZ}$. We have $\chi_n=u\chi'_n$ because
\[
\chi_n(g_n')=\chi_n(g_n^u)=u\chi_n(g_n)=u\cdot\gindex{\Gamma}{\Gamma_n}^{-1}=u\chi'_n(g'_n)
\] 
and in particular
\[
\tateiso{n}{h} = u^{h}\tateiso[\prime]{n}{h},
\]
showing that the definitions of $\HHm{i}{\nsu}$ and of $\HHm{i}{\nsd}$ are independent of the choice of $\topgen$.
\end{remark}

Arguing by induction, one can show that $\HHm{i}{\nsd_{m,n}}$ and $\HHm{i}{\nsu_{n,m}}$ are $\ZZ_p[\splitG]$-homomorphisms: we already observed this in the cases $i=-1, 0$. To prove the assertion for $i\in\ZZ$, let $\kappa\colon\DD\to \ZZ_p^\times$ be the homomorphism defined by 
\[
g\topgen g^{-1} = \topgen^{\kappa(g)} \qquad \text{for $g\in\DD$}. 
\]
For all $g\in\DD$, denote by $g^*$ the corresponding conjugation automorphism in cohomology. Then (see \cite[Proposition~I.1.5.3]{NeuSchWin08})
\[
g^* \circ \tateiso{n}{} = \left(\cup (g^* \circ \varkappa_n)\right)\circ g^* = \kappa(g)^{-1}\tateiso{n}{}\circ g^*
\]
since
\[
g^* \circ \varkappa_n = g^*\circ \cch{1}(\chi_n)= \cch{1}(g^*\circ \chi_n) = \cch{1}(g^*\circ \chi_n) 
\]
and $g^*\circ \chi_n: \gamma_0 \mapsto g^{-1}\gamma_0g$. 
Therefore, on the one hand  we can write 
\begin{align*}
\HHm{i+2}{\nsu_{n,m}}\circ g^*&= \HHm{i+2}{\nsu_{n,m}}\circ g^*\circ\tateiso{n}{}\circ\tateiso{n}{-1} \\
&= \kappa(g)^{-1}\HHm{i+2}{\nsu_{n,m}}\circ \tateiso{n}{}\circ g^* \circ\tateiso{n}{-1}\\
&=\kappa(g)^{-1}\tateiso{m}{}\circ\HHm{i}{\nsu_{n,m}}\circ  g^* \circ\tateiso{n}{-1}.
\end{align*}
On the other hand we have
\begin{align*}
g^*\circ\HHm{i+2}{\nsu_{n,m}}&=g^*\circ\HHm{i+2}{\nsu_{n,m}}\circ\tateiso{n}{}\circ\tateiso{n}{-1}\\
&=g^*\circ\tateiso{m}{}\circ\HHm{i}{\nsu_{n,m}}\circ\tateiso{n}{-1}\\
&=\kappa(g)^{-1}\tateiso{m}{}\circ g^* \circ\HHm{i}{\nsu_{n,m}}\circ\tateiso{n}{-1}.
\end{align*}
We deduce that 
\[
g^*\circ\HHm{i}{\nsu_{n,m}} = \HHm{i}{\nsu_{n,m}}\circ g^* \Longleftrightarrow g^*\circ\HHm{i+2}{\nsu_{n,m}} = \HHm{i+2}{\nsu_{n,m}}\circ g^*
\]
showing the inductive step. A similar proof works when $\nsu$ is replaced by $\nsd$.

\subsubsection{Comparison with other maps}\label{subsec:comparison_inf} 
When $i\geq 1$, other maps $\HH{i}{G_n}{B_n}\to \HH{i}{G_m}{B_m}$ exist for $m\geq n\geq 0$. These are obtained by composing the inflation map
\[
\infl= \infl^i\colon \HH{i}{G_n}{B_m^{G_{m,n}}}\longrightarrow \HH{i}{G_m}{B_m}
\]
with the map induced by $\nsu_{n,m}\colon B_n\to B_m^{G_m}$ in cohomology
\[
j^* = j_{n,m}^{*,i}\colon \HH{i}{G_n}{B_n}\longrightarrow \HH{i}{G_n}{B_m^{G_{m,n}}},
\]
and they are denoted $\infl^i\circ \nsu^{*,i}_{n,m}$ or simply $\infl \circ \nsu^*$. In what follows we analyse the relationship of these maps with the ones defined in \S\ref{subsubsec:arbitrarydegree}.

\begin{remark}\label{rmk:cp}
For every $[b]\in B_n[\algnm{G_n}]/I_{G_n}B_n=\HH{-1}{G_n}{B_n}$, there is an explicit description of the cup product $[b]\cup \varkappa_n\in \HH{1}{G_n}{B_n}$ as the class $\xi^{(b)}_n$ of the $1$-cocycle sending the generator $g_n\in G_n$ to $b$. This is probably well-known, but we provide a proof because we were unable to find an explicit reference in the literature. By~\cite[Proposition~I.1.4.8]{NeuSchWin08}, the cup product $[b]\cup \varkappa_n$ is the class of the $2$-cocycle
\begin{equation}\label{cases:eta} 
g_n\longmapsto \sum_{\sigma\in G_n}\sigma(b)\otimes\sigma\widetilde{\varkappa}_n(\sigma^{-1},g_n)=\sum_{\sigma\in G_n}\sigma(b)\otimes\widetilde{\varkappa}_n(\sigma^{-1},g_n)
\end{equation}
where $\widetilde{\varkappa}_n$ is any $2$-cochain representing $\varkappa_n$. To compute this expression explicitly, recall that $\varkappa_n=\cch{1}(\chi_n)$ where $\chi_n\colon G_n\to\QQ/\ZZ$ is the character sending $g_n$ to $\gindex{\Gamma}{\Gamma_n}^{-1}$. The definition of the connecting homomorphism $\cch{1}\colon \HH{1}{G_n}{\QQ/\ZZ}\to\HH{2}{G_n}{\ZZ}$ is given explicitly in~\cite[proof of Theorem~I.1.3.2]{NeuSchWin08} in terms of homogeneous cocycles: it is the connecting homomorphism $\sch\colon \kernel\partial_{\QQ/\ZZ}\to\cokernel{\partial_{\ZZ}}$ obtained by applying the snake lemma to the diagram on~\cite[page~27]{NeuSchWin08}. By invoking~\cite[page~14]{NeuSchWin08}, we can replace homogeneous cochains by inhomogeneous ones and transform the diagram in \loccit into the commutative diagram of exact rows
\begin{equation*}
\begin{split}
\xymatrix@C=3em{
&\mathscr{C}^1(G_n,\ZZ)\ar@{->}[1,0]_{\partial^2_{\ZZ}}\ar@{->}[0,1]&\mathscr{C}^1(G_n,\QQ)\ar@{->}[1,0]_{\partial^2_{\QQ}}\ar@{->}[0,1]&\mathscr{C}^1(G_n,\QQ/\ZZ)\ar@{->}[1,0]_{\partial^2_{\QQ/\ZZ}}\ar@{->}[0,1]&0\\
0\ar@{->}[0,1]&\mathscr{C}^2(G_n,\ZZ)\ar@{->}[0,1]&\mathscr{C}^2(G_n,\QQ)\ar@{->}[0,1]&\mathscr{C}^2(G_n,\QQ/\ZZ)}
\end{split}\end{equation*}
In this notation, $\chi_n\in\kernel{\partial^2_{\QQ/\ZZ}}$ and, by definition of the connecting homomorphism in the snake lemma, $\cch{1}(\chi_n)$ is the class of $\widetilde{\varkappa}_n=\partial^2_{\QQ}(\chi_n')\in \cokernel{\partial^2_{\ZZ}}$, where $\tilde{\chi}_n\colon G_n\to\QQ$ is any $1$-cochain lifting $\chi_n$. As lift, choose the cochain $\tilde{\chi}_n\colon g_n^a\mapsto a\gindex{\Gamma}{\Gamma_n}^{-1}$ for $1\leq a\leq \gindex{\Gamma}{\Gamma_n}$. By the explicit description of $\partial^1_{\QQ}$ on inhomogeneous $1$-cochains given in~\cite[page~14]{NeuSchWin08}, we obtain
\[
\widetilde{\varkappa}_n(g_n^a,g_n)=\tilde{\chi}_n(g_n)-\tilde{\chi}_n(g_n^{a+1})+\tilde{\chi}_n(g_n^a)=\begin{cases}
1&\text{ if }a= \gindex{\Gamma}{\Gamma_n}\\
0&\text{ if }1\leq a\leq \gindex{\Gamma}{\Gamma_n}-1
\end{cases}
\]
It follows that that the sum in~\eqref{cases:eta} contains only the summand corresponding to $\sigma=g_n^{\gindex{\Gamma}{\Gamma_n}}=\id{G_n}$ and this term equals $b$, showing that $[b]\cup\varkappa_n=\xi_n^{(b)}$.
\end{remark}

\begin{lemma}\label{lemma:HHm(j)_and_inf}
Let $\calB$ be a normic system and let $i\geq 1$ be an integer. Then $\HHm{i}{\nsu}$ coincides with $\infl^{i} \circ \nsu^\ast$.
\end{lemma}
\begin{proof}
If $i$ is odd \textup{(}\rsp~even\textup{)}, write $i=2h-1$ and set $\varepsilon=-1$ \textup{(}\rsp~$i=2h$ and $\varepsilon=0$\textup{)}.
By definition of $\HHm{i}{\nsu}$ and since $\tateiso{m}{h}$ is an isomorphism, the statement of the lemma is equivalent to the commutativity of the following diagram, for all $m\geq n\geq 0$:
\begin{equation}\begin{split}\label{eq:nsuinfl}
	\xymatrix@C 4em@R 4em{
		\HH{\varepsilon}{G_m}{B_m}\ar@{->}[0,1]^{\tateiso{m}{h}} &\HH{i}{G_m}{B_m}\\
		\HH{\varepsilon}{G_n}{B_n}\ar@{->}[0,1]^{\tateiso{n}{h}}\ar@{->}[-1,0]^{\HHm{\varepsilon}{\nsu}} &\HH{i}{G_n}{B_n}\ar@{->}[-1,0]_{\infl^{i}\circ \nsu^\ast}
	}
\end{split}\end{equation}
We claim that commutativity of~\eqref{eq:nsuinfl} for all $i\geq 1$ follows once we can prove it commutes for $i=\varepsilon+2$. Indeed, for arbitrary $i\geq 3$ we have
\begin{align*}
\HHm{i}{\nsu}&=\tateiso{m}{h} \circ \HHm{\varepsilon}{\nsu_{n,m}} \circ \tateiso{n}{-h}\\
&=\tateiso{m}{h-1} \circ \tateiso{m}{} \circ \HHm{\varepsilon}{\nsu_{n,m}} \circ \tateiso{m}{-1} \circ \tateiso{n}{-h+1}\\
\text{(by assumption)}&=\tateiso{m}{h-1} \circ \infl^{\varepsilon+2}\circ \nsu^\ast \circ \tateiso{n}{-h+1}\\
\text{(compatibility of $\cup$ and $\infl^i\circ \nsu^\ast$)}&=\infl^i\circ \nsu^\ast.
\end{align*}
Hence, it now remains to show that \eqref{eq:nsuinfl} is commutative for $i=1,2$
 
Consider first the case $i=1$. Given the definition of $\HHm{-1}{\nsu}$ and the explicit description of the cup product described in Remark \ref{rmk:cp}, commutativity of~\eqref{eq:nsuinfl} is equivalent to
\begin{equation*}
\infl^{1}\circ\nsu^{\ast}(\xi^{(b)}_n)=\xi^{(j(b))}_{m}
\end{equation*}	
and this follows by the definition of the map $\infl^{1}\circ\nsu^{\ast}$ on cocycles (see \cite[\S~I.5]{NeuSchWin08}).

Suppose now that $i=2$. We first show that
\begin{equation*}
\xymatrix@C 4em @R 4em{
	\HHN{0}{G_m}{B_m}\ar@{->}[0,1]^{\tateiso{m}{}} &\HH{2}{G_m}{B_m}\\
	\HHN{0}{G_n}{B_n}\ar@{->}[0,1]^{\tateiso{n}{}} \ar@{->}[-1,0]^{\gindex{\Gamma_n}{\Gamma_m} \nsu_{n,m}}&\HH{2}{G_n}{B_n}\ar@{->}[-1,0]_{\infl^{2}\circ\nsu^{\ast}}
}
\end{equation*} 
commutes. Note that in the above diagram we have ordinary cohomology on the left-hand side (as opposed to Tate cohomology): this allows us to consider inflation in degree $0$, which is simply the identity: 
\[
\infl^{0}\colon \HHN{0}{G_n}{B_m^{G_{m,n}}}=B_m^{G_m}\longrightarrow \HHN{0}{G_m}{B_m}=B_m^{G_m}.
\]  
Let $x\in \HHN{0}{G_n}{B_n}$: the compatibility of the cup product with inflation yields
\begin{align*}
\infl^{2}\circ\nsu^{\ast}(x \cup \varkappa_n)
&=\infl^{2}(\nsu_{n,m}(x)\cup \varkappa_n)\\
&= \infl^{0}\bigl(\nsu_{n,m}(x)\bigr)\cup \infl^{2}(\varkappa_n)\\
&= \nsu_{n,m}(x) \cup \bigl(\gindex{\Gamma_n}{\Gamma_m}\varkappa_m\bigr)\\
&= \gindex{\Gamma_n}{\Gamma_m}\nsu_{n,m}(x) \cup \varkappa_m.
\end{align*}
This shows that the above diagram commutes. To show that \eqref{eq:nsuinfl} commutes for $i=2$ and conclude the proof, we just need to observe that the cup product on Tate cohomology is compatible with that on ordinary cohomology. This can be checked, for instance, using the definition on cochains of the two cup products (compare \mbox{\cite[($\ast$), \S~I.4]{NeuSchWin08}} and \cite[proof of Proposition~I.1.4.6]{NeuSchWin08}).
\end{proof}

\begin{remark}\label{rmk:h0surjnorm}
Assume that $\nsd_{m,n}$ is surjective. Then the map $\HHm{0}{\nsu_{n,m}}$ factors through ``multiplication by $\gindex{\Gamma_n}{\Gamma_m}$'', \ie there is a map $\hatnsu{n,m}\colon \HH{0}{G_n}{B_n} \to \HH{0}{G_m}{B_m}$ such that $\HHm{0}{\nsu_{n,m}}=\gindex{\Gamma_n}{\Gamma_m}\hatnsu{n,m}$. This is defined as
\[
\hatnsu{n,m}(y\mod{\algnm{G_n}}B_n)=\nsu_{n,m}(y)\mod{\algnm{G_m}B_m} \qquad \text{for all $y\in B_n^{G_n}$.}
\] 
To see that this is well defined, let $b_n\in B_n^{G_n}$. Then it is clear that $\nsu_{n,m}(b)\in B_m^{G_m}$. Suppose now that $b_n= \algnm{G_n}(b'_n)\in \algnm{G_n}B_n$. Since $\nsd_{m,n}$ is surjective, there exists $b'_m\in B_m$ such that $\nsd_{m,n}(b'_m)=b'_n$. Then
\begin{align*}
\nsu_{n,m}(b_n) &= \nsu_{n,m} \circ\algnm{G_n}(b'_n)\\ 
&= \nsu_{n,m} \circ\algnm{G_n}\circ \nsd_{m,n}(b'_m)\\
& =\nsu_{n,m}\circ\nsu_{0,n}\circ\nsd_{n,0}\circ \nsd_{m,n}(b'_m)\\
&= \nsu_{0,m}\circ\nsd_{m,0} (b'_m) = \algnm{G_m}(b_m)
\end{align*}
showing that $\hatnsu{n,m}$ is well-defined. By definition of $\HHm{0}{j}$, the relation $\HHm{0}{\nsu_{n,m}}=\gindex{\Gamma_n}{\Gamma_m}\hatnsu{n,m}$ holds. Using~$\hatnsu{n,m}$, Lemma~\ref{lemma:HHm(j)_and_inf} can be restated, under the assumption that $\nsd_{m,n}$ is surjective, as the commutativity of the following diagram :
\[
\xymatrix@C 4em @R 4em{
	\HH{0}{G_m}{B_m}\ar@{->}[0,1]^{\tateiso{m}{}}  &\HH{2}{G_m}{B_m}\\
	\HH{0}{G_n}{B_n}\ar@{->}[0,1]^{\tateiso{n}{}}\ar@{->}[-1,0]^{\gindex{\Gamma_n}{\Gamma_m}\hatnsu{n,m}}&\HH{2}{G_n}{B_n}\ar@{->}[-1,0]_{\HHm{2}{\nsu}}
}
\] 
\end{remark}

\begin{remark}
One can wonder whether results similar to that of Lemma~\ref{lemma:HHm(j)_and_inf} hold when $j$ is replaced by $k$. Assume that $\calB$ satisfies~\ref{cond:inj} and~\ref{cond:surj}: under this assumption, the inverse maps $(\nsu_{n,m})^{-1}\colon B_m^{G_{m,n}}\to B_n$ and $\HHm{i}{\nsd_{m,n}}$ are defined and we have maps 
\[
(\nsu^*)^{-1}=(\nsu^{*,i}_{m,n})^{-1}\colon\HH{i}{G_m}{B_m^{G_{m,n}}} \longrightarrow \HH{i}{G_m}{B_n}
\]
induced by $(\nsu_{m,n})^{-1}$. Write $i=2h$ or $i=2h-1$, according to whether $i$ is even or odd. Then we claim that 
\begin{equation}\label{eq:kdefl}
(\nsu^{*,i}_{m,n})^{-1}\circ \defl^i =\gindex{\Gamma_n}{\Gamma_m}^{-h}\HHm{i}{\nsd_{m,n}} \qquad\text{ for all }i\leq 0
\end{equation}
where $\defl$ denotes the deflation map
\[
\defl=\defl^i\colon \HH{i}{G_m}{B_m} \longrightarrow \HH{i}{G_n}{B_m^{G_{m,n}}}
\]
as defined by Weiss (see~\cite{Wei59}). Equality~\eqref{eq:kdefl} is clear for $i=0,-1$ since in these cases $h=0$ and we have
\begin{align*}
\defl^0(x \mod{\algnm{G_{m}}B_m}) &= x \mod{\algnm{G_n}B_m^{G_{m,n}}} \quad \text{ for }x\in B_m^{G_m}
\intertext{as well as}
\defl^{-1}(x \mod{I_{G_{m}}B_m}) &= \algnm{G_{m,n}}x \mod{I_{G_n}B_m^{G_{m,n}}} \quad \text{ for }x\in B_m[\algnm{G_m}]
\end{align*}
(see \cite[(1) and Proposition 1]{Wei59}). For arbitrary $i$, we argue by backward induction: assume that~\eqref{eq:kdefl} holds for a given $i\leq 0$. Then for $x\in \HH{i-2}{G_m}{B_m}$
\begin{align*}
(\nsu^{*,i-2}_{n,m})^{-1}\circ \defl^{i-2}(x) \cup \varkappa_n &=(\nsu_{n,m}^{*,i})^{-1} (\defl^{i-2}(x) \cup \varkappa_n)\\ 
& = (\nsu_{n,m}^{*,i})^{-1}\circ \defl^{i}(x\cup \infl^2(\varkappa_n)) \quad \text{(by \cite[Theorem 3]{Wei59})}\\
&=  (\nsu_{n,m}^{*,i})^{-1}\circ \defl^{i}(x\cup \gindex{\Gamma_n}{\Gamma_m}\varkappa_m)\\ 
&= \gindex{\Gamma_n}{\Gamma_m}^{-(h-1)}\HHm{i}{\nsd_{m,n}}(x\cup \varkappa_m) \\
&= \gindex{\Gamma_n}{\Gamma_m}^{-(h-1)}\HHm{i-2}{\nsd_{m,n}}(x)\cup \varkappa_n.
\end{align*}
Since the cup product with $\varkappa_n$ is an isomorphism,~\eqref{eq:kdefl} holds for $i-2$.
\end{remark}

\subsubsection{Functorial behaviour}
We now come to the main results of this section, describing the functorial behaviour of the maps defined in \S\ref{subsubsec:arbitrarydegree}. First of all, observe that if $f=(f_n)\colon \calB \to \calB'$ is a morphism of normic systems then, for every $i\in\ZZ$, there are maps
\[
f^\ast=f_n^{\ast,i}\colon \HH{i}{G_n}{B_n} \longrightarrow \HH{i}{G_n}{B'_n}.
\]
The diagrams
\begin{equation}\label{diag:f*_functorial}
\begin{array}{lcr}
\xymatrix@C=4em @R3em{
	\HH{i}{G_{m}}{B_{m}}\ar@{->}[0,1]^{f^\ast}&\HH{i}{G_{m}}{B_{m}'}\\
	\HH{i}{G_{n}}{B_{n}}\ar@{->}[0,1]^{f^\ast}\ar@{->}[-1,0]^{\HHm{i}{\nsu}}&\HH{i}{G_{n}}{B_{n}'}\ar@{->}[-1,0]^{\HHm{i}{\nsu}}
}
&\parbox[t][4.7em][c]{6em}{\begin{center}\text{and}\end{center}}&
\xymatrix@C=4em @R3em{
	\HH{i}{G_m}{B_m}\ar@{->}[1,0]^{\HHm{i}{\nsd}}\ar@{->}[0,1]^{f^\ast}&\HH{i}{G_m}{B'_m}\ar@{->}[1,0]^{\HHm{i}{\nsd}}\\
	\HH{i}{G_n}{B_n}\ar@{->}[0,1]^{f^\ast}&\HH{i}{G_n}{B'_n}
}
\end{array}
\end{equation}
are commutative (in the rightmost one we are assuming that $\HHm{i}{\nsd}$ is defined): to show this, it is enough to consider the cases $i=-1,0$ because the isomorphism $\tateiso{n}{}$ commutes with $f_n^\ast$ . Since $f_n$ commutes with  $\nsu_{n,m}$ and $\nsd_{m,n}$, the cases $i=-1,0$ are easily proved by direct inspection. The behaviour of the maps $\HHm{i}{j}$ and $\HHm{i}{k}$ with respect to connecting homomorphisms is given by the following two propositions. 

\begin{proposition}\label{prop:long_cohom_up}
Let $0\to\calB'\to\calB\to\calB''\to 0$ be an exact sequence of normic systems. Then, for all $m\geq n\geq 0$, the following diagram with exact rows commutes
\begin{equation*}
\xymatrix@C=2.25em{
0\ar@{->}[0,1]&\HHN{0}{G_m}{B'_m}\ar@{->}[0,1]&\HHN{0}{G_m}{B_m}\ar@{->}[0,1]&\HHN{0}{G_m}{B''_m}\ar@{->}[0,1]^{\cch{0}}&\HH{1}{G_m}{B'_m}\ar@{->}[0,1]&\cdots\\
0\ar@{->}[0,1]&\HHN{0}{G_n}{B'_n}\ar@{->}[-1,0]^{\nsu'}\ar@{->}[0,1]&\HHN{0}{G_n}{B_n}\ar@{->}[-1,0]^{\nsu}\ar@{->}[0,1]&\HHN{0}{G_n}{B''_n}\ar@{->}[-1,0]^{\nsu''}\ar@{->}[0,1]^{\cch{0}}&\HH{1}{G_n}{B'_n}\ar@{->}[-1,0]^{\HHm{1}{\nsu'}}\ar@{->}[0,1]&\cdots\\
&	\cdots\ar@{->}[0,1]&\HH{1}{G_m}{B_m}\ar@{->}[0,1]&\HH{1}{G_m}{B''_m}\ar@{->}[0,1]^{\cch{1}}&
	\HH{2}{G_m}{B'_m}\ar@{->}[0,1]&\HH{2}{G_m}{B_m}\\
&	\cdots\ar@{->}[0,1]&\HH{1}{G_n}{B_n}\ar@{->}[-1,0]^{\HHm{1}{\nsu}}\ar@{->}[0,1]&\HH{1}{G_n}{B''_n}\ar@{->}[-1,0]^{\HHm{1}{\nsu''}}\ar@{->}[0,1]^{\cch{1}}&
	\HH{2}{G_n}{B'_n}\ar@{->}[-1,0]^{\HHm{2}{\nsu'}}\ar@{->}[0,1]&\HH{2}{G_n}{B_n}\ar@{->}[-1,0]^{\HHm{2}{\nsu}}
}
\end{equation*}
\end{proposition}
\begin{proof}
By Lemma~\ref{lemma:HHm(j)_and_inf} we can replace $\HHm{1}{j'}$ with $\infl\circ(j')^\ast$, and likewise for the other vertical maps. The result follows from functoriality of inflation and of the maps $j^\ast,(j')^\ast,(j'')^\ast$.
\end{proof}

The above proposition holds for the full long exact sequence of ordinary cohomology extended in every positive degree but we will only need it for the truncated exact sequences up to degree $2$, as in the statement. The following proposition, on the other hand, only holds in degrees smaller or equal than $2$, in general. This is because we assume that $\calB''$ only satisfies \ref{cond:inj} (instead of both~\ref{cond:inj} and~\ref{cond:surj}, as for $\calB$ and $\calB'$). Such formulation of the statement reflects how this result will be used for arithmetic applications (see for instance the exact sequence \ref{eq:long_cohom_DS} which is the basic ingredient of Theorem \ref{thm:i_cinque_dell'apocalisse}). 

\begin{proposition}\label{prop:long_cohom_down}
Let $0\to\calB'\to\calB\to\calB''\to 0$ be an exact sequence of normic systems and assume that $\calB$ and $\calB'$ satisfy~\ref{cond:inj} and~\ref{cond:surj} and that $\calB''$ satisfies~\ref{cond:inj}. Then, for all $m\geq n\geq 0$, the following diagram with exact rows commutes
\[
\xymatrix@C=2.25em{
	0\ar@{->}[0,1]&\HHN{0}{G_m}{B'_m}\ar@{->}[1,0]^{\nsd'}\ar@{->}[0,1]&\HHN{0}{G_m}{B_m}\ar@{->}[1,0]^{\nsd}\ar@{->}[0,1]&\HHN{0}{G_m}{B''_m}\ar@{->}[1,0]^{\nsd''}\ar@{->}[0,1]^{\cch{0}}&\HH{1}{G_m}{B'_m}\ar@{->}[1,0]^{\HHm{1}{\nsd'}}\ar@{->}[0,1]&\cdots\\
	0\ar@{->}[0,1]&\HH{0}{G_n}{B'_n}\ar@{->}[0,1]&\HHN{0}{G_n}{B_n}\ar@{->}[0,1]&\HHN{0}{G_n}{B''_n}\ar@{->}[0,1]^{\cch{0}}&\HH{1}{G_n}{B'_n}\ar@{->}[0,1]&\cdots\\
&	\cdots\ar@{->}[0,1]&\HH{1}{G_m}{B_m}\ar@{->}[1,0]^{\HHm{1}{\nsd}}\ar@{->}[0,1]&\HH{1}{G_m}{B''_m}\ar@{->}[1,0]^{\HHm{1}{\nsd''}}\ar@{->}[0,1]^{\cch{1}}&
	\HH{2}{G_m}{B'_m}\ar@{->}[1,0]^{\HHm{2}{\nsd'}}\ar@{->}[0,1]&\HH{2}{G_m}{B_m}\ar@{->}[1,0]^{\HHm{2}{\nsd}}\\
&	\cdots\ar@{->}[0,1]&\HH{1}{G_n}{B_n}\ar@{->}[0,1]&\HH{1}{G_n}{B''_n}\ar@{->}[0,1]^{\cch{1}}&
	\HH{2}{G_n}{B'_n}\ar@{->}[0,1]&\HH{2}{G_n}{B_n}
}
\]
\end{proposition}
\begin{proof}
First of all, the existence of vertical maps in degree greater than $0$ is guaranteed by the assumptions on the normic systems. 

The commutativity of the two leftmost squares follows immediately from the fact that $\calB'\to\calB$ and $\calB\to\calB''$ are morphisms of normic systems. 

The commutativity of the third square is equivalent to $\tateiso{n}{-1} \circ \cch{0} \circ \nsd'' = \HHm{-1}{\nsd'} \circ \tateiso{m}{-1} \circ \cch{0}$. Fix $b''\in \HHN{0}{G_m}{B_m''}$ and let $b\in B_m$ be an element mapping to $b''$. By the definition of the connecting homomorphism, $\cch{0}(b'')$ is represented by the cocycle $g_m\mapsto (g_m-1)b:=b'$. Similarly, since $\nsd_{m,n}(b)$ maps to $\nsd''_{m,n}(b'')$, the element $\cch{0}\circ \nsd_{m,n}''(b'')$ is represented by the cocycle defined by $g_n\mapsto (g_n-1)\nsd_{m,n}(b)$ and, in turn, this is the map $g_n\mapsto\nsd_{m,n}'(b')$ because the morphisms $\nsd_{m,n}'$ are $G_m$-equivariant. In the notation of Remark~\ref{rmk:cp}, we thus need to show that
\begin{equation*}
\tateiso{n}{-1}\bigl(\xi^{(\nsd_{m,n}'(b'))}_{n}\bigr) = \HHm{-1}{\nsd'} \circ \tateiso{m}{-1} (\xi^{(b')}_{m}).
\end{equation*}
Using the description of $\tateiso{m}{}$ given in Remark~\ref{rmk:cp}, this amounts to 
\begin{equation*}
[\nsd_{m,n}'(b')] = \HHm{-1}{\nsd'} ([b'])
\end{equation*}
which follows from the definition of $\HHm{-1}{\nsd'}$.

As for the remaining part of the diagram, since the definitions of $\HHm{1}{\nsd}$ and $\HHm{2}{\nsd}$ rely on cup products which are isomorphisms compatible with long exact Tate cohomology sequences, it is enough to show that the shifted diagram
\begin{equation}\label{diag:shifted}\begin{split}
\xymatrix{
\HH{-1}{G_m}{B'_m}\ar@{->}[1,0]^{\HHm{-1}{\nsd}}\ar@{->}[0,1]&\HH{-1}{G_m}{B_m}\ar@{->}[1,0]^{\HHm{-1}{\nsd}}\ar@{->}[0,1]&\HH{-1}{G_m}{B''_m}\ar@{->}[1,0]^{\HHm{-1}{\nsd}}\ar@{->}[0,1]^{\cchT{-1}}&
\HH{0}{G_m}{B'_m}\ar@{->}[1,0]^{\HHm{0}{\nsd}}\ar@{->}[0,1]&\HH{0}{G_m}{B_m}\ar@{->}[1,0]^{\HHm{0}{\nsd}}\\
\HH{-1}{G_n}{B'_n}\ar@{->}[0,1]&\HH{-1}{G_n}{B_n}\ar@{->}[0,1]&\HH{-1}{G_n}{B''_n}\ar@{->}[0,1]^{\cchT{-1}}&\HH{0}{G_n}{B'_n}\ar@{->}[0,1]&\HH{0}{G_n}{B_n}
}\end{split}\end{equation}
is commutative. The two leftmost and the fourth squares are commutative as discussed in~\eqref{diag:f*_functorial}. As for the third square, let $b''\in B''_m[\algnm{G_m}]$ and pick $b\in B_m$ mapping to $b''$. Then, the image of $\nsd_{m,n}(b)$ in $B''_n$ is $\nsd_{m,n}''(b'')$ and
\[
\cchT{-1} \circ \HHm{-1}{\nsd}([b'']) = \cchT{-1}([\nsd_{m,n}''(b'')])= \algnm{G_n}([\nsd_{m,n}(b)]).
\] 
On the other hand,
\begin{align*}
\HHm{0}{\nsd}\circ \cchT{-1}([b''])&=\HHm{0}{\nsd}([\algnm{G_m}b])\\
&= (\nsu_{n,m})^{-1}([\algnm{G_m}b])\\
&= [(\nsu_{n,m})^{-1}\circ \nsu_{0,m}\circ\nsd_{m,0}(b)]\\
&=[\nsu_{0,n}\circ \nsd_{n,0}\circ \nsd_{m,n}(b)]=\algnm{G_n}([\nsd_{m,n}(b)])
\end{align*}
showing that the third square in~\eqref{diag:shifted} is commutative.
\end{proof}

\begin{remark} 
Observe that although in \S\ref{subsubsub:HH0} we have considered Tate cohomology of a normic system and defined morphisms $\HHm{0}{\nsu}$ and $\HHm{0}{\nsd}$, Propositions~\ref{prop:long_cohom_up} and~\ref{prop:long_cohom_down} start with \emph{ordinary} cohomology in degree $0$. The appearance of ordinary cohomology is somehow exceptional, because we will mainly be working with Tate cohomology. For instance, in the following section, we will see that for each $i \in \ZZ$, the $i$th Tate cohomology groups of a normic system form a \emph{double system}, a notion to be defined \ibid.
\end{remark}

\subsection{The category of double systems}\label{subsec:double_sys}
In the previous section, we have considered Tate cohomology groups of a normic system, showing that (under some conditions) these can be simultaneously endowed with the structure of direct systems and of inverse systems of $\ZZ_p[\splitG]$-modules with respect to the maps defined in \S\ref{subsubsec:arbitrarydegree}, where $\splitG$ is the group appearing in~\eqref{eq:seq_gal_groups}. In this section, we focus on the study of such ``double systems''.
\begin{definition}\label{def:double_system} 
In the setting of Definition~\ref{def:normic_system}, a $(\Gamma, \{\Gamma_n\}_{n\in\mathbb{N}}, \DD)$-double system (or simply a $\DD$-double system or a double system if the groups are understood) $\calX=(X_n,\dsu_{n,m},\dsd_{m,n})_{m\geq n\geq 0}$ is a collection of $\ZZ_p[\splitG$]-modules $X_n$ which are $\ZZ_p$-torsion and such that 
\begin{itemize}
\item $(X_n,\dsu_{n,m})$ is a direct system;
\item $(X_n,\dsd_{m,n})$ is an inverse system;
\item $\dsu_{n,m}\circ\dsd_{m,n}=\dsd_{m,n}\circ\dsu_{n,m}=\gindex{\Gamma_n}{\Gamma_m}$ for all $m\geq n\geq 0$.
\end{itemize}
The morphisms $\Hom_{\dblsys{\DD}}(\calX,\calX')$ between two double systems $\calX$ and $\calX'$
are collections of $\ZZ[\splitG]$-homomorphisms $f_n\colon X_n\to X'_n$ which are both morphisms of direct and of inverse systems. This defines the category $\dblsys{\DD}$ of $\DD$-double systems.
\end{definition}
The results in~\S\ref{subsec:normic_sys} yield the following

\begin{proposition}\label{prop:cohomds}
Let $\calB$ be a normic system. If $\calB$ satisfies~\ref{cond:inj}, then
\begin{equation}
\sysHH{i}{\calB}=\bigl(\HH{i}{G_n}{B_n},\HHm{i}{\nsu_{n,m}},\HHm{i}{\nsd_{m,n}}\bigr)_{m\geq n\geq 0}
\end{equation}
is a double system for every odd integer $i$. If condition~\ref{cond:surj} is also satisfied, then $\sysHH{i}{\calB}$ is a double system for all $i\in\ZZ$.
\end{proposition}

\begin{remark} 
The category $\dblsys{\DD}$ is abelian: the proof  is analogous to that of Lemma~\ref{lemma:normsys_abelian}, by replacing the requirement $\vec{v}_{\ell,i}\cev{v}_{i,\ell}=\sum_{a=1}^{\gindex{\Gamma_\ell}{\Gamma_i}}\topgen^{ap^\ell}$ in~\eqref{eq:def_R} with $\vec{v}_{\ell,i}\cev{v}_{i,\ell}=\gindex{\Gamma_\ell}{\Gamma_i}$ for $\ell\leq i$. Alternatively, one can observe that $\dblsys{\DD}$ is a full subcategory of $\normsys$ closed under direct sums, kernels and cokernels and it is therefore abelian, thanks to~\cite[Proposition~5.92]{Rot09}.
\end{remark}

Let $\calX=(X_n,\dsu_{n,m},\dsd_{m,n})_{m\geq n\geq 0}$ be a double system, and write $\dsu_n\colon X_n\to\varinjlim X_n$ for the direct limit of the maps $\dsu_{n,m}$ (for $m\geq n$). Similarly, write $\dsd_n\colon\varprojlim X_n\to X_n$ for the inverse limit of $\dsd_{m,n}$. Let
\begin{equation}\label{eq:def_capzero_n}
\capzero{X_n}=\kernel\bigr(\dsu_n\colon X_n\longrightarrow\varinjlim X_n\bigl)\subseteq X_n\qquad\text { and }\qquad\capbar{X_n}=X_n/\capzero{X_n}
\end{equation}
as well as
\begin{equation}\label{eq:def_nuczero_n}
\nuczero{X_n}=\im\bigr(\dsd_n\colon \varprojlim X_n\longrightarrow X_n\bigl)\subseteq X_n\qquad\text{ and }\qquad
\nucbar{X_n}=\cokernel\bigr(\dsd_n\colon \varprojlim X_n\longrightarrow X_n\bigl)=X_n/\nuczero{X_n}.
\end{equation}
The next proposition shows that $\capzero{\calX}=(\capzero{X_n})$, $\nuczero{\calX}=(\nuczero{X_n})$, as well as $\capbar{\calX}=(\capbar{X_n})$ and $\nucbar{\calX}=(\nucbar{X_n})$ inherit the structure of double systems from $\calX$. It also shows that the assignments $\calX\mapsto \capzero{\calX}$, $\calX\mapsto \nuczero{\calX}$, $\calX\mapsto \capbar{\calX}$ and $\calX\mapsto\nucbar{\calX}$ are functorial.
\begin{proposition}\label{prop:univ_capzero} 
Let $\calX=(X_n,\dsu_{n,m},\dsd_{m,n})_{m\geq n\geq 0}$ be a double system. The restrictions of $\dsd_{m,n}$ and $\dsu_{n,m}$ endow $\capzero{\calX}$ and $\nuczero{\calX}$---and thus $\capbar{\calX}$ and $\nucbar{\calX}$---with a structure of double system. These double systems satisfy:
\begin{enumerate}[label=\textit{\roman*}\textup{)}]
\item  \label{point:prop_capzero:functor} The formation of $\capzero{\calX}$ and $\capbar{\calX}$ is functorial, and the functor $\calX\mapsto \capzero{\calX}$ is left exact. Similarly, the formation of $\nuczero{\calX}$ and $\nucbar{\calX}$ is functorial, and the functor $\calX\mapsto \nucbar{\calX}$ is right exact.
\item \label{point:prop_capzero:sono_inj}  Write $\capbar{\calX}=(\capbar{X}_n,\capbar{\dsu}_{n,m},\capbar{\dsd}_{m,n})_{m\geq n\geq 0}$ and $\nuczero{\calX}=(\nuczero{X}_n,\nuczero{\dsu}_{n,m},\nuczero{\dsd}_{m,n})_{m\geq n\geq 0}$. Then $\capbar{\dsu}_{n,m}$ is injective and $\nuczero{\dsd}_{m,n}$ is surjective for all $m\geq n\geq 0$.
\item \label{point:prop_capzero:univ} For every double system $\calX'=(X_n',\dsu_{n,m}',\dsd_{m,n}')_{m\geq n\geq 0}$ such that $\dsu_{n,m}'$ is injective for $m\geq n\geq 0$, and every morphism $f\colon \calX\to \calX'$, there exists a unique map $\capbar{f}\colon \capbar{\calX}\to \calX'$ making the following diagram commute
\[\xymatrix{
\calX\ar@{->}[0,2]^{f}\ar@{->>}[1,1]&&\calX'\\
&\capbar{\calX}\ar@{->}[-1,1]_{\capbar{f}}
}\]
Analogously, given a double system $\calX''=(X_n',\dsu_{n,m}'',\dsd_{m,n}'')_{m\geq n\geq 0}$ such that $\dsd_{m,n}''$ is surjective for $m\geq n\geq 0$ \textup{(}and so $\calX''=\nuczero{\calX''}$\textup{)}, every morphism $g\colon \calX''\to \calX$ factors through the subsystem $\nuczero{\calX}$, making the following diagram commute
\[\xymatrix{
\calX''\ar@{->}[0,2]^{g}\ar@{->}[1,1]_(.4){\nuczero{g}}&&\calX\\
&\nuczero{\calX}\ar@{^(->}[-1,1]
}\]
\end{enumerate}
\end{proposition}
\begin{proof} 
We claim that the restriction of $\dsu_{n,m}$ (\rsp~of $\dsd_{m,n}$) maps $\capzero{X_n}$ to $\capzero{X_m}$ (\rsp$\capzero{X_m}$ to $\capzero{X_n}$). This is obvious for $\dsu_{n,m}$, thanks to the relation $\dsu_m\circ\dsu_{n,m}=\dsu_n$. Concerning $\dsd_{m,n}$, the relation $\dsu_{n,m}\circ\dsd_{m,n}=\gindex{\Gamma_n}{\Gamma_m}$ shows that for all $x\in\capzero{X_m}$ it holds
\[
\dsu_n\bigl(\dsd_{m,n}(x)\bigr)=\dsu_m\circ\dsu_{n,m}\circ\dsd_{m,n}(x)=\dsu_m\bigl(\gindex{\Gamma_n}{\Gamma_m}x\bigr)=0 
\]
and hence $\dsd_{m,n}(x)\in\capzero{X_n}$. The claim implies that, defining $\capzero{\dsu_{n,m}}$ (\rsp $\capzero{\dsd_{m,n}}$) to be the restriction of $\dsu_{n,m}$ (\rsp of $\dsd_{m,n}$) to $\capzero{X_n}$, we obtain a double system $\capzero{\calX}$ which is a subobject of $\calX$. Moding out $\calX$ by $\capzero{\calX}$ defines the double system $\capbar{\calX}$. The claim that $\nuczero{\calX}$ is a double system which is a subobject of $\calX$ such that $\nucbar{\calX}=\calX/\nuczero{\calX}$ (in the category $\dblsys{\DD}$) is analogous: the fact that the restriction $\nuczero{\dsd_{m,n}}$ to $\nuczero{X_m}$ takes values in $\nuczero{X_n}$ is automatic, and the relation $\dsu_{n,m}\circ\dsd_{m,n}=\gindex{\Gamma_n}{\Gamma_m}$ ensures $\nuczero{\dsu_{n,m}}(\nuczero{X_n})\subseteq \nuczero{X_m}$.

Functoriality as well as the claimed exactness properties of $\calX\mapsto\capzero{\calX}$ and of $\calX\mapsto\nucbar{\calX}$ are straightforward: given a morphism $\calX\to \calX'$, the rightmost square of the diagram~\eqref{diag:functor_capbar} below is commutative for all $n\geq 0$, implying the existence of the leftmost vertical morphism (making the whole diagram commute):
\begin{equation}\label{diag:functor_capbar}\begin{split}
\xymatrix{
0\ar@{->}[0,1]&\capzero{X_n}\ar@{->}[0,1]\ar@{->}[1,0]&X_n\ar@{->}[1,0]\ar@{->}[0,1]&\varinjlim X_n\ar@{->}[1,0]\\
0\ar@{->}[0,1]&\capzero{X_n'}\ar@{->}[0,1]&X_n'\ar@{->}[0,1]&\varinjlim{X_n'}
}\end{split}
\end{equation}
This shows that $\calX\mapsto \capzero{\calX}$ is functorial. Moreover, as discussed in~\cite[\S1.6]{Gro57}, a morphism of double systems is injective (\rsp surjective) if an only if each of its components has this property. Hence, if $\calX\hookrightarrow\calX'$ is injective,~\eqref{diag:functor_capbar} shows that each component of the induced morphism $\capzero{\calX}\to\capzero{\calX'}$ is injective, so $\capzero{\calX}\hookrightarrow\capzero{\calX}$ is again injective. This is well-known to be equivalent to the exactness on the left of the functor $\calX\mapsto\capzero{\calX}$ (since~$\dblsys{\DD}$ is abelian), establishing the first part of~\ref{point:prop_capzero:functor}. Analogously, a surjection $\calX\twoheadrightarrow\calX'$ induces a commutative diagram of exact rows
\[\xymatrix{
\varprojlim X_n\ar@{->}[0,1]\ar@{->}[1,0]&X_n\ar@{->>}[1,0]\ar@{->}[0,1]&\nucbar{X_n}\ar@{->}[0,1]\ar@{->}[1,0]&0\\
\varprojlim X_n'\ar@{->}[0,1]&X_n'\ar@{->}[0,1]&\nucbar{X_n'}\ar@{->}[0,1]&0
}\]
showing that each component of the induced morphism $\capzero{\calX}\to\capzero{\calX'}$ is surjective, finishing the proof of~\ref{point:prop_capzero:functor}.

To prove~\ref{point:prop_capzero:sono_inj}, observe that $\capbar{\dsu}_{m,n}$ is injective for all $m\geq n\geq 0$ by construction: indeed, for all $[x]=x\pmod{\capzero{X_n}}\in\capbar{X}_n$, the definition of $\capbar{\dsu}_{n,m}$ is $\capbar{\dsu}_{n,m}[x]=\dsu_{n,m}(x)\pmod{\capzero{X_m}}$. Then $\capbar{\dsu}_{n,m}[x]=0$ implies $\dsu_{n,m}(x)\in \capzero{X_m}$, which means $\dsu_m\dsu_{n,m}(x)=\dsu_{n}(x)=0$ and hence $x\in \capzero{X_n}$ and $[x]=0$, establishing the injectivity of $\capbar{\dsu}_{n,m}$. The proof of the surjectivity of $\nuczero{\dsd}_{m,n}$ is analogous.

We are left with the proof of~\ref{point:prop_capzero:univ}. Let $\calX'=(X_n',\dsu_{n,m}',\dsd_{m,n}')$ be such that $\dsu_{n,m}'$ is injective for $m\geq n\geq 0$, and let $f\colon \calX\to \calX'$. Write $f_\infty\colon \varinjlim X_n\to\varinjlim X_n'$ for the direct limit of the maps $f_n$. The injectivity of all $\dsu_{n,m}'$ implies that $\dsu_n'\colon X_n'\to \varinjlim X_n'$ is also injective. It follows that $f_n(\capzero{X_n})=0$, because
\[
\dsu_n'\bigl(f_n(x)\bigr)=f_\infty\bigl(\dsu_n(x)\bigr)=f_\infty(0)=0\Longrightarrow f_n(x)=0\qquad\text{ for all }x\in\capzero{X_n}.
\]
In particular, for all $n$, there exists a unique map $\capbar{f}_n\colon X_n/\capzero{X_n}=\capbar{X_n}\to X_n'$ making the following diagram commute
\[\xymatrix{
X_n\ar@{->}[0,2]^{f_n}\ar@{->>}[1,1]&&X_n'\\
&\capbar{X_n}\ar@{->}[-1,1]
_{\capbar{f}_n}
}\]
We need to check that the collection $\capbar{f}_n$ defines a (necessarily unique) morphism $\capbar{f}\in\Hom_{\dblsys{\DD}}(\capbar{\calX},\calX')$, namely that
\begin{align}
\dsd'_{m,n}\circ \capbar{f}_m&=\capbar{f}_n\circ\capbar{\dsd}_{m,n} \label{eq:communiv1}\\ 
\dsu'_{n,m}\circ \capbar{f}_n&=\capbar{f}_m\circ \capbar{\dsu}_{n,m}.\label{eq:communiv2}
\end{align}
This is tantamount to showing the commutativity of the right oblique squares in the following diagrams:
\begin{equation*}
\begin{array}{lcr}
\xymatrix{
X_m\ar@{->}[2,0]_{\dsd_{m,n}}\ar@{->>}[1,1]\ar@{->}[0,2]^(.4){f_m}&&X_m'\ar@{->}[2,0]^{\dsd_{m,n}'}\\
&\capbar{X_m}\ar@{->}[-1,1]^(.45){\capbar{f_m}}\ar@{->}[2,0]&&\\
X_n\ar@{->>}[1,1]\ar@{->}[0,2]^(.35){f_n}|\hole|\hole&&X_n'&\\
&\capbar{X_n}\ar@{->}[-1,1]_{\capbar{f_n}}
}
&\parbox[t][10.5em][c]{6em}{\begin{center}\text{and}\end{center}}&
\xymatrix{
X_m\ar@{->>}[1,1]\ar@{->}[0,2]^(.4){f_m}&&X_m'\\
&\capbar{X_m}\ar@{->}[-1,1]^(.45){\capbar{f_m}}&&\\
X_n\ar@{->}[-2,0]^{\dsu_{n,m}}\ar@{->>}[1,1]\ar@{->}[0,2]^(.35){f_n}|\hole|\hole&&X_n'\ar@{->}[-2,0]_{\dsu_{n,m}'}&\\
&\capbar{X_n}\ar@{->}[-2,0]\ar@{->}[-1,1]_{\capbar{f_n}}}
\end{array}\end{equation*}

Since the left oblique squares commute, and so do the straight squares and triangles, one obtains that \eqref{eq:communiv1} (\rsp~\eqref{eq:communiv2}) holds after composing with $X_m \to \capbar{X_m}$ (\rsp~$X_n\to \capbar{X_n}$). Since $X_m \to \capbar{X_m}$ and $X_n \to \capbar{X_n}$ are surjective, the required commutativity follows. The argument for the universal property enjoyed by $\nuczero{\calX}$ is analogous. This finishes the proof of the proposition.
\end{proof}
From now on, we drop the $\capzero{\empty}$ from the transition morphisms $\capzero{\dsd_{m,n}}$ and $\capzero{\dsu_{n,m}}$ of $\capzero{\calX}$, simply writing $\dsd_{m,n}$ and $\dsu_{n,m}$ for their restrictions to $\capzero{X_m}$ and to $\capzero{X_n}$, respectively. The same convention is adopted for the other three double systems $\capbar{\calX},\nuczero{\calX}$ and $\nucbar{\calX}$.
\begin{remark}
We are grateful to the referee for pointing out that Proposition~\ref{prop:univ_capzero}-\ref{point:prop_capzero:univ} can be interpreted as saying that $\calX\mapsto\capbar{\calX}$ and $\calX\mapsto\nuczero{\calX}$, suitably restricted to certain subcategories of $\dblsys{\DD}$, are adjoint to a forgetful functor. Since we have no use for this in our work, we omit the details.
\end{remark}
\begin{remark}\label{rmk:injsurj_imply_iso}
If $\calX$ is a double system such that $\dsd_{m,n}$ is injective for all $m\geq n\gg 0$, then $\calX=\capzero{\calX}$; similarly, if~$\dsu_{n,m}$ is surjective for all $m\geq n\gg 0$ and $X_n$ has bounded exponent for every $n$, then $\calX=\nucbar{\calX}$. To see the first implication, let $x\in X_n$ and pick $m\geq n$ such that $\gindex{\Gamma_n}{\Gamma_m}x=0$. Writing $\gindex{\Gamma_n}{\Gamma_m}=\dsd_{m,n}\circ\dsu_{n,m}$ and using that $\dsd_{m,n}$ is injective, we deduce $\dsu_{n,m}(x)=0$, hence $\dsu_n(x)=0$ and $x\in\capzero{X_n}$. This shows $\calX=\capzero{\calX}$. Similarly, suppose that $\dsu_{n,m}$ is surjective for $m\geq n\gg 0$, and let $x\in\nuczero{ X_n}$. Pick $m\geq n$ such that $\gindex{\Gamma_n}{\Gamma_m}X_n=0$ and let $y\in X_m$ be such that $\dsd_{m,n}(y)=x$, which exists because $x\in\nuczero{X_n}$. Surjectivity of $\dsu_{n,m}$ guarantees the existence of an element $z\in X_n$ such that $\dsu_{n,m}(z)=y$, and $x=\dsd_{m,n}\circ\dsu_{n,m}(z)=\gindex{\Gamma_n}{\Gamma_m}z=0$. Then $\nuczero{\calX}=0$ and $\calX=\nucbar{\calX}$.
\end{remark}
Before stating the next lemma, observe that if $\calX=(X_n)$ is a double system such that the exponents of $X_n$ are bounded independently of $n$, then \begin{equation}\label{eq:nuc_sub_cap}
\nuczero{\calX}\subseteq\capzero{\calX}.
\end{equation}
Indeed, for $n\geq 0$ and $x\in \nuczero{X_n}$, one can write $x=\dsd_{n}(\xi)$ for some $\xi\in\varprojlim X_m$. The element $\xi$ is torsion, because every $X_m$ has bounded exponent, and we pick $m$ such that $\gindex{\Gamma_n}{\Gamma_m}\xi =0$. Then
\[
\dsu_{n}(x)=\dsu_{m}\circ\dsu_{n,m}\circ\dsd_{m,n}\circ\dsd_{m}(\xi)=\dsu_{m}\circ\dsd_{m}\bigl(\gindex{\Gamma_n}{\Gamma_m}\xi\bigr)=0
\]
showing $x\in\capzero{X_n}$, and thus~\eqref{eq:nuc_sub_cap}.
\begin{lemma}\label{lemma:eqB} 
Let $\calX=(X_n,\dsu_{n,m},\dsd_{m,n})_{m\geq n\geq 0}$ be a $\DD$-double system such that the orders $\gorder{X_n}$ are bounded independently of~$n$. The following conditions are equivalent:
\begin{enumerate}[label=$\bddsys{\DD}$-$\mathbf{\arabic*})$]
\item\label{point:lemma_eqB:pi} $\dsd_{m,n}\colon \capzero{X_m}\to \capzero{X_n}$ is an isomorphism for all $m\geq n\gg 0$\textup{;}
\item\label{point:lemma_eqB:phi} $\dsu_{n,m}\colon \nucbar{X_n}\to \nucbar{X_m}$ is an isomorphism for all $m\geq n\gg 0$\textup{;}
\item\label{point:lemma_eqB:cap=nuc} $\nuczero{X_n}=\capzero{X_n}$, or equivalently $\nucbar{X_n}=\capbar{X_n}$, for all $n\gg 0$.
\end{enumerate}
When any of the above conditions is met, the orders $\gorder{X_n}$ are eventually constant, and both $\dsd_{m,n}\colon \nuczero{X_m}\to \nuczero{X_n}$ and $\dsu_{n,m}\colon \capbar{X_n}\to \capbar{X_m}$ are isomorphisms for all $m\geq n\gg 0$.
\end{lemma}
\begin{proof} 
To prove that~\ref{point:lemma_eqB:pi} implies~\ref{point:lemma_eqB:cap=nuc}, observe that it is enough to show $\capzero{X_n}\subseteq \nuczero{X_n}$ for all $n\gg 0$, in light of~\eqref{eq:nuc_sub_cap}. By~\ref{point:lemma_eqB:pi} we can find $n$ such that $\dsd_{m,n}\colon \capzero{X_m}\xrightarrow{\sim} \capzero{X_n}$ for all $m\geq n$, and let $x\in\capzero{X_n}$. For all~$m\geq n$ there exists a unique $x_m\in\capzero{X_m}$ such that $\dsd_{m,n}(x_m)=x$. It follows that $\xi=\lim x_m\in\varprojlim X_m$ exists and $x=\dsd_{n}(\xi)$, showing $x\in\nuczero{X_n}$.

Assuming~\ref{point:lemma_eqB:cap=nuc}, Proposition~\ref{prop:univ_capzero}-\ref{point:prop_capzero:sono_inj} implies that $\dsu_{n,m}\colon\nucbar{X_n}\to\nucbar{X_m}$ is injective for all $m\geq n\gg 0$. It follows that $n\mapsto \gorder{\nucbar{X_n}}$ is an non-decreasing sequence of bounded, positive integers and thus eventually constant. As an injection of finite groups of the same order is an isomorphism, this proves the implication \ref{point:lemma_eqB:cap=nuc}$\Rightarrow$\ref{point:lemma_eqB:phi}.

Finally, assume~\ref{point:lemma_eqB:phi} and let $n_0$ be such that $\dsu_{n,m}$ is an isomorphism for all $m\geq n\geq n_0$. Consider the shifted double system
\[
\calX'=\bigr(X_n'=X_{n+n_0},\dsu_{n,m}'=\dsu_{n+n_0,m+n_0},\dsd_{m,n}'=\dsd_{m+n_0,n+n_0}\bigl)_{m\geq n\geq 0}.
\]
Then $\calX'$ satisfies $\nuczero{X_n'}=\nuczero{X_{n+n_0}}$ and $\capzero{X_n'}=\capzero{X_{n+n_0}}$ because $\varprojlim X_n'=\varprojlim X_n$ and $\varinjlim X_n=\varinjlim X_n'$. It follows that $\nucbar{X_n'}=\nucbar{X_{n+n_0}}$ and $\capbar{X_n'}=\capbar{X_{n+n_0}}$ for all $n\geq 0$, so $\calX'$ satisfies~\ref{point:lemma_eqB:phi} for all $m\geq n\geq 0$. By Proposition~\ref{prop:univ_capzero}-\ref{point:prop_capzero:univ}, the natural surjection $\calX'\twoheadrightarrow \nucbar{\calX'}$ factors as
\[\xymatrix{
\calX'\ar@{->>}[0,2]\ar@{->>}[1,1]&&\nucbar{\calX'}\\
&\capbar{\calX'}\ar@{->>}[-1,1]
}\]
On the other hand, the inclusion~\eqref{eq:nuc_sub_cap} induces a surjection $\nucbar{\calX'}\twoheadrightarrow \capbar{\calX'}$. Since all groups $\capbar{X_n'}$ and $\nucbar{X_n'}$ are finite, the existence of surjections $\capbar{X_n'}\twoheadrightarrow \nucbar{X_n'}$ and $\nucbar{X_n'}\twoheadrightarrow\capbar{X_n'}$ implies that these surjections are isomorphisms, proving~\ref{point:lemma_eqB:cap=nuc} for the system $\calX'$. Upon reindexing, this implies~\ref{point:lemma_eqB:cap=nuc} for $\calX$.

To establish the statement concerning the stabilisation of the orders $\gorder{X_n}$, write $\gorder{X_n}=\gorder{\capzero{X_n}}\cdot\gorder{\capbar{X_n}}$. Thanks to~\ref{point:lemma_eqB:cap=nuc} this coincides with $\gorder{\capzero{X_n}}\cdot\gorder{\nucbar{X_n}}$, and both the orders $\gorder{\capzero{X_n}}$ and $\gorder{\nucbar{X_n}}$ are eventually constant because of~\ref{point:lemma_eqB:pi} and~~\ref{point:lemma_eqB:phi}. Replacing in the isomorphisms in~\ref{point:lemma_eqB:pi} and~~\ref{point:lemma_eqB:phi} the equalities of~\ref{point:lemma_eqB:cap=nuc} yields the final statement.
\end{proof}
\begin{definition}\label{def:bddsys}
Let $\bddsys{\DD}$ be the full subcategory of $\dblsys{\DD}$ whose objects are the $\DD$-double systems of finite $p$-groups whose order is bounded independently of $n$ and that satisfy any of the equivalent conditions~\ref{point:lemma_eqB:pi}--\ref{point:lemma_eqB:cap=nuc} of Lemma~\ref{lemma:eqB}. In particular, the order of the components of the system is in fact eventually constant, independently of $n\gg 0$.
\end{definition}

\begin{remark}
In Example~\ref{ex:asym_not_B}, we construct a double system of finite $p$-groups of constant order which is not an object of $\bddsys{\DD}$.
\end{remark}
\begin{example}\label{ex:thanks_to_referee}
Suppose that $\Gamma$ is filtered by $\Gamma_n=\Gamma^{p^n}$ for all $n\geq 0$ and consider the double system $\calX=(\ZZ/p,\id{},0)$, where all components $X_n$ coincide with $\ZZ/p$, the ascending morphisms are the identity and the descending ones are the $0$ map. One computes that $\varprojlim \calX=0$ and $\varinjlim \calX=\ZZ/p$, so that $\capzero{X_n}=\nuczero{X_n}=0$, and conditions~\ref{point:lemma_eqB:pi} and~\ref{point:lemma_eqB:cap=nuc} (and thus~\ref{point:lemma_eqB:phi}) are satisfied, showing that $\calX$ is an object of $\bddsys{\DD}$.

Similarly, one can consider the system $\calY=(\ZZ/p,0,\id{})$, that is again an object of $\bddsys{\DD}$ because $\capzero{Y_n}=\nuczero{Y_n}=Y_n$, so $\capbar{Y_n}=0$ and conditions~\ref{point:lemma_eqB:cap=nuc} and~\ref{point:lemma_eqB:phi} are satisfied (and so is~\ref{point:lemma_eqB:pi}, then).
\end{example}
\begin{example}\label{example:shift}
In this example, denote for simplicity the index $\gindex{\Gamma_n}{\Gamma_m}$ by $\filindex{n}{m}$ and write $\filindex[\empty]{}{n}$ for $\gindex{\Gamma}{\Gamma_n}$; similarly, write $g$ for $g(\Gamma)$. In case the filtration $\{\Gamma_n\}_{n\geq 0}$ coincides with $\{\Gamma^{p^n}\}_{n\geq 0}$, these become $\filindex{n}{m}=p^{m-n}$, $\filindex[\empty]{}{n}=p^n$ and $g=0$, and the reader may keep this particular case in mind.

Consider the $\DD$-double system $\calZ(\Gamma)$ (denoted simply $\calZ$ if the filtered group $\Gamma$ is understood) defined as
\[
\calZ(\Gamma)=\bigr(\ZZ/\filindex[\empty]{}{n},\cdot \filindex{n}{m}\colon\ZZ/\filindex[\empty]{}{n}\hookrightarrow\ZZ/\filindex[\empty]{}{m},\pr\bigl)_{m\geq n\geq 0}.
\]
Here, for $m\geq n$, $\cdot \filindex{n}{m}$ sends $x \mod{\filindex[\empty]{}{n}}\ZZ$ to $\filindex{n}{m}x \mod{\filindex[\empty]{}{m}}\ZZ$ and $\pr\colon\ZZ/\filindex[\empty]{}{m}\to\ZZ/\filindex[\empty]{}{n}$ denotes the canonical projection (we omit the indexes in $\pr$ to avoid notation overload); all components have trivial action of $\splitG$. Let $n_0\geq 1$: define the shift $\calZ_{[n_0]}(\Gamma)=\calZ_{[n_0]}$ of $\calZ$ as the double system whose components are
\[
(\calZ_{[n_0]})_n=\begin{cases}
				\ZZ/\filindex[\empty]{}{n}p^{-n_0}&\text{ for all }n \geq  g+n_0\\
				0&\text{ for all }n < g+n_0
				\end{cases}
\]
with transition morphisms $\dsu_{n,n+1}=\cdot p$ and $\dsd_{n+1,n}=\pr$ for $n\geq n_0$, and the trivial maps otherwise. Observe that the requirement $n\geq g+n_0$ ensures that $\filindex[\empty]{}{n}p^{-n_0}\in\NN$ and $\filindex[\empty]{}{n+1}=p\filindex[\empty]{}{n}$. We intend to show that, although $\calZ_{[n_0]}$ is not isomorphic to $\calZ$ (none of the components are, in fact, and since kernels and cokernels are computed component-wise in $\dblsys{\DD}$, the systems cannot be isomorphic), their ``difference'' lies in $\bddsys{\DD}$ in a precise sense.

Define a map $\iota:\calZ_{[n_0]}\to \calZ$ by $\iota_n=0$ if $n<g+n_0$ and  $\iota_n= \cdot p^{n_0}\colon(\calZ_{[n_0]})_n=\ZZ/\filindex[\empty]{}{n}p^{-n_0}\hookrightarrow \ZZ/\filindex[\empty]{}{n}=(\calZ)_n$ if $n\geq g+n_0$. One easily checks compatibility of $\iota_n$ with transition maps, so $\iota$ defines an injection $\iota\colon \calZ_{[n_0]}\hookrightarrow \calZ$ in the category $\dblsys{\DD}$. The cokernel of $\iota$ is the double system of the cokernels, hence
\[
(\cokernel\iota)_n=	\begin{cases}
					\ZZ/p^{n_0}&\text{ if }n\geq g+n_0\\
					\ZZ/\filindex[\empty]{n}{}&\text{ if }n< g+n_0
					\end{cases}
\]
with transition morphisms
\[
\dsu_{n,n+1}=p
\qquad\text{ and }\qquad
\dsd_{n+1,n}=	\begin{cases}
				\id{\empty}\colon\ZZ/p^{n_0}\longrightarrow\ZZ/p^{n_0}&\text{ if }n\geq g+n_0\\
				\pr\colon\ZZ/\filindex[\empty]{n+1}{}\longrightarrow\ZZ/\filindex[\empty]{n}{}&\text{ if }n< g+n_0
				\end{cases}
\]
Since the cokernel is non-trivial, the morphism $\iota$ is not an isomorphism in $\dblsys{\DD}$. Yet, the direct limit $\varinjlim(\cokernel\iota)$ is trivial, because for all $n\geq g+n_0$ and all $x\in(\cokernel\iota)_n$, $\dsu_{n,2n}(x)=p^{n}x=0$. In particular, $(\capzero{\cokernel\iota})_n=(\cokernel\iota)_n$ for $n\geq n_0$ and thus condition~\ref{point:lemma_eqB:pi} is fulfilled, showing $\cokernel\iota\in\bddsys{\DD}$.

Rather than considering the injection $\iota\colon\calZ_{[n_0]}\to\calZ$ one can look at the morphism $\phi\colon \calZ\to\calZ_{[n_0]}$ defined component-wise by the projections $(\calZ)_n=\ZZ/\filindex[\empty]{n}{}\to\ZZ/\filindex[\empty]{n}{}p^{-n_0}=(\calZ_{[n_0]})_n$, for $n\geq g+n_0$. One can verify that the double system $(\kernel\phi,\dsu',\dsd')$ satisfies $\varprojlim \kernel\phi=0$ and $\dsu'=\id{\empty}$ for all $n\gg 0$, showing that $\kernel\phi\in\bddsys{\DD}$ by condition~\ref{point:lemma_eqB:phi}.
\end{example}

The following proposition provides a class of double systems belonging to $\bddsys{\DD}$.

\begin{proposition}\label{prop:trivialbdd} 
Let $\calX=(X_n,\dsu_{n,m},\dsd_{m,n})_{m\geq n\geq 0}$ be a double system of \emph{finite} $\ZZ_p[\splitG]$-modules such that either $\dsu_{n,m}$ is an isomorphism for $m\geq n\gg 0$, or $\dsd_{m,n}$ is an isomorphism for $m\geq n\gg 0$: then $\calX\in\bddsys{\DD}$.
\end{proposition}
\begin{proof}
Observe that if $\dsu_{n,m}$ (\rsp $\dsd_{m,n}$) is eventually an isomorphism, then the orders $\gorder{X_n}$ are bounded independently of $n$ and the direct limit $\dsu_n$ (\rsp the inverse limit $\dsd_{n}$) is also an isomorphism for $n\gg 0$. In particular $\capzero{X_n}=0$ (\rsp $\nucbar{X_n}=0$) for $n\gg 0$ and condition~\ref{point:lemma_eqB:pi} (\rsp~condition~\ref{point:lemma_eqB:phi}) of Lemma~\ref{lemma:eqB} is satisfied, so $\calX\in\bddsys{\DD}$.
\end{proof}

\begin{remark}
Proposition~\ref{prop:trivialbdd} shows that, for a double system of finite $\ZZ_p[\splitG]$-modules, the condition of being in $\bddsys{\DD}$ only depends upon the underlying direct, or inverse, system. This will be useful for our arithmetic applications, especially at the end of the proof of Theorem~\ref{thm:i_cinque_dell'apocalisse}.
\end{remark}

We recall the following result concerning quotient categories, for which we refer to~\cite[Ch.~III,\S1]{Gab62}:

\begin{defthm*} A full subcategory $\mathbf{C}$ of an abelian category $\mathbf{A}$ is called thick \textup{(}\emph{épaisse} in French\textup{)} if for every exact sequence
\[
0\longrightarrow X'\longrightarrow X\longrightarrow X''\longrightarrow 0
\]
in $\mathbf{A}$, the object $X$ belongs to $\mathbf{C}$ if and only if both $X'$ and $X''$ belong to $\mathbf{C}$. For every thick subcategory $\mathbf{C}\subseteq\mathbf{A}$ there exists a quotient category $\mathbf{A}/\mathbf{C}$, which is abelian and comes with a full localising functor
\[
\mathbf{A}\longrightarrow \mathbf{A}/\mathbf{C}.
\]
The objects of $\mathbf{A}/\mathbf{C}$ are the objects of $\mathbf{A}$. Given two objects $X,X'$ in $\mathbf{A}$ and a morphism $f\colon X\to X'$, its image \textup{(}still denoted $f$\textup{)} in $\Hom_{\mathbf{A}/\mathbf{C}}(X,X')$ is a monomorphism \textup{(}\rsp epimorphism, 	isomorphism\textup{)} if and only if $\kernel{f}\in\mathbf{C}$ \textup{(}\rsp $\cokernel{f}\in\mathbf{C}$, $\kernel{f}$ and $\cokernel{f}\in\mathbf{C}$ \textup{)}. If this is the case, we say that $f$ is a $\mathbf{C}$-monomorphism \textup{(}\rsp $\mathbf{C}$-epimorphism, $\mathbf{C}$-isomorphism\textup{)}.
\end{defthm*}

The main technical result of this section is the following proposition.

\begin{proposition}\label{prop:is_thick}
The category $\bddsys{\DD}$ is a thick subcategory of $\dblsys{\DD}$.
\end{proposition}
\begin{proof}
Let $0\to \calX'\to\calX\to\calX''\to 0$ be an exact sequence in $\dblsys{\DD}$. For every $n$ we have $\gorder{X_n}=\gorder{X_n'}\cdot\gorder{X_n''}$, so the order $\gorder{X_n}$ is bounded if and only if $\gorder{X_n'}$ and $\gorder{X_n''}$ are bounded. 

By Proposition~\ref{prop:univ_capzero}-\ref{point:prop_capzero:functor}, for every $m\geq n\geq 0$, there are commutative diagrams
\begin{equation}\label{diag:thick}\begin{split}
\xymatrix{
	0\ar@{->}[0,1]&\capzero{X_m'}\ar@{->}[1,0]^{\dsd_{m,n}'}\ar@{->}[0,1]&\capzero{X_m}\ar@{->}[1,0]^{\dsd_{m,n}}\ar@{->}[0,1]&\capzero{X_m''}\ar@{->}[1,0]^{\dsd_{m,n}''}\ar@{->}[0,1]&0\\
	0\ar@{->}[0,1]&\capzero{X_n'}\ar@{->}[0,1]&\capzero{X_n}\ar@{->}[0,1]&\capzero{X_n''}\ar@{->}[0,1]&0
}\end{split}
\end{equation}
whose rows are exact, except possibly on the right.

Suppose that $\calX'$ and $\calX''$ are objects in $\bddsys{\DD}$, and fix $c\geq g(\Gamma)$ such that $p^c X_n=p^cX_n'=p^cX_n''=0$ for all~$n\geq 0$. We claim that since $\calX''\in\bddsys{\DD}$, rows of \eqref{diag:thick} are exact on the right for $n$ large enough. For any~$n\geq 0$ large enough, \ref{point:lemma_eqB:pi} ensures that $\dsd_{\ell,n}''\colon \capzero{X_\ell''}\to\capzero{X_n''}$ is an isomorphism for all $\ell\geq n$; we fix $\ell\geq n+c$, so $\gindex{\Gamma_n}{\Gamma_\ell}X_\ell=p^{\ell-n}X_\ell=0$. Let $[x]\in\capzero{X_n''}$. There exists $[y]\in \capzero{X_\ell''}$ such that $[x]=\dsd_{\ell,n}''[y]$ and if $y\in X_\ell$ is a lift of $[y]$, then $x=\dsd_{\ell,n}(y)\in X_n$ satisfies $x\pmod {X_n'}=[x]$. Moreover, $\dsu_n(x)=\dsu_\ell\dsu_{n,\ell}\dsd_{\ell,n}(y)=\dsu_\ell (p^{\ell-n}y)=0$, showing $x\in \capzero{X_n}$. It follows that $\capzero{X_n}\to\capzero{X_n''}$ is surjective for all $n\gg 0$. In particular,~\eqref{diag:thick} has exact rows and the snake lemma shows that since both $\dsd_{m,n}'$ and $\dsd_{m,n}''$ are isomorphisms, the same holds for $\dsd_{m,n}$. Hence, if both $\calX',\calX''$ are in $\bddsys{\DD}$, then also $\calX$ is in $\bddsys{\DD}$. 

Conversely, suppose that $\calX\in\bddsys{\DD}$ and that~$n\geq 0$ is such that $\dsd_{m,n}\colon\capzero{X_m}\to\capzero{X_n}$ is an isomorphism for all~$m\geq n$. Even when the rows in~\eqref{diag:thick} are not exact on the right, the restrictions $\dsd'_{m,n}\colon\capzero{X_m'}\to\capzero{X_n'}$ are injective. It follows that $n\mapsto \gorder{\capzero{X_n'}}$ is a non-increasing sequence of positive integers and thus eventually constant. As an injection of finite groups of the same order is an isomorphism, this shows $\calX'\in\bddsys{\DD}$. In order to show that $\calX''\in\bddsys{\DD}$, consider the following commutative diagram, analogous to~\eqref{diag:thick}:
\begin{equation}\label{diag:thick_dual}\begin{split}
\xymatrix{
\nucbar{X_m'}\ar@{->}[0,1]&\nucbar{X_m}\ar@{->}[0,1]&\nucbar{X_m''}\ar@{->}[0,1]&0\\
\nucbar{X_n'}\ar@{->}[-1,0]_{\dsu_{n,m}'}\ar@{->}[0,1]&\nucbar{X_n}\ar@{->}[0,1]\ar@{->}[-1,0]_{\dsu_{n,m}}&\nucbar{X_n''}\ar@{->}[0,1]\ar@{->}[-1,0]_{\dsu_{n,m}''}&0
}\end{split}\end{equation}
Exactness of the rows follows from Proposition~\ref{prop:univ_capzero}-\ref{point:prop_capzero:functor} and shows that $\dsu_{n,m}''$ is surjective, since $\dsu_{n,m}\colon\nucbar{X_n}\to \nucbar{X_m}$ is an isomorphism by~\ref{point:lemma_eqB:cap=nuc}. Therefore, the orders $\gorder{\nucbar{X_n''}}$ are non-increasing, which implies that $\dsu_{n,m}''$ restricted to $\nucbar{X_n''}$ is actually an isomorphism for $m\geq n\gg 0$, yielding that $\calX''\in\bddsys{\DD}$.
\end{proof}
A consequence of the above proposition is that it is possible to form the quotient category $\dblsys{\DD}/\bddsys{\DD}$, and to speak about double systems that are $\bddsys{\DD}$-isomorphic. For instance, the systems $\calZ$ and $\calZ_{[n_0]}$ from Example~\ref{example:shift} are $\bddsys{\DD}$-isomorphic, two $\bddsys{\DD}$-isomorphisms being the maps $\iota$ and $\phi$ constructed \ibid. To mark the difference between morphisms in $\dblsys{\DD}$ and morphisms in $\quotsys{\DD}$, we denote the latter by dashed arrows, thus writing
\[
f\colon\calX\longdashrightarrow \calX'
\]
for a map $f\in\Hom_{\quotsys{\DD}}(\calX,\calX')$.

Gabriel's theory of quotient categories provides an algebraic framework to study objects, typically occurring in Iwasawa theory, whose sizes is ``eventually equal''. Since this notion is pivotal in our study, we single it out in the following

\begin{definition}\label{def:eventually_prop}
We say that two sequences $\{x_n\}_{n\in\NN}$ and $\{y_n\}_{n\in\NN}$ of positive natural numbers are eventually proportional if there exists $n_0\geq 0$ and a rational number $c\in\QQ$, independent of $n$, such that
\[
x_n=c y_n\qquad\text{ for all }n\geq n_0.
\]
In this case we write $x_n\bddeq y_n$. It is clear that $\bddeq$ is an equivalence relation.
\end{definition}

We isolate the following easy consequence of Lemma~\ref{lemma:eqB}.

\begin{corollary}\label{cor:bddcong_implica_bellezze} 
Let $\calX$ and $\calX'$ be two double systems of \emph{finite} $\ZZ_p[\splitG]$-modules which are $\bddsys{\DD}$-isomorphic, then $\gorder{X_n}\bddeq\gorder{X_n'}$. 
\end{corollary}
\begin{proof}
Since $\calX\bddcong\calX'$, there exists a morphism $f\in\Hom_{\dblsys{\DD}}(\calX,\calX')$ such that both $\kernel(f)$ and $\cokernel(f)$ are in $\bddsys{\DD}$. Lemma~\ref{lemma:eqB} implies that $\gorder{X_n}\bddeq\gorder{X_n'}$. 
\end{proof}

\begin{example}\label{ex:asym_not_B}
A slight modification of Example~\ref{ex:thanks_to_referee} yields the following. Suppose, as \ibid, that $\Gamma$ is filtered by $\Gamma_n=\Gamma^{p^n}$ for all $n\geq 0$ and consider the double system $\calX=(\ZZ/p,0,0)$, where all components $X_n$ coincide with~$\ZZ/p$ and all transition morphisms are the $0$ map. In particular, $\varprojlim \calX=\varinjlim \calX=0$ and therefore $X_n=\capzero{X_n}=\nucbar{X_n}$. It follows that both $\capzero{X_n}$ and $\nucbar{X_n}$ are non-zero, but the transition morphisms $\dsd_{m,m}\colon \capzero{X_m}\to\capzero{X_n}$ and $\dsu_{n,m}\colon\nucbar{X_n}\to\nucbar{X_m}$ are the zero map. Thus, neither of conditions~\ref{point:lemma_eqB:pi} nor~~\ref{point:lemma_eqB:phi} is satisfied, showing that $\calX$ is not an object in $\bddsys{\DD}$, and in particular it is not $\bddsys{\DD}$-isomorphic to the zero system $O=(O_n=0,0,0)$ (since $\bddsys{\DD}$ is thick). Yet, their orders are eventually proportional, showing that being $\bddsys{\DD}$-isomorphic is a stronger notion than having eventually proportional orders.
\end{example}

For the arithmetic applications relevant to this work, it is enough to restrict to a certain subcategory of $\dblsys{\DD}$, which we now define. Observe that every $\ZZ_p$-torsion $\ZZ_p[\splitG]$-module $M$ can be decomposed as $M=M_\mathrm{div}\oplus M_\mathrm{red}$, where $M_\mathrm{div}$ is the maximal $\ZZ_p$-divisible submodule of $M$ and $M_\mathrm{red}$ is the maximal $\ZZ_p$-reduced submodule. This decomposition is in fact one of $\ZZ_p[\splitG]$-modules, because $\splitG$ acts by automorphisms on the $\ZZ_p$-module underlying~$M$.

For a finite group $H$, write $\catMod[{[H]}]{\cofg}$ for the abelian category of $\ZZ_p$-torsion $\ZZ_p[H]$-modules $M$ such that~$M_\mathrm{red}$ is finite and such that $M_\mathrm{div}\cong (\QQ_p/\ZZ_p)^{\lambda}$ for some $\lambda\in\NN$. The subcategory $\catMod[{[H]}]{fin}\subseteq\catMod[{[H]}]{\cofg}$ of finite modules (those with $\lambda=0$) is clearly a thick subcategory, so that the quotient category $\catMod[{[H]}]{\cofg}/\catMod[{[H]}]{fin}$ is well-defined. When $H$ is the trivial group, we drop it from the notation and simply write $\catMod[]{\cofg}$ and $\catMod[]{fin}$.

\begin{definition}
For a torsion $\ZZ_p$-module $M$, we set
\[
\corank M = \dim_{\QQ_p} \left(\Hom_{\ZZ_p}(M,\QQ_p/\ZZ_p) \otimes_{\ZZ_p} \QQ_p\right).
\]
\end{definition}
The function $\corank$ is well-defined on objects of the quotient category $\catMod{\cofg}/\catMod{fin}$, because they coincide with objects of $\catMod{\cofg}$. Moreover, using~\cite[Ch.~III, \S1, Corollaire~1]{Gab62}, one sees that it is multiplicative in short exact sequences because $\corank M=0$ for every $M$ in the thick subcategory $\catMod{fin}$. We will refer at the number $\corank M$ as the \emph{corank of $M$}.

We can finally define the subcategory $\dblsys[\cofg]{\DD}\subseteq \dblsys{\DD}$ mentioned in the Introduction:

\begin{definition}\label{def:dblsys_cofg} 
Let $\dblsys[\cofg]{\DD}$ be the full subcategory of $\dblsys{\DD}$ whose objects $\calX=(X_n,\dsu_{n,m},\dsd_{m,n})_{m\geq n\geq 0}$ satisfy the following conditions:
\begin{enumerate}[label=(\textup{\cofg}~\Alph*), itemindent=\parindent, labelindent=!]
\item Every $X_n$ is a finite group such that $\gindex{\Gamma}{\Gamma_n}X_n=0$ for all $n\gg 0$;\label{point:def_cofg_torsion}
\item $\varinjlim \calX$ belongs to $\catMod{\cofg}$.\label{point:def_cofg_limit} 
\end{enumerate}
Similarly, let $\bddsys[\cofg]{\DD}$ be the full subcategory of $\bddsys{\DD}$ consisting of objects satisfying~\ref{point:def_cofg_torsion}.
\end{definition} 
Observe that $\dblsys[\cofg]{\DD}$ is an abelian category because it is an additive full subcategory of the abelian category $\dblsys{\DD}$ which is closed under kernels and cokernels. Moreover, we claim that $\bddsys[\cofg]{\DD}$ is a thick subcategory of $\dblsys[\cofg]{\DD}$. That every object in $\bddsys[\cofg]{\DD}$ lies in $\dblsys[\cofg]{\DD}$ is obvious, because condition~\ref{point:def_cofg_torsion} holds by definition and condition~\ref{point:def_cofg_limit} follows from condition~\ref{point:lemma_eqB:phi} of Lemma~\ref{lemma:eqB} (which shows that the direct limit is actually finite). This implies, in particular, that objects in $\bddsys[\cofg]{\DD}$ can equivalently be defined as those in $\dblsys[\cofg]{\DD}$ which moreover satisfy either of the conditions of Lemma~\ref{lemma:eqB}. The claim that it is thick is also immediate: given an exact sequence
\[
0\longrightarrow \calX'\longrightarrow\calX\longrightarrow\calX''\longrightarrow 0
\]
in $\dblsys[\cofg]{\DD}$, each of the elements $\calX',\calX,\calX''$ lies in $\bddsys[\cofg]{\DD}$ if and only it lies in $\bddsys{\DD}$, because~\ref{point:def_cofg_torsion} is satisfied automatically by the requirement that the objects in the sequence be in $\dblsys[\cofg]{\DD}$. Therefore $\bddsys[\cofg]{\DD}\subseteq\dblsys[\cofg]{\DD}$ is thick because $\bddsys{\DD}\subseteq\dblsys{\DD}$ is thick, and we can consider the quotient category $\quotsys[\cofg]{\DD}$, which is again abelian.

We now introduce the step-torsion functor $\stepfunctor\colon \catMod{\cofg}\to \dblsys[\cofg]{\DD}$. As in Example~\ref{example:shift}, denote the index $\gindex{\Gamma_n}{\Gamma_m}$ by $\filindex{n}{m}$ and simply write $\filindex[\empty]{}{n}$ for $\gindex{\Gamma}{\Gamma_n}$; similarly, write $g$ for $g(\Gamma)$. Given a $\ZZ_p$-module $M$ and $n\in\mathbb{N}$, denote by $\tors[{\filindex[\empty]{}{n}}]{M}$ the submodule of $M$ consisting of elements annihilated by $\filindex[\empty]{}{n}$. Define $\stepfunctor$ to be the functor that associates to every $\ZZ_p[\splitG]$-module $M\in\catMod{\cofg}$ the $\DD$-double system
\[
\stepfunctor(M)=\bigl(M_n=\tors[{\filindex[\empty]{}{n}}]{M},\dsu_{n,m}=\operatorname{id}_M\colon \tors[{\filindex[\empty]{}{n}}]{M}\hookrightarrow \tors[{\filindex[\empty]{}{m}}]{M},\dsd_{m,n}=\filindex{n}{m}\colon \tors[{\filindex[\empty]{}{m}}]{M}\to \tors[{\filindex[\empty]{}{n}}]{M}\bigr)_{m\geq n\geq 0}
\]
and to every morphism $u\colon M\to M'$ the morphism $\stepfunctor(u)=(u_n=u\vert_{\tors[{\filindex[\empty]{}{n}}]{M}}\colon \tors[{\filindex[\empty]{}{n}}]{M}\to \tors[{\filindex[\empty]{}{n}}]{M'})$. The fact that $\stepfunctor$ takes values in $\dblsys[\cofg]{\DD}$ is immediate: condition~\ref{point:def_cofg_torsion} in Definition~\ref{def:dblsys_cofg} is satisfied by construction, and condition~\ref{point:def_cofg_limit} follows by writing $\varinjlim \stepfunctor(M) = \varinjlim \tors[{\filindex[\empty]{}{n}}]{M} = M$ (since $M$ is $\ZZ_p$-torsion).

It is easy to see that $\stepfunctor$ sends the thick subcategory $\catMod{fin}\subseteq\catMod{\cofg}$ of finite modules to $\bddsys[\cofg]{\DD}$: indeed, given a finite module $M\in\catMod{fin}$, one has $\tors[{\filindex[\empty]{}{n}}]{M}=M$ for all $n\gg 0$, and therefore the transition morphisms $\dsu_{n,m}$ are the identity $\id{M}\colon M\to M$ for all $m\geq n\gg 0$. Proposition~\ref{prop:trivialbdd} shows $\stepfunctor(M)\in\bddsys{\DD}$, and since $\stepfunctor(M)\in\dblsys[\cofg]{\DD}$ by construction, we obtain $\stepfunctor(M)\in\bddsys[\cofg]{\DD}$, as claimed. It follows that $\stepfunctor$ induces a functor between quotient categories, still denoted by the same symbol
\begin{equation}\label{eq:emb_quoz}
\stepfunctor\colon\catMod{\cofg}/\catMod{fin}\longrightarrow
\quotsys[\cofg]{\DD}.
\end{equation}
Related to $\stepfunctor$ is the limit functor
\[
\varinjlim\colon \dblsys[\cofg]{\DD}\longrightarrow \catMod{\cofg}
\]
taking a double system and associating to it the direct limit of the underlying direct system: it takes values in $\catMod{\cofg}$ by definition of the category $\dblsys[\cofg]{\DD}$. By composing it with $\stepfunctor$, we obtain the endofunctor
\[
\limemb=\stepfunctor\circ\varinjlim\colon\dblsys[\cofg]{\DD}\longrightarrow\dblsys[\cofg]{\DD}.
\]
Given any object $\calX$ in $\bddsys[\cofg]{\DD}$, the direct limit $\varinjlim \calX$ is finitely generated and of finite exponent, by Lemma~\ref{lemma:eqB}; it is therefore finite, and the direct limit functor factors  as
\begin{equation}\label{eq:lim_quoz}
\varinjlim\colon \quotsys[\cofg]{\DD}\longrightarrow \catMod{\cofg}/\catMod{fin}
\end{equation}
(see~\cite[Ch.~III, \S1, Corollaire~2]{Gab62}). Combining~\eqref{eq:lim_quoz} with~\eqref{eq:emb_quoz} shows that $\limemb$ extends to an endofunctor, still denoted by the same symbol, 
\[
\limemb\colon\quotsys[\cofg]{\DD}\longrightarrow\quotsys[\cofg]{\DD}.
\]

\begin{corollary}\label{cor:ha_vinto_caputo} 
Let $\calX$ be a double system in $\dblsys[\cofg]{\DD}$, and suppose that $\calX\bddcong\limemb(\calX)$. Then $\gorder{X_n}\bddeq p^{\lambda_\calX n}$, where $\lambda_{\calX}$ is the corank of $\varinjlim\calX$. Moreover, given an exact sequence
\[
0\longdashrightarrow \calX'\longdashrightarrow \calX \longdashrightarrow \calX'' \longdashrightarrow 0
\]
in $\quotsys[\cofg]{\DD}$ of double systems satisfying the above assumption, the relation
\[
\lambda_{\calX}=\lambda_{\calX'}+\lambda_{\calX''}
\]
holds.
\end{corollary}
\begin{proof} 
As observed in Corollary~\ref{cor:bddcong_implica_bellezze}, the existence of a $\bddsys{\DD}$-isomorphism as in the statement implies that $\gorder{X_n}\bddeq \gorder{\tors[{\filindex[\empty]{}{n}}]{(\varinjlim\calX)}}$. On the other hand, the explicit structure of the direct limit of the groups $X_n$ yields
\[
\gorder{\tors[{\filindex[\empty]{}{n}}]{(\varinjlim\calX)}} \bddeq p^{\lambda_{\calX} n}.
\]

Given an exact sequence as in the statement,~\cite[Ch.~III, \S1, Corollaire~1]{Gab62} yields an exact sequence
\[
0\longrightarrow \calY' \longrightarrow \calY \longrightarrow \calY'' \longrightarrow 0
\]
in $\dblsys[\cofg]{\DD}$ such that the diagram 
\begin{equation*}\xymatrix@C=4em{
0\ar@{-->}[0,1]&\calY' \ar@{-->}[0,1]\ar@{-->}[1,0]^{\bddcong}& \calY \ar@{-->}[0,1]\ar@{-->}[1,0]^{\bddcong}& \calY'' \ar@{-->}[0,1]\ar@{-->}[1,0]^{\bddcong}& 0\\
0\ar@{-->}[0,1]&\calX' \ar@{-->}[0,1]& \calX \ar@{-->}[0,1]& \calX'' \ar@{-->}[0,1]& 0
}\end{equation*}
commutes. It yields, again through Corollary~\ref{cor:bddcong_implica_bellezze}, the chain of eventual proportionalities
\[
p^{\lambda_{\calX}n}\bddeq\gorder{X_n}\bddeq \gorder{Y_n}=\gorder{Y_n'}\cdot\gorder{Y_n''}\bddeq\gorder{X_n'}\cdot\gorder{X_n''}\bddeq p^{(\lambda_{\calX'}+\lambda_{\calX''})n}
\]
finishing the proof.
\end{proof}

We conclude this section by analysing the exactness of the functor $\limemb$. 
\begin{proposition}\label{prop:Sexact} 
Given an exact sequence
\[
0\longrightarrow M'\overset{v}{\longrightarrow} M\overset{u}{\longrightarrow} M''\longrightarrow 0
\]
in $\catMod{\cofg}$, the sequence
\begin{equation}\label{eq:exactness_stepfunctor}
0\longdashrightarrow\stepfunctor (M')\overset{\stepfunctor (v)}{\longdashrightarrow}\stepfunctor (M)\overset{\stepfunctor (u)}{\longdashrightarrow} \stepfunctor (M'')\longdashrightarrow 0
\end{equation}
is exact in $\quotsys[\cofg]{\DD}$. In particular, $\limemb$ is an exact functor when regarded as
\[
\limemb\colon \quotsys[\cofg]{\DD}\longrightarrow \quotsys[\cofg]{\DD}.
\]
\end{proposition}
\begin{proof} 
Observe first that $\stepfunctor$ is left exact, since so is the functor $M\mapsto\tors[{\filindex[\empty]{}{n}}]{M}$ for all $n\geq 0$. Therefore we have an exact sequence
\[
0\longrightarrow\stepfunctor (M')\overset{\stepfunctor (v)}{\longrightarrow}\stepfunctor (M)\overset{\stepfunctor (u)}{\longrightarrow} \stepfunctor (M'')\longrightarrow \cokernel{\stepfunctor(u)}\longrightarrow 0.
\] 
We claim that $\cokernel{\stepfunctor(u)}\in\bddsys[\cofg]{\DD}$: this will imply that~\eqref{eq:exactness_stepfunctor} is right exact. 

For $n\in \mathbb{N}$, let $\sch_n\colon\tors[{\filindex[\empty]{}{n}}]{M''}\to M'/{\filindex[\empty]{}{n}}M'$ be the connecting homomorphism of the snake lemma corresponding to multiplication by ${\filindex[\empty]{}{n}}$. The formation of the snake exact sequence being natural, for $m\geq n$ we get commutative diagrams 
\begin{equation}\label{diag:snake_per_caputo}
\begin{array}{lcr}
\xymatrix@C=4em @R4em{
	\tors[{\filindex[\empty]{}{m}}]{M''}\ar@{->}[0,1]^{\sch_m}&M'/{\filindex[\empty]{}{m}}M'\\
	\tors[{\filindex[\empty]{}{n}}]{M''}\ar@{->}[0,1]^{\sch_n}{\ar@{}[-1,0]^(.15){}="a"\ar@{}[-1,0]^(.9){}="b"\ar@{^(->}"a";"b"}&M'/{\filindex[\empty]{}{n}}M'\ar@{->}[-1,0]_{\dsu_{n,m}=\cdot {\filindex{n}{m}}}
}
&\parbox[t][4.7em][c]{6em}{\begin{center}\text{and}\end{center}}&
\xymatrix@C=4em @R4em{
	\tors[{\filindex[\empty]{}{m}}]{M''}\ar@{->}[0,1]^{\sch_m}\ar@{->}[1,0]^{ {\filindex{n}{m}}}&M'/{\filindex[\empty]{}{m}}M'\ar@{->}[1,0]^{\dsd_{m,n}=\operatorname{mod}{\filindex[\empty]{}{n}}M'}\\
	\tors[{\filindex[\empty]{}{n}}]{M''}\ar@{->}[0,1]^{\sch_n}&M'/{\filindex[\empty]{}{n}}M'
}
\end{array}
\end{equation}
The commutativity of~\eqref{diag:snake_per_caputo} together with the exactness of the snake sequence shows that the cokernel $\cokernel{\stepfunctor(u)}$ is a subobject of the double system $\calX=(M'/{\filindex[\empty]{}{n}}M',\dsu_{n,m},\dsd_{m,n})$. Writing $M'=M'_\mathrm{div}\oplus M'_\mathrm{red}$, we obtain
\[
\calX=\bigl(M'/{\filindex[\empty]{}{n}}M'=M'_\mathrm{red}/{\filindex[\empty]{}{n}}M'_\mathrm{red},\dsu_{n,m}=\cdot {\filindex{n}{m}},\dsd_{m,n}\colon M'_\mathrm{red}/{\filindex[\empty]{}{m}}M'_\mathrm{red}\twoheadrightarrow M'_\mathrm{red}/{\filindex[\empty]{}{n}}M'_\mathrm{red}\bigr)_{m\geq n\geq 0}.
\]
Since, by assumption, $M'_\mathrm{red}$ is finite, the quotient $M'_\mathrm{red}/{\filindex[\empty]{}{n}}M'_\mathrm{red}$ coincides with $M'_\mathrm{red}$ for all $n\gg 0$ and the projection $\dsd_{m,n}$ is the identity. Thus, this system satisfies the assumptions of Proposition~\ref{prop:trivialbdd}, showing that $\calX\in\bddsys{\DD}$ and therefore $\cokernel{\stepfunctor(u)}\in\bddsys[\cofg]{\DD}$. This concludes the proof of the exactness on the right of \eqref{eq:exactness_stepfunctor}.

The final assertion concerning $\limemb\colon\quotsys[\cofg]{\DD}\to\quotsys[\cofg]{\DD}$ follows from the previous point, combined with exactness of $\varinjlim\colon\quotsys[\cofg]{\DD}\to\catMod{\cofg}/\catMod{fin}$.
\end{proof}

Given any double system $\calX=(X_n,\dsu_{n,m},\dsd_{m,n})$ in $\dblsys[\cofg]{\DD}$, for each $n\geq 0$ the map $\dsu_n\colon X_n\to\varinjlim \calX$ factors through $\natr_{\calX,n}\colon X_n\to\tors[{{\filindex[\empty]{}{n}}}]{(\varinjlim\calX)}=\limemb(\calX)_n$. We claim that the diagrams \begin{equation}\label{diag:limemb_functorial}
\begin{array}{lcr}
\xymatrix@C=4em@R=3em{
	X_m\ar@{->}[0,1]^{\natr_{\calX,m}}&\limemb(\calX)_m\\
	X_n\ar@{->}[0,1]^{\natr_{\calX,n}}\ar@{->}[-1,0]^{\dsu_{n,m}}&\limemb(\calX)_n{\ar@{}[-1,0]^(.15){}="a"\ar@{}[-1,0]^(.9){}="b"\ar@{^(->}"a";"b"}
}
&\parbox[t][3.75em][c]{7em}{\begin{center}\text{and}\end{center}}&
\xymatrix@C=4em@R=3em{
	X_m\ar@{->}[1,0]_{\dsd_{m,n}}\ar@{->}[0,1]^{\natr_{\calX,m}}&\limemb(\calX)_m\ar@{->}[1,0]^{{\filindex{n}{m}}}\\
	X_n\ar@{->}[0,1]^{\natr_{\calX,n}}&\limemb(\calX)_n
}
\end{array}
\end{equation}
are commutative, and hence determine a morphism $\natr_{\calX}\colon\calX\to\limemb(\calX)$ of double systems. The commutativity of the left square in~\eqref{diag:limemb_functorial} follows from the fact that the direct limit is taken with respect to the morphism $\dsu_{n,m}$. Concerning that on the right, let $x\in X_m$: we need to show that, in the direct limit, the equality
\[
{\filindex{n}{m}}x=\dsd_{m,n}x\in\varinjlim \calX
\]
holds. By definition of direct limit, this follows from the equality ${\filindex{n}{m}}x=\dsu_{m,n}\bigl(\dsd_{m,n}x\bigr)$, which is a consequence of the relation $\dsu_{m,n}\circ\dsd_{n,m}={\filindex{n}{m}}$. The next corollary shows that this construction is functorial:
\begin{corollary}\label{cor:natural_transformation} 
The collection $\natr=(\natr_{\calX})_{\calX\in\dblsys[\cofg]{\DD}}$ defines a natural transformation $\natr\colon \operatorname{id}_{\quotsys[\cofg]{\DD}}\to\limemb$ of endofunctors of the abelian category $\quotsys[\cofg]{\DD}$.
\end{corollary}
\begin{proof} 
Let $f\colon\calX\to\calX'$ be a morphism of double systems. We need to prove the commutativity of
\begin{equation}\label{diag:nat_trans}\begin{split}
\xymatrix@C=4em{
	\calX\ar@{->}[0,1]^{f}\ar@{->}[1,0]_{\natr_{\calX}}&\calX'\ar@{->}[1,0]^{\natr_{\calX'}}\\
	\limemb(\calX)\ar@{->}[0,1]^{\limemb(f)}&\limemb(\calX')
}\end{split}\end{equation}
Fix $n\geq 0$: by construction, $\limemb(\calX)_n=\tors[{\filindex[\empty]{}{n}}]{(\varinjlim\calX)}$ and $\limemb(\calX')_n=\tors[{\filindex[\empty]{}{n}}]{(\varinjlim\calX')}$. Thus, commutativity of~\eqref{diag:nat_trans} translates in the requirement that the following diagram commutes, for all $n\geq 0$:
\begin{equation*}
\xymatrix@C=4em{
	X_n\ar@{->}[0,1]^{f_n}\ar@{->}[1,0]_{\natr_{\calX,n}}&X_n'\ar@{->}[1,0]^{\natr_{\calX',n}}\\
	\tors[{\filindex[\empty]{}{n}}]{\varinjlim\calX'}\ar@{->}[0,1]^{\varinjlim f}&\tors[{\filindex[\empty]{}{n}}]{\varinjlim\calX'}
}
\end{equation*}
This is clear.
\end{proof}

\begin{example}\label{example:Q/Z}
We show that the hypothesis of Corollary~\ref{cor:ha_vinto_caputo} holds for the double systems $\calZ$ and $\calZ_{[n_0]}$ from Example~\ref{example:shift}, where we fix $n_0\geq 1$ in what follows. Since $\varinjlim \calZ=\QQ_p/\ZZ_p$ and each component $(\calZ)_n$ is finite, the system $\calZ$ actually lies in the subcategory $\dblsys[\cofg]{\DD}$. We start by computing $\limemb(\calZ)$. By definition, $\limemb(\calZ)=\stepfunctor\varinjlim\calZ$, and we find
\[
\limemb(\calZ)_n=\tors[{\filindex[\empty]{}{n}}]{\QQ_p/\ZZ_p}=\Bigl(\frac{1}{{\filindex[\empty]{}{n}}}\ZZ\Bigr)/\ZZ\cong \ZZ/{\filindex[\empty]{}{n}}=(\calZ)_n
\]
with the obvious transition morphisms. In particular, the evaluation $\natr_{\calZ}$ of the natural transformation $\natr$ at $\calZ$ is actually an isomorphism in $\dblsys[\cofg]{\DD}$, not only in $\quotsys[\cofg]{\DD}$.

Concerning $\calZ_{[n_0]}$ we have seen in Example~\ref{example:shift} that $\calZ_{[n_0]}\bddcong\calZ$, and clearly $\calZ_{[n_0]}$ is an element of $\dblsys[\cofg]{\DD}$ as well. By first applying $\limemb$ to this isomorphism and  then composing the resulting isomorphism with $\natr_{\calZ}$, we find
\[
\limemb(\calZ_{[n_0]})\bddcong \limemb(\calZ)\bddcong\calZ\bddcong\calZ_{[n_0]}\quad\text { in }\quad\quotsys[\cofg]{\DD}
\]
showing that $\natr_{\calZ_{[n_0]}}$ is an isomorphism as well. Note, though, that this is not an isomorphism before taking the quotient by $\bddsys[\cofg]{\DD}$. Indeed the direct limits $\varinjlim\calZ$ and $\varinjlim\calZ_{[n_0]}$ are both isomorphic to $\QQ_p/\ZZ_p$, showing that
\[
\limemb(\calZ_{[n_0]})\cong \calZ\not\cong \calZ_{[n_0]}\quad\text { in }\quad\dblsys[\cofg]{\DD}.
\]
\end{example}
\section{The double system of the cohomology of units}\label{sec:cohom_units}
\subsection{Arithmetic set-up}\label{subsection:arithmetic_set-up}
In this section, we apply the algebraic results from Section~\ref{sec:algprel} to an arithmetic setting arising in Iwasawa theory. Fix a number field $k$ and a finite normal extension $F/k$, and set $\splitG=\Gal{F}{k}$. Consider a $\ZZ_p$-extension $L_\infty/F$, which is also normal over $k$: set $\DD=\Gal{L_\infty}{k}$ and $\Gamma=\Gal{L_\infty}{F}$, so that there is an exact sequence
\begin{equation}\tag{\ref{eq:seq_gal_groups}}
1 \longrightarrow \Gamma \longrightarrow \DD \longrightarrow  \splitG\longrightarrow 1
\end{equation}
as in Section~\ref{subsec:normic_sys}. Endow $\Gamma$ with the filtration $\Gamma_n=\Gamma^{p^n}$ for all $n\geq 0$, so that, with notation as in~\loccit, $g(\Gamma)=0$ and define the groups $G_n, G_{m,n}$ and $\splitG[n]$ accordingly. Denote by $L_n$ the subfield of $L_\infty$ fixed by $\Gamma_n$, so that $L_\infty=\bigcup L_n$, $F=L_0$ and $[L_m:L_n]=p^{m-n}$.

We refer the reader to \cite{Tat67} and to \cite{Ser67} for the global and local class field theory results that we use. Consider the following commutative diagram of Galois modules with exact rows and columns:
	\begin{equation}\begin{aligned}\label{diag:magic}
	\xymatrix{
	&1\ar@{->}[1,0]&1\ar@{->}[1,0]&1\ar@{->}[1,0]&\\
	1\ar@{->}[0,1]&\unit{L_n}\ar@{->}[0,1]\ar@{->}[1,0]&L_n^\times\ar@{->}[0,1]\ar@{->}[1,0]&\princ{L_n}\ar@{->}[0,1]\ar@{->}[1,0]&1\\
	1\ar@{->}[0,1]&\idunit{L_n}\ar@{->}[0,1]\ar@{->}[1,0]&\idele{L_n}\ar@{->}[0,1]\ar@{->}[1,0]&\ideal{L_n}\ar@{->}[0,1]\ar@{->}[1,0]&1\\
	1\ar@{->}[0,1]&\q{L_n}\ar@{->}[0,1]\ar@{->}[1,0]&\idclg{L_n}\ar@{->}[0,1]\ar@{->}[1,0]&\clsyl{L_n}\ar@{->}[0,1]\ar@{->}[1,0]&1\\
	&1&1&1
	}\end{aligned}\end{equation}
Here we use the following notation for a number field $M$: $M^\times$ denotes the group $(M\setminus \{0\})\otimes\ZZ_p$, $\ideal{M}$ is the group of fractional ideals tensored with $\ZZ_p$, $\princ{M}$ is the group of principal ideals tensored with~$\ZZ_p$, $\clsyl{M}$ is the quotient $\ideal{M}/\princ{M}$ (\ie the $p$-Sylow of the class group of $M$), $\idele{M}$ is the idèle group tensored with $\ZZ_p$, $\idclg{M}$ is the quotient $\idele{M}/M^\times$ (\ie the idèle class group tensored with $\ZZ_p$), $\unit{M}$ is the group of units tensored with~$\ZZ_p$, $\idunit{M}$ is the group of idèles of valuation $0$ at every finite place tensored with $\ZZ_p$ and $\q{M}$ is defined by the exactness of the diagram. Every term in~\eqref{diag:magic} belongs to a $\DD$-normic system, in the sense of Definition~\ref{def:normic_system}, where the transition maps are induced by the inclusion
\[
\arext{L_m}{L_n}\colon L_n\longhookrightarrow L_m
\]
and the (arithmetic) norm
\[
\arnm{L_m}{L_n}\colon L_m^\times \longrightarrow L_n^\times.
\]
The collection of all diagrams~\eqref{diag:magic} for $n\geq 0$ is a commutative diagram in the category $\normsys$. All these normic systems satisfy \ref{cond:inj}, except for the system $(\clsyl{L_n})$. As for \ref{cond:surj}, it only holds for the four normic systems in the upper left square and for the idèle class group.

Together with the global objects introduced above, two local normic systems will be needed. For every prime $\primeid$ of $F$, fix a prime $\Primeid{\infty}$ in $L_\infty$ above $\primeid$. This choice determines a decomposition group $\DD(\primeid)\subseteq \DD$ as well as a decomposition group $\Gamma(\primeid)=\Gamma\cap\DD(\primeid)\subseteq \Gamma$. The prime $\Primeid{\infty}\cap L_n$ of $L_n$ will be denoted $\Primeid{n}$ and its decomposition groups in $G_n$ and in $\splitG[n]$ will be denoted by $G_n(\primeid)$ and $\splitG[n](\primeid)$, respectively. When there are only finitely many primes above $\primeid$ in $L_\infty$ (this will be the only case of interest to us), $\Gamma(\primeid)$ is isomorphic to $\ZZ_p$ and we can endow $(L_{\Primeid{n}}^\times)$ and $(\unit{L_{\Primeid{n}}})$ with a structure of $\DD(\primeid)$-normic system as follows (for a local field $E$ we follow the same convention explained above, so $E^\times:=(E\setminus \{0\})\otimes\ZZ_p$ and $\unit{E}$ is the group of units tensored with $\ZZ_p$). Endow $\Gamma(\primeid)$ with the filtration $\bigl(\Gamma(\primeid)_n=\Gamma(\primeid)\cap\Gamma^{p^n}\bigr)_{n\in\NN}$ so that $\Gamma(\primeid)/\Gamma(\primeid)_n:=G(\primeid)_n=G_n(\primeid)$: this filtration satisfies the condition of \S\ref{subsec:normic_sys} with $g\bigl(\Gamma(\primeid)\bigr)=g$, where $p^g$ is the (finite) number of primes in $L_\infty$ above $\primeid$---in other words, $L_{g}/F$ is the maximal subextension in $L_\infty/F$ where $\primeid$ splits completely. The collections $(L_{\Primeid{n}}^\times)$ and $(\unit{L_{\Primeid{n}}})$, endowed with the natural extension and norm maps, denoted $\arext{\Primeid{m}}{\Primeid{n}}$ and $\arnm{\Primeid{m}}{\Primeid{n}}$, are then $\DD(\primeid)$-normic systems. We denote these normic systems by $\unit{\pallino[\Primeid]}$ and $L_{\pallino[\Primeid]}^\times$, respectively. They satisfy both~\ref{cond:inj} and~\ref{cond:surj}.
Before passing to the analysis of double systems arising from Tate cohomology, we establish one lemma on the residual action on the cohomology of local units. First, we introduce the following notation.
\begin{notation}\label{notation:ram_and_ind}
Given a number field $M_2$ and an extension $M_1/M_2$, denote by $\Ram{M_1}{M_2}$ the set of primes of~$M_2$ which ramify in $M_1/M_2$. Given any prime $\primeid\subseteq\rint{F}$, denote by $\overline{\primeid}$ the prime $\overline{\primeid}=\primeid\cap\rint{k}$ of $k$ and by $p^{\nprimes}$ the number of primes in $L_\infty$ above it. Finally, let $\ramdeg{m}{n}$ (\rsp $\indeg{m}{n}$, $\locdeg{m}{n}=\ramdeg{m}{n}\cdot\indeg{m}{n}$) denote the ramification index (\rsp the inertia degree, the extension degree) of $L_{\Primeid{m}}/L_{\Primeid{n}}$. To lighten the typesetting of some formul\ae\ we also introduce the notation $\logindeg$ to denote $\ln_p(\indeg[\empty]{\infty}{})$ where $\indeg[\empty]{\infty}{}$ is the inertia degree of $\primeid$ in $L_\infty/F$. When the prime $\primeid$ needs to be specified, we write $\logindeg[\primeid]$ instead.
\end{notation}

The following general lemma will be useful to analyse the action of $\splitG$ on the $G_n$-cohomology of the idelic units.

\begin{lemma}\label{lemma:cohomaction}
Let $D$ be a finite group and let $G<H<D$ be subgroups of $D$. For every $H$-module $M$ and every $i\in\ZZ$, there is an isomorphism
\[
\HH{i}{G}{M\otimes_{\ZZ[H]}\ZZ[D]} \cong \HH{i}{G}{M}\otimes_{\ZZ[H]}\ZZ[D]
\]
of $D$-modules.
\end{lemma}
\begin{proof}
Given any $\ZZ[G]$-module $P$, consider the $\ZZ[D]$-homomorphism
\begin{equation*}\begin{split}
\tau\colon\Hom_{\ZZ[G]}(P, M)\otimes_{\ZZ[H]}\ZZ[D] \longrightarrow &\Hom_{\ZZ[G]}(P, M\otimes_{\ZZ[H]}\ZZ[D]) \\
\sum_i (f_i \otimes m_i)\bigr)\longmapsto
&\Bigl(z \longmapsto  \sum_i (f_i(z) \otimes m_i)\Bigr).
\end{split}\end{equation*}
Since $\ZZ[D]$ is $\ZZ[H]$-free, \cite[Chap.I, \S2, n°9, Proposition~10]{Bou61}, shows that $\tau$ is an isomorphism whenever $P$ is projective and finitely presented.
Now let $P^\bullet\to \ZZ$ be resolution of the trivial $\ZZ[G]$-module $\ZZ$ made of finitely presented projective modules. Then
\begin{align}
\HHN{i}{G}{M\otimes_{\ZZ[H]}\ZZ[D]} &= H^i\bigl(\Hom_{\ZZ[G]}(P^\bullet, M\otimes_{\ZZ[H]}\ZZ[D])\bigr) \nonumber\\ 
&\cong H^i\bigl(\Hom_{\ZZ[G]}(P^\bullet,M)\otimes_{\ZZ[H]}\ZZ[D]\bigr) \label{eq:tau}\\
&\cong H^i\bigl(\Hom_{\ZZ[G]}(P^\bullet,M)\bigr)\otimes_{\ZZ[H]}\ZZ[D] \label{eq:univcoeff}\\
&=\HHN{i}{G}{M}\otimes_{\ZZ[H]}\ZZ[D].\nonumber
\end{align}
The isomorphism~\eqref{eq:tau} is induced by $\tau$ and~\eqref{eq:univcoeff} is an isomorphism because $\ZZ[D]$ is $\ZZ[H]$-flat. In particular, we get an isomorphism of $\ZZ[D]$-modules 
\[
\HHN{i}{G}{M\otimes_{\ZZ[H]}\ZZ[D]} \cong \HHN{i}{G}{M}\otimes_{\ZZ[H]}\ZZ[D].
\]
To conclude the proof we argue inductively by backward dimension shifting. We start the induction at $i=1$: we have just proved that for every $\ZZ[H]$-module $N$ there is a $\ZZ[D]$-isomorphism
\[
\HH{1}{G}{N\otimes_{\ZZ[H]}\ZZ[D]}\cong \HH{1}{G}{N}\otimes_{\ZZ[H]}\ZZ[D].
\] 
To go from $i$ to $i-1$, apply the functor $-\otimes_\ZZ M$ to the augmentation map $\ZZ[H]\to \ZZ$. We obtain an exact sequence of $\ZZ[H]$-modules   
\[
0\longrightarrow M' \longrightarrow M\otimes_\ZZ \ZZ[H] \longrightarrow M \longrightarrow 0
\]
for a suitable $M'$, where $M\otimes_\ZZ \ZZ[H]$ is induced, hence cohomologically trivial. We deduce that there is an isomorphism $\HH{i-1}{G}{M} \cong \HH{i}{G}{M'}$ of $\ZZ[H]$-modules,inducing a $\ZZ[D]$-isomorphism
\[
\HH{i-1}{G}{M}\otimes_{\ZZ[H]}\ZZ[D] \cong \HH{i}{G}{M'}\otimes_{\ZZ[H]}\ZZ[D].
\]  
We also have an exact sequence of $\ZZ[D]$-modules
\[
0\longrightarrow M'\otimes_{\ZZ[H]}\ZZ[D] \longrightarrow M\otimes_\ZZ \ZZ[D] \longrightarrow M\otimes_{\ZZ[H]}\ZZ[D] \longrightarrow 0
\]
yielding an isomorphism of $\ZZ[D]$-modules
\[
\HH{i-1}{G}{M\otimes_{\ZZ[H]}\ZZ[D]} \cong \HH{i}{G}{M'\otimes_{\ZZ[H]}\ZZ[D]}.
\] 
Since, by the inductive hypothesis, $\HH{i}{G}{M'\otimes_{\ZZ[H]}\ZZ[D]}\cong \HH{i}{G}{M'}\otimes_{\ZZ[H]}\ZZ[D]$ as $\ZZ[D]$-modules, we obtain a $\ZZ[D]$-isomorphism
\samepage{
\[
\HH{i-1}{G}{M\otimes_{\ZZ[H]}\ZZ[D]}\cong \HH{i-1}{G}{M}\otimes_{\ZZ[H]}\ZZ[D],
\]
concluding the proof.
}
\end{proof}

\begin{lemma}\label{lemma:action_D_on_coom}
Fix $n\geq 0$. The decomposition
\[
\idunit{L_n} = \prod_{\Otherprime[n]\subset \rint{L_n}}\unit{\Otherprime[n]}
\] 
together with Shapiro's lemma induce isomorphisms of $\splitG$-modules
\begin{equation}\label{eq:iso_shapiro}
\HH{i}{G_n}{\idunit{L_n}}\cong\bigoplus_{\overline{\primeid}\in\Ram{L_\infty}{k}}\bigoplus_{\primeid\mid\overline{\primeid}}\HH{i}{G_n(\primeid)}{\unit{\Primeid{n}}}
\end{equation}
for all $i\in \ZZ$. In particular, there are $\splitG$-isomorphisms
\begin{equation}\label{eq:cohom_trivial_D_action}
\HH*{i}{G_n}{\idunit{L_n}}\cong (\ZZ/\ramdeg{n}{0})^{\Ram{L_\infty}{F}}\qquad\text{ for }i=1,2
\end{equation}
where $\splitG$ acts on $(\ZZ/\ramdeg{n}{0})^{\Ram{L_\infty}{F}}$ through its action on the set $\Ram{L_\infty}{F}$.
\end{lemma}
\begin{proof}
We start by analysing the cohomology of local units. Let $\primeid\in \Ram{L_\infty}{F}$ be a prime. Denoting by $v_{\Primeid{n}}$ the $\Primeid{n}$-adic valuation, we have an exact sequence of $\splitG[n](\primeid)$-modules
\begin{equation}\label{seq:exact_of_val}
0\longrightarrow \unit{\Primeid{n}}\longrightarrow L_{\Primeid{n}}^\times\xrightarrow{v_{\Primeid{n}}}\ZZ_p\longrightarrow 0.
\end{equation}
Taking cohomology of the sequence we get, through Hilbert Theorem 90, a $\splitG[n](\primeid)$-equivariant identification 
\begin{equation}\label{eq:h1loc}
\cch{0}_{\Primeid{n}}\colon\ZZ/\ramdeg{n}{0}=\HHN{0}{G_n}{\ZZ}/\im (v_{\Primeid{n}})\cong \HH{1}{G_n(\primeid)}{\unit{\Primeid{n}}}.
\end{equation}
A similar result holds for $i=2$: local class field theory shows that the cup-product with the fundamental class induces a $\splitG[n](\primeid)$-equivariant isomorphism $\ZZ/\locdeg{n}{0}=\HH{0}{G_n(\primeid)}{\ZZ}\cong \HH{2}{G_n(\primeid)}{L_{\Primeid{n}}^\times}$. This restricts to a $\splitG[n](\primeid)$-equivariant isomorphism
\begin{equation}\label{eq:h2loc}
(\cup \fundclass{L_n}{k})\colon \ZZ/\ramdeg{n}{0}\cong\HH{2}{G_n(\primeid)}{\unit{\Primeid{n}}}.
\end{equation} 
Observe that \eqref{eq:h1loc} and \eqref{eq:h2loc} show that $\splitG[n](\primeid)$ acts trivially on $\HH{i}{G_n(\primeid)}{\unit{\Primeid{n}}}$ for $i=1,2$.

We now proceed with the proof of the lemma. For a prime $\primeid$ of $F$ (\rsp a prime $\overline{\primeid}$ of $k$), set 
\[
{\semiloc[\primeid]{L_n}} = \prod_{\Otherprime[n]\mid \primeid} \unit{\Otherprime[n]}\qquad \text{(\rsp }{\semiloc[\overline{\primeid}]{L_n}} = \prod_{\Otherprime[n]\mid \overline{\primeid}} \unit{\Otherprime[n]}\text{)}
\]
where the product is taken over the prime ideals of $L_n$ dividing $\primeid$ (\rsp dividing $\overline{\primeid}$). Then $\idunit{L_n}$ decomposes as 
\[
\idunit{L_n} = \prod_{\overline{\primeid}} \semiloc[\overline{\primeid}]{L_n}
\]
as $\splitG[n]$-modules and, since the local units in unramified extension are cohomologically trivial, this induces an isomorphism of $\splitG[n]$-modules
\[
\HH{i}{G_n}{\idunit{L_n}} = \bigoplus_{\overline{\primeid}\in \Ram{L_\infty}{k}}\HH{i}{G_n}{\semiloc[\overline{\primeid}]{L_n}}.
\]
In particular, to establish~\eqref{eq:iso_shapiro} it suffices to work component-wise, and we focus on $\HH{i}{G_n}{\semiloc[\overline{\primeid}]{L_n}}$ for a fixed prime $\overline{\primeid}\in\Ram{L_\infty}{k}$ from now on. Note that 
\[
\semiloc[\overline{\primeid}]{L_n} = \bigoplus_{\primeid\mid\overline{\primeid}} \semiloc[\primeid]{L_n}
\]
as $G_n$-modules and this gives a $G_n$-isomorphism 
\begin{equation}\label{eq:innominato}
\HH{i}{G_n}{\semiloc[\overline{\primeid}]{L_n}} = \bigoplus_{\primeid\mid\overline{\primeid}} \HH{i}{G_n}{\semiloc[\primeid]{L_n}}.
\end{equation}
To analyse the $\splitG[n]$-action, observe that $\semiloc[\primeid]{L_n}$ is in fact a $H_n(\primeid)$-module where $H_n(\primeid)=\langle G_n, \splitG[n](\primeid)\rangle\subseteq \splitG[n]$ is the subgroup generated by $G_n$ and $\splitG[n](\primeid)$, and
\[
\semiloc[\overline{\primeid}]{L_n} \cong  \semiloc[\primeid]{L_n}\otimes_{\ZZ[H_n(\primeid)]}\ZZ[\splitG[n]]
\]
as $\ZZ[\splitG[n]]$-modules. In particular, Lemma~\ref{lemma:cohomaction} yields an isomorphism of $\splitG[n]$-modules
\[
\HH{i}{G_n}{\semiloc[\overline{\primeid}]{L_n}} \cong \HH{i}{G_n}{\semiloc[\primeid]{L_n}}\otimes_{\ZZ[H_n(\primeid)]}\ZZ[\splitG[n]]=\bigoplus_{\primeid\mid\overline{\primeid}} \HH{i}{G_n}{\semiloc[\primeid]{L_n}}
\] 
that coincides with~\eqref{eq:innominato}.

Now observe that Shapiro's isomorphism factors as 
\[
\HH{i}{G_n}{\semiloc[\primeid]{L_n}}=\HH{i}{G_n}{\prod_{\Otherprime\mid\primeid}\unit{\Otherprime}}
\xrightarrow{\res_{G_n(\primeid)}}
\HH{i}{G_n(\primeid)}{\prod_{\Otherprime\mid\primeid}\unit{\Otherprime}}={\bigoplus}_{\Otherprime\mid\primeid}\HH{i}{G_n(\primeid)}{\unit{\Otherprime}}
\xrightarrow{\pr_{\Primeid{n}}}
\HH{i}{G_n(\primeid)}{\unit{\Primeid{n}}}
\] 
(see~\cite[\S7.2]{Tat67}) and thus it is $\splitG[n](\primeid)$-equivariant. Since $\splitG[n](\primeid)$ acts trivially on $\HH{i}{G_n(\primeid)}{\unit{\Primeid{n}}}$ in light of the discussion in the first part of the proof, the same holds for $\HH{i}{G_n}{\semiloc[\primeid]{L_n}}$ and therefore the whole $H_n(\primeid)$ acts trivially on $\HH{i}{G_n}{\semiloc[\primeid]{L_n}}$. Hence we have
\[
\HH{i}{G_n}{\semiloc[\primeid]{L_n}}\otimes_{\ZZ[H_n(\primeid)]}\ZZ[\splitG[n]] = \HH{i}{G_n}{\semiloc[\primeid]{L_n}}\otimes_{\ZZ}\ZZ[\splitG/\splitG(\primeid)] \cong \HH{i}{G_n(\primeid)}{\unit{\Primeid{n}}}\otimes_{\ZZ}\ZZ[\splitG/\splitG(\primeid)]
\]
as $\splitG[n]$-modules, because $\splitG[n]/H_n(\primeid)\cong \splitG/\splitG(\primeid)$. This concludes the proof of the first assertion of the lemma.

The existence of the isomorphisms~\eqref{eq:cohom_trivial_D_action} follows from \eqref{eq:h1loc} and \eqref{eq:h2loc}.
\end{proof}

\subsection{Double systems of Tate cohomology}\label{subsec:cohom_units}
Given any $\DD$-normic system $\calB=(B_{L_n})$ as in \S\ref{subsec:normic_sys}, we equivalently denote its direct (\rsp inverse) limit as $\varinjlim \calB$ or $\varinjlim B_{L_n}$ or $B_{L_\infty}$ (\rsp $\varprojlim \calB$ or $\varprojlim B_{L_n}$), according to convenience. If $\calB$ satisfies~\ref{cond:inj} (\rsp \ref{cond:inj} and \ref{cond:surj}), then for all odd $i\in\ZZ$ (\rsp all $i\in\ZZ$) we can consider the collection
\[
\sysHH{i}{\calB}=\bigl(\HH{i}{G_n}{B_{L_n}},\HHm{i}{\nsu_{n,m}},\HHm{i}{\nsd_{m,n}}\bigl)_{m\geq n\geq 0}
\]
which is $\DD$-a double system in the language of \S\ref{subsec:double_sys} by Proposition \ref{prop:cohomds}. Since any cohomology group of $G_n$ is annihilated by $p^n$, for every $\DD$-normic system $\calB$ satisfying~\ref{cond:inj} (\rsp both~\ref{cond:inj} and~\ref{cond:surj}), condition~\ref{point:def_cofg_torsion} is satisfied for $\sysHH{i}{\calB}$ for $i$ odd (\rsp for all $i$). Moreover there are canonical isomorphisms
\[
\HHN{i}{\Gamma}{B_{L_\infty}}\cong \varinjlim\HH{i}{G_n}{B_{L_n}}\qquad \text{ for }i\geq 1,
\]
where the limit is taken with respect to the maps $\infl\circ \nsu^*$. In particular, when $i\geq 1$, condition~\ref{point:def_cofg_limit} for $\sysHH{i}{\calB}$ is equivalent to the condition that $\HHN{i}{\Gamma}{B_{L_\infty}}$ belong to $\catMod{\cofg}$, because the direct systems $\bigl(\HH{i}{G_m}{B_{L_m}},\HHm{i}{\nsu}\bigr)$ and $\bigl(\HH{i}{G_m}{B_{L_m}},\infl\circ \nsu^\ast\bigr)$ coincide by Lemma~\ref{lemma:HHm(j)_and_inf}. Hence, when this is the case, the double system $\sysHH{i}{\calB}$ actually lies in $\dblsys[\cofg]{\DD}$ and can be regarded as an object of the quotient category $\quotsys[\cofg]{\DD}$. When considering a local normic systems $\calB_{\pallino[\Primeid]}\in\{\unit{\pallino[\Primeid]},L_{\pallino[\Primeid]}^\times\}$, the cohomology we consider is relative to the decomposition groups $G_n(\primeid)$. To stress this difference, we denote by a subscript $\primeid$ the corresponding $\DD(\primeid)$-double system
\[
\sysHH[\primeid]{i}{\calB_{\pallino[\Primeid]}}=\bigl(\HH{i}{G_n(\primeid)}{B_{\Primeid{n}}},\arext{\Primeid{m}}{\Primeid{n}},\arnm{\Primeid{m}}{\Primeid{n}})_{m\geq n\geq 0}.
\]

The starting point of our study is the exact sequence of double systems
\begin{equation}\label{eq:long_cohom_DS}
\sysHH{1}{\unit{\pallino}}\longrightarrow\sysHH{1}{\idunit{\pallino}}\longrightarrow\sysHH{1}{\q{\pallino}}\longrightarrow\sysHH{2}{\unit{\pallino}}\longrightarrow\sysHH{2}{\idunit{\pallino}}
\end{equation}
whose existence follows from Propositions~\ref{prop:long_cohom_up} and~\ref{prop:long_cohom_down}. Let $\calB\in\{\unit{\pallino},\idunit{\pallino},\q{\pallino}\}$ be any of the normic systems occurring in~\eqref{eq:long_cohom_DS}. It is well-known that $\HHN{i}{\Gamma}{\varinjlim \calB}$ belongs to $\catMod{\cofg}$ when $i=1,2$ (see, for instance, \cite[\S 2]{Iwa83}), so~\eqref{eq:long_cohom_DS} is an exact sequence in $\dblsys[\cofg]{\DD}$. 

We first analyse the $\DD$-double systems $\sysHH{i}{\idunit{\pallino}}$, and we single out their relationship with the $\DD(\primeid)$-double systems $\sysHH[\primeid]{i}{\unit{\pallino[\Primeid]}}$ in the following proposition.

\begin{proposition}\label{prop:shapiro}
For each prime $\overline{\primeid}\in\Ram{L_\infty}{k}$, the sum $\bigoplus_{\primeid\mid\overline{\primeid}}\sysHH[\primeid]{i}{\unit{\pallino[\Primeid]}}$ is a $\DD$-double system, with the $\splitG$-action discussed in Lemma~\ref{lemma:action_D_on_coom}. Moreover the isomorphisms~\eqref{eq:iso_shapiro} induce isomorphisms of $\DD$-double systems 
\begin{equation}\label{eq:idunits_sum_locunits}
\sysHH{i}{\idunit{\pallino}}\cong\bigoplus_{\overline{\primeid}\in \Ram{L_\infty}{k}}\bigoplus_{\primeid\mid\overline{\primeid}}\sysHH[\primeid]{i}{\unit{\pallino[\Primeid]}} \qquad\text{ for all }i\in\ZZ.
\end{equation}
\end{proposition}
\begin{proof}
To prove the first assertion we need to show the relation
\begin{equation}\label{eq:locunitsds}
\HHm{i}{\nsu_{n,m}}\circ\HHm{i}{\nsd_{m,n}}=\HHm{i}{\nsd_{m,n}}\circ\HHm{i}{\nsu_{n,m}}=\gindex{\Gamma_n}{\Gamma_m}\qquad\text{ for all }m\geq n\geq 0
\end{equation}
for the transition morphisms of the system $\bigoplus_{\primeid\mid\overline{\primeid}}\HH{i}{G_n(\primeid)}{\unit{\Primeid{n}}}$. Let $\indextr[\primeid]$ be the first index above which any of the local $\ZZ_p$-extensions $L_{\Primeid{\infty}}/F_{\primeid}$ (for any  $\primeid\mid\overline{\primeid}$) is totally ramified. For every natural number $n\geq \indextr[\primeid]$, we have $\gindex{\Gamma_n}{\Gamma_m}=p^{m-n}=\gindex{\Gamma(\primeid)_n}{\Gamma(\primeid)_m}$, and \eqref{eq:locunitsds} is satisfied since every $\sysHH[\primeid]{i}{\unit{\pallino[\Primeid]}}$ is a $\DD(\primeid)$-double system. When $n\leq \indextr[\primeid]$, the group $\bigoplus_{\primeid\mid\overline{\primeid}}\HH{i}{G_n(\primeid)}{\unit{\Primeid{n}}}$ vanishes because the cohomology of local units in unramified extensions is trivial. Therefore in this case both $\HHm{i}{\nsu_{n,m}}\circ\HHm{i}{\nsd_{m,n}}$ and $\HHm{i}{\nsd_{m,n}}\circ\HHm{i}{\nsu_{n,m}}$ are trivial endomorphisms and \eqref{eq:locunitsds} is equivalent to saying that $\gindex{\Gamma_n}{\Gamma_m}$ annihilates $\bigoplus_{\primeid\mid\overline{\primeid}}\HH{i}{G_n(\primeid)}{\unit{\Primeid{n}}}$ and $\bigoplus_{\primeid\mid\overline{\primeid}}\HH{i}{G_m(\primeid)}{\unit{\Primeid{m}}}$. This is obvious for $\bigoplus_{\primeid\mid\overline{\primeid}}\HH{i}{G_n(\primeid)}{\unit{\Primeid{n}}}=0$. As for the groups at level $m$, observe that each group $\HH{i}{G_m(\primeid)}{\unit{\Primeid{m}}}$ is annihilated by $\gorder{G_{m,\indextr[\primeid]}(\primeid)}=\gindex{\Gamma(\primeid)_{\indextr[\primeid]}}{\Gamma(\primeid)_m}=\ramdeg{m}{\indextr[\primeid]}=\ramdeg{m}{0}$: we have shown this in Lemma~\ref{lemma:action_D_on_coom} for degrees $1$ and $2$, and periodicity of Tate cohomology yields the same result for all~$i\in\ZZ$. Thus, since $\gindex{\Gamma(\primeid)_{\indextr[\primeid]}}{\Gamma(\primeid)_m}$ divides $\gindex{\Gamma_n}{\Gamma_m}$ because $n\leq \indextr[\primeid]$, it follows that $\gindex{\Gamma_n}{\Gamma_m}$ also acts as $0$ on $\bigoplus_{\primeid\mid\overline{\primeid}}\HH{i}{G_m(\primeid)}{\unit{\Primeid{m}}}$. 

Having established that the right-hand side of~\eqref{eq:idunits_sum_locunits} is a $\DD$-double system, the fact that~\eqref{eq:idunits_sum_locunits} is an isomorphism will be an immediate consequence of Lemma~\ref{lemma:action_D_on_coom} once we prove that the isomorphisms~\eqref{eq:iso_shapiro} are compatible with the ascending and descending morphisms, hence defining isomorphisms of double systems. Since both the ascending and descending morphisms are relative to the extension $L_\infty/F$ (rather than $L_\infty/k$), we can and do fix a prime $\primeid\subseteq\rint{F}$ and check the compatibility for the relevant component. Given the explicit definition of Shapiro's lemma in Tate cohomology, the required compatibilities correspond to the commutativity, for all~$i\in\ZZ$, of the diagrams\footnote{We adopt the typesetting convention of denoting by $\Otherprime$ the generic prime dividing $\primeid$ in $L_m$ and denoting by $\otherprime$ the generic prime dividing $\primeid$ in $L_n$. We adopt the notation $\arnm{\Otherprime}{\otherprime}$ and $\arext{\Otherprime}{\otherprime}$ to denote the corresponding norm and extension maps.}
\begin{equation}\begin{split}\label{diag:shapiro_up}
\xymatrix@C4.5em@R3em{
\HH{i}{G_m}{\displaystyle{\prod_{\Otherprime\mid\primeid}\unit{\Otherprime}}}\ar@{->}[0,1]^(.35){\res_{G_m(\primeid)}}&\HH{i}{G_m(\primeid)}{\displaystyle{\prod_{\Otherprime\mid\primeid}\unit{\Otherprime}}}=\prod_{\Otherprime\mid\primeid}\HH{i}{G_m(\primeid)}{\unit{\Otherprime}}\ar@{->}[0,1]^(.65){\pr_{\Primeid{m}}}&\HH{i}{G_m(\primeid)}{\unit{\Primeid{m}}}\\
\HH{i}{G_n}{\displaystyle{\prod_{\otherprime\mid\primeid}\unit{\otherprime}}}\ar@{->}[0,1]^(.35){\res_{G_n(\primeid)}}\ar@{->}[-1,0]_{\HHm{i}{\arext{L_m}{L_n}}}&\HH{i}{G_n(\primeid)}{\displaystyle{\prod_{\otherprime\mid\primeid}\unit{\otherprime}}}=\prod_{\otherprime\mid\primeid}\HH{i}{G_n(\primeid)}{\unit{\otherprime}}\ar@{->}[0,1]^(.65){\pr_{\Primeid{n}}}&\HH{i}{G_n(\primeid)}{\unit{\Primeid{n}}}\ar@{->}[-1,0]_{\HHm{i}{\arext{\Primeid{m}}{\Primeid{n}}}}
}\end{split}\end{equation}
and
\begin{equation}\begin{split}\label{diag:shapiro_down}
\xymatrix@C4.5em@R3em{
\HH{i}{G_m}{\displaystyle{\prod_{\Otherprime\mid\primeid}\unit{\Otherprime}}}\ar@{->}[0,1]^(.35){\res_{G_m(\primeid)}}\ar@{->}[1,0]^{\HHm{i}{\arnm{L_m}{L_n}}}&\HH{i}{G_m(\primeid)}{\displaystyle{\prod_{\Otherprime\mid\primeid}\unit{\Otherprime}}}=\prod_{\Otherprime\mid\primeid}\HH{i}{G_m(\primeid)}{\unit{\Otherprime}}\ar@{->}[0,1]^(.65){\pr_{\Primeid{m}}}
&\HH{i}{G_m(\primeid)}{\unit{\Primeid{m}}}\ar@{->}[1,0]^{\HHm{i}{\arnm{\Primeid{m}}{\Primeid{n}}}}\\
\HH{i}{G_n}{\displaystyle{\prod_{\otherprime\mid\primeid}\unit{\otherprime}}}\ar@{->}[0,1]^(.35){\res_{G_n(\primeid)}}&\HH{i}{G_n(\primeid)}{\displaystyle{\prod_{\otherprime\mid\primeid}\unit{\otherprime}}}=\prod_{\otherprime\mid\primeid}\HH{i}{G_n(\primeid)}{\unit{\otherprime}}\ar@{->}[0,1]^(.65){\pr_{\Primeid{n}}}&\HH{i}{G_n(\primeid)}{\unit{\Primeid{n}}}
}\end{split}
\end{equation}
All terms in the bottom lines of~\eqref{diag:shapiro_up} and~\eqref{diag:shapiro_down} vanish if $n\leq \nprimes$, and there is nothing to check: hence we assume, from now on, that $m\geq n\geq \nprimes$. Observe that Shapiro's isomorphism is compatible with the cup-products defining $\HHm{i}{\nsu}$ and $\HHm{i}{\nsd}$ for arbitrary $i$. More precisely
\[
\res_{G_n(\primeid)}(x \cup \varkappa_n) = \res_{G_n(\primeid)}(x) \cup \res_{G_n(\primeid)}(\varkappa_n)
\]
and $\res_{G_n(\primeid)}(\varkappa_n)$ corresponds to the element of $\Hom(G_n(\primeid),\QQ_p/\ZZ_p)$ defined by 
\begin{equation}\label{eq:lcdn1}
\res_{G_n(\primeid)}(\chi_n)\colon\topgen^{p^{\nprimes}}\longmapsto p^{\nprimes}/\gindex{\Gamma}{\Gamma_n}=\gindex{\Gamma}{\Gamma(\primeid)}/\gindex{\Gamma}{\Gamma_n}.
\end{equation}
On the other hand, by definition of $\HHm{i}{\nsu}$ and $\HHm{i}{\nsd}$ for the local system, they coincide with cup product with the class 
$\varkappa_n^{\primeid}\in\HH{2}{G_n(\primeid)}{\ZZ_p}$ that corresponds to the element of $\Hom(G_n(\primeid),\QQ_p/\ZZ_p)$ defined by 
\begin{equation}\label{eq:lcdn2}
\chi_n^\primeid\colon\topgen^{p^{\nprimes}}\longmapsto 1/\gindex{\Gamma(\primeid)}{\Gamma_n(\primeid)}.
\end{equation}
Since $n\geq \nprimes$, we have $\Gamma_n(\primeid) = \Gamma_n$ so that \eqref{eq:lcdn1} and \eqref{eq:lcdn2} agree. Therefore it is enough to check commutativity of both~\eqref{diag:shapiro_up} and of~\eqref{diag:shapiro_down} for $i=-1$ and for $i=0$. 

Starting with $i=-1$, let $\xi\in \HH{i}{G_n}{\prod_{\otherprime\mid\primeid}\unit{\otherprime}}$ and $\Xi\in \HH{i}{G_m}{\prod_{\Otherprime\mid\primeid}\unit{\Otherprime}}$. By~\cite[Proposition~7]{AtiWal67}, restriction in degree $-1$ is represented by multiplication by the algebraic norm, hence commutativity of the
diagrams~\eqref{diag:shapiro_up} and~\eqref{diag:shapiro_down} are equivalent to the equalities 
\begin{align}
\HHm{-1}{\arext{\Primeid{m}}{\Primeid{n}}}\pr_{\Primeid{n}}\bigl(\algnm{G_n/G_n(\primeid)}\xi\bigr)&=\pr_{\Primeid{m}}\Bigl(\algnm{G_m/G_m(\primeid)}\bigl(\HHm{-1}{\arext{L_m}{L_n}} \xi\bigr)\Bigr)\label{eq:shapiro_up_mem_one}
\intertext{and}
\HHm{-1}{\arnm{\Primeid{m}}{\Primeid{n}}}\pr_{\Primeid{m}}\bigl(\algnm{G_m/G_m(\primeid)}\Xi\bigr)&=\pr_{\Primeid{n}}\Bigl(\algnm{G_n/G_n(\primeid)}\bigl(\HHm{-1}{\arnm{L_m}{L_n}} \Xi\bigr)\Bigr).\label{eq:shapiro_down_mem_one}
\end{align}
To check~\eqref{eq:shapiro_up_mem_one}, suppose that $\xi$ can be represented by $(\xi_\otherprime)\in \prod_{\otherprime\mid\primeid}\unit{\otherprime}$. Thanks to the assumption $n\geq \nprimes$, the restriction of Galois automorphisms induces a bijection between the sets $G_m/G_m(\primeid)$ and $G_n/G_n(\primeid)$: denote this bijection by a bar, writing $\overline{\sigma}$ for the generic element of $G_m/G_m(\primeid)$ that is the restriction of $\sigma\in G_n/G_n(\primeid)$. Then, with a slight abuse of notation consisting in identifying $\xi$ with the tuple of $\xi_{\otherprime}$'s,
\begin{align*}
\HHm{-1}{\arext{\Primeid{m}}{\Primeid{n}}}\pr_{\Primeid{n}}\bigl(\algnm{G_n/G_n(\primeid)}\xi\bigr) & = \HHm{-1}{\arext{\Primeid{m}}{\Primeid{n}}}\pr_{\Primeid{n}}\left(\Bigl(\;\prod_{\overline{\sigma}\in G_n/G_n(\primeid)}\overline{\sigma}(\xi_{\overline{\sigma}^{-1}\otherprime})\Bigr)_{\otherprime\mid\primeid}\right)\\
& = \HHm{-1}{\arext{\Primeid{m}}{\Primeid{n}}}\prod_{\overline{\sigma}\in G_n/G_n(\primeid)}\overline{\sigma}(\xi_{\overline{\sigma}^{-1}\Primeid{n}})\\
& = \prod_{\overline{\sigma}\in G_n/G_n(\primeid)}\arext{\Primeid{m}}{\Primeid{n}}\overline{\sigma}(\xi_{\overline{\sigma}^{-1}\Primeid{n}})\\
& = \prod_{\sigma\in G_m/G_m(\primeid)}\sigma(\xi_{\sigma^{-1}\Primeid{n}})\\
& = \pr_{\Primeid{m}}\biggl(\Bigl(\prod_{\sigma\in G_m/G_m(\primeid)}\sigma(\arext{\sigma^{-1}\Otherprime}{\sigma^{-1}\otherprime}\xi_{\sigma^{-1}\otherprime})\Bigr)_{\Otherprime\mid\primeid}\biggr)\\
& = \pr_{\Primeid{m}}\Bigl(\algnm{G_m/G_m(\primeid)}\bigl(\HHm{-1}{\arext{L_m}{L_n}} \xi\bigr)\Bigr)
\end{align*}
establishing~\eqref{eq:shapiro_up_mem_one}. 

The computation for~\eqref{eq:shapiro_down_mem_one} is analogous, by replacing the equality $\overline{\sigma}\circ\arext{\sigma^{-1}\Primeid{m}}{\sigma^{-1}\Primeid{n}}=\arext{\Primeid{m}}{\Primeid{n}}\circ \sigma$ with $\overline{\sigma}\circ\arnm{\sigma^{-1}\Primeid{m}}{\sigma^{-1}\Primeid{n}}=\arnm{\Primeid{m}}{\Primeid{n}}\circ \sigma$ for all $\sigma\in G_m/G_m(\primeid)$ corresponding to $\overline{\sigma}\in G_n/G_n(\primeid)$.

The case $i=0$ can be handled similarly. This time, restriction reduces to projecting a class modulo $\algnm{G_n}$ (\rsp modulo $\algnm{G_m}$) to the class modulo $\algnm{G_n(\primeid)}$ (\rsp modulo $\algnm{G_m(\primeid)}$). Given an element $b\mod{\algnm{G_n}}\in \HH{0}{G_n}{\prod_{\otherprime\mid\primeid}\unit{\otherprime}}$ (\rsp an element $B \mod{\algnm{G_m}}\in \HH{0}{G_m}{\prod_{\Otherprime\mid\primeid}\unit{\Otherprime}}$, the commutativity of diagrams~\eqref{diag:shapiro_up} and~\eqref{diag:shapiro_down} can be written as
\begin{align}
\HHm{0}{\arext{\Primeid{m}}{\Primeid{n}}}\pr_{\Primeid{n}}\bigl(b\mod{\algnm{G_n(\primeid)}}\bigr)&=\pr_{\Primeid{m}}\bigl(\HHm{0}{\arext{L_m}{L_n}} (b\mod{\algnm{G_m(\primeid)}})\bigr)\label{eq:shapiro_up_mem_zero}
\intertext{and}
\HHm{0}{\arnm{\Primeid{m}}{\Primeid{n}}}\pr_{\Primeid{m}}\bigl(B\mod{\algnm{G_m(\primeid)}}\bigr)&=\pr_{\Primeid{n}}\bigl(\HHm{0}{\arnm{L_m}{L_n}}( B\mod{\algnm{G_m(\primeid)}})\bigr).\label{eq:shapiro_down_mem_zero}
\end{align}
Thanks to the assumption $n\geq \nprimes$, the sets $\{\Otherprime\mid\primeid\}$ and $\{\otherprime\mid\primeid\}$ are in bijection; moreover, the index $\gindex{\Gamma_n}{\Gamma_m}$ equals $p^{m-n}$. Thus, writing $b=(b_\otherprime)_{\otherprime\mid\primeid}$, the relation~\eqref{eq:shapiro_up_mem_zero} reduces to the chain of equalities
\begin{align*}
\HHm{0}{\arext{\Primeid{m}}{\Primeid{n}}}\pr_{\Primeid{n}}\bigl(b\mod{\algnm{G_n(\primeid)}}\bigr)&=\HHm{0}{\arext{\Primeid{m}}{\Primeid{n}}}\bigl(b_{\Primeid{n}}\mod{\algnm{G_n(\primeid)}}\bigr)\\
&=\bigl(\arext{\Primeid{m}}{\Primeid{n}}b_{\Primeid{n}}^{p^{m-n}}\mod{\algnm{G_m(\primeid)}}\bigr)\\
&=\pr_{\Primeid{m}}\bigl(\arext{\Otherprime}{\otherprime}b_{\otherprime}^{p^{m-n}}\mod{\algnm{G_m(\primeid)}}\bigr)_{\Otherprime\mid\primeid}\\
&=\pr_{\Primeid{m}}\bigl(\HHm{0}{\arext{L_m}{L_n}} (b\mod{\algnm{G_m(\primeid)}})\bigr).
\end{align*}
Equation~\eqref{eq:shapiro_down_mem_zero} can be checked similarly.
\end{proof}

We will use the above proposition to obtain an explicit description of the double system of the cohomology of local units. First, we need a lemma.

\begin{lemma}\label{lemma:fundclasscommute}
Let $\primeid\in\Ram{L_\infty}{F}$ be a prime of $F$ and set $\loclev{}=F_\primeid$. For every $s\in \NN$, set $\loclev{s}=L_{\Primeid{s}}$ and let $\fundclass{\loclev{s}}{\loclev{\empty}}\in \HHN{2}{G_{s}(\primeid)}{\loclev{s}^\times}$ be the fundamental class of $\loclev{s}/\loclev{}$. Then, for $m\geq n\geq 0$, the diagrams
\[\begin{array}{lcr}
\xymatrix@C=4pc{
	\HH{0}{G_m(\primeid)}{\ZZ_p}\ar@{->}[0,1]^{(\cup\fundclass{\loclev{m}}{\loclev{}})}
	&\HH{2}{G_m(\primeid)}{\loclev{m}^\times}\\
	\HH{0}{G_n(\primeid)}{\ZZ_p}\ar@{->}[0,1]^{(\cup \fundclass{\loclev{n}}{\loclev{}}) }\ar@{->}[-1,0]_{\cdot\locdeg{m}{n}}
	&\HH{2}{G_n(\primeid)}{\loclev{n}^\times}\ar@{->}[-1,0]_{\HHm{2}{\arext{\loclev{m}}{\loclev{n}}}}
}
&\parbox[t][3.6em][c]{3.3em}{\begin{center}\text{and}\end{center}}&
\xymatrix@C=4pc{
	\HH{0}{G_m(\primeid)}{\ZZ_p}\ar@{->}[0,1]^{(\cup\fundclass{\loclev{m}}{\loclev{}})}
	&\HH{2}{G_m(\primeid)}{\loclev{m}^\times}\\
	\HH{0}{G_n(\primeid)}{\ZZ_p}\ar@{->}[0,1]^{(\cup \fundclass{\loclev{n}}{\loclev{}}) }\ar@{<-}[-1,0]
	&\HH{2}{G_n(\primeid)}{\loclev{n}^\times}\ar@{<-}[-1,0]_{\HHm{2}{\arnm{\loclev{m}}{\loclev{n}}}}
}
\end{array}\]
commute, where the map $\HH{0}{G_m(\primeid)}{\ZZ_p}\to\HH{0}{G_n(\primeid)}{\ZZ_p}$ is the canonical projection and $\locdeg{m}{n}=[E_m:E_n]$ is the degree of the local extension.
\end{lemma}
\begin{proof}
For the first diagram, the proof is similar to that of Lemma~\ref{lemma:HHm(j)_and_inf} for $i=2$, replacing $\varkappa$ by the fundamental class and using that
\[
\HHm{2}{\arext{\loclev{m}}{\loclev{n}}}(\fundclass{\loclev{n}}{\loclev{}})=\infl^2\circ \arext{\loclev{m}}{\loclev{n}}^*(\fundclass{\loclev{n}}{\loclev{}}) = \locdeg{m}{n}\fundclass{\loclev{m}}{\loclev{}}
\]	
(see~\cite[Ch.~XI~\S3]{Ser62}, note that Serre's $\operatorname{Inf}$ is our $\infl\circ\arext{\loclev{m}}{\loclev{n}}^*$).

To analyse the second diagram, let $\chi_s^\primeid\in\HH{1}{G_s(\primeid)}{\QQ_p/\ZZ_p}$ be the element introduced in~\eqref{eq:lcdn2} and set $\varkappa^\primeid_s=\cch{1}\chi_s^\primeid\in\HH{2}{G_s(\primeid)}{\ZZ_p}$, for $s\geq 0$. As discussed in~\S\ref{subsubsec:arbitrarydegree}, the relation
\begin{equation}\label{eq:chiproj}
\chi_n^{\primeid}(\pr_{m,n}\sigma) =  \chi^\primeid_m (\sigma)^{\locdeg{m}{n}}
\end{equation}
holds for all $\sigma\in G_m(\primeid)^\mathrm{ab}\otimes \ZZ_p=G_m(\primeid)$
where $\pr_{m,n}\colon G_m(\primeid)\to G_n(\primeid)$ is the surjective map induced by the restriction to $\loclev{n}$. Given $a\in\ZZ_p$, set $\overline{y} = \tateiso[\primeid]{m}{-1}\bigl((a \mod{\gorder{G_m(\primeid)}}) \cup \fundclass{\loclev{m}}{k}\bigr)\in \HH{0}{G_m(\primeid)}{\loclev{m}^\times}$ so that
\begin{equation}\label{eq:ymdef}
\overline{y}\cup \varkappa^\primeid_m = (a \mod{\gorder{G_m(\primeid)}}) \cup \fundclass{\loclev{m}}{\loclev{}}.
\end{equation}
If $y\in(\loclev{m}^\times)^{G_m(\primeid)}=\loclev{}^\times$ represents $\overline{y}$, then it also represents $\HHm{0}{\arnm{\loclev{m}}{\loclev{n}}}(\overline{y})\in\HH{0}{G_m(\primeid)}{\loclev{m}^\times}$. Denote by $(\;\cdot\;,\loclev{m}/\loclev{})\colon \loclev{}^\times \to G_m(\primeid)^{\mathrm{ab}}\otimes\ZZ_p$ the reciprocity map of local class field theory. Then, with notation as in~\cite[Ch.~XI]{Ser62},
\begin{align*}
\invJP\bigl[\HHm{2}{\arnm{\loclev{m}}{\loclev{n}}}\bigl((a \mod{\gorder{G_m(\primeid)}}) \cup
\fundclass{\loclev{m}}{\loclev{\empty}}
\bigr)\bigr] &= \invJP\bigl[\tateiso[\primeid]{n}{}\circ\HHm{0}{\arnm{\loclev{m}}{\loclev{n}}}(\overline{y}_m)\bigr]\\
&=\invJP\bigl[\HHm{0}{\arnm{\loclev{m}}{\loclev{n}}}(\overline{y})\cup\cch{1}\chi^{\primeid}_n\bigr]\\
\text{(by \cite[Ch.~XI, Proposition 2]{Ser62})}&=\langle\chi^{\primeid}_n,\bigl( y,\loclev{n}/{\loclev{}}\bigr)\rangle\\ 
\text{(by \cite[p.178, (4)]{Ser62})}&=\langle\chi^{\primeid}_n,\pr_{m,n}\bigl((y,\loclev{m}/\loclev{})\bigr)\rangle\\
\text{(by~\eqref{eq:chiproj})}&=\langle\chi^{\primeid}_m,\bigl(y,\loclev{m}/\loclev{}\bigr)^{p^{m-n}}\rangle\\
\text{(by~\cite[Ch.~XI, Proposition 2]{Ser62})}&=\locdeg{m}{n}\invJP\bigl[\overline{y}\cup\cch{1}\chi_m^{\primeid}\bigr]\\
\text{(by~\eqref{eq:ymdef})}&=\locdeg{m}{n}\invJP\bigl[(a \mod{\gorder{G_m(\primeid)}}) \cup \fundclass{\loclev{m}}{\loclev{}}\bigr]\\
&=\invJP\bigl[(a \mod{\gorder{G_m(\primeid)}}) \cup (\locdeg{m}{n}\fundclass{\loclev{m}}{\loclev{}})\bigr]\\
\text{(by~\cite[Ch.~XI,~\S3]{Ser62})}&=\invJP\bigl[(a \mod{\gorder{G_m(\primeid)}}) \cup \infl^2\circ \arext{\loclev{m}}{\loclev{n}}^\ast(\fundclass{\loclev{n}}{\loclev{}})\bigr]\\
&=\invJP\bigl[\infl^{2}\circ\arext{\loclev{m}}{\loclev{n}}^\ast\bigl((a \mod{\gorder{G_n(\primeid)}}) \cup \fundclass{\loclev{n}}{\loclev{}}\bigr)\bigr]\\
&=\invJP\bigl[(a \mod{\gorder{G_n(\primeid)}}) \cup \fundclass{\loclev{n}}{\loclev{}}\bigr]
\end{align*}
where the last equality follows from the definition of the absolute $\invJP$ as the limit, relative to inflation, of $\invJP$ at finite levels.
Since $\invJP$ is injective, we conclude that
\[
\HHm{2}{\arnm{\loclev{m}}{\loclev{n}}}\bigl((a \mod{\gorder{G_m(\primeid)}}) \cup \fundclass{\loclev{m}}{\loclev{\empty}}\bigr) = (a \mod{\gorder{G_n(\primeid)}}) \cup \fundclass{\loclev{n}}{\loclev{\empty}}
\]
showing that the second diagram is commutative.
\end{proof}

\begin{remark}
Consider the normic system ${{\ZZ_p}_{\pallino[\empty]}}=\bigl(\ZZ_p,\id{\empty},\locdeg{m}{n}\bigr)_{m\geq n\geq 0}$. Then in the left (\rsp in the right) diagram of the statement of Lemma~\ref{lemma:fundclasscommute}, the left vertical map is in fact $\HHm{0}{\id{\empty}}$ (\rsp $\HHm{0}{\locdeg{m}{n}}$). In particular, the lemma can be reformulated as saying that the double systems $\sysHH[\primeid]{0}{{{\ZZ_p}_{\pallino[\empty]}}}$ and $\sysHH[\primeid]{2}{L_{\pallino[\Primeid]}^\times}$ are isomorphic, an isomorphism being the cup-product with the fundamental classes.
\end{remark}

\begin{corollary}\label{cor:indunits_12}
The isomorphisms~\eqref{eq:cohom_trivial_D_action} define isomorphisms of $\DD$-double systems 
\begin{equation}\label{eq:idunits_12_permutation}
\sysHH{i}{\idunit{\pallino}}\bddcong\calZ^{\Ram{L_\infty}{F}}\qquad \text{ for }i=1,2
\end{equation}
where $\calZ=\calZ(\Gamma)$ denotes the double system introduced in Example~\ref{example:shift} and where $\splitG$ acts on $\calZ^{\Ram{L_\infty}{F}}$ through its action on the set $\Ram{L_\infty}{F}$. In particular, there are $\bddsys{\DD}$-isomorphisms
\[
\sysHH{i}{\idunit{\pallino}}\bddcong\limemb\bigl(\sysHH{i}{\idunit{\pallino}}\bigr)\bddcong\calZ^{\Ram{L_\infty}{F}}\qquad \text{ for }i=1,2.
\]
\end{corollary}
\begin{proof} 
For a prime $\overline{\primeid}\in\Ram{L_\infty}{k}$, the maps~\eqref{eq:cohom_trivial_D_action} arise from isomorphisms of $\splitG$-modules 
\begin{equation}\label{eq:finlevlocunits}
\bigoplus_{\primeid\mid\overline{\primeid}}\HH{i}{G_n(\primeid)}{\unit{\Primeid{n}}}\cong \bigoplus_{\primeid\mid\overline{\primeid}}\ZZ/\ramdeg{n}{0}
\end{equation}
induced by ~\eqref{eq:h1loc} and \eqref{eq:h2loc}. Considering the groups on the right-hand side as the terms of the $\DD$-double system $\oplus_{\primeid\mid\overline{\primeid}}\calZ(\Gamma)_{[\nprimes+\logindeg]}$, we are going to prove that
\begin{equation}\label{eq:ddsyslocunits}
\bigoplus_{\primeid\mid\overline{\primeid}}\sysHH[\primeid]{i}{\unit{\pallino[\Primeid]}} \cong \bigoplus_{\primeid\mid\overline{\primeid}}\calZ(\Gamma)_{[\nprimes+\logindeg]}
\end{equation}
as $\DD$-double systems: this is enough to prove the main assertion of the corollary because, arguing as in Example~\ref{example:shift}, we see that 
\[
\bigoplus_{\primeid\mid\overline{\primeid}}\calZ(\Gamma)_{[\nprimes+\logindeg]} \bddcong \bigoplus_{\primeid\mid\overline{\primeid}}\calZ(\Gamma)
\]
where $\splitG$ acts by conjugation on the set $\{\primeid\mid \overline{\primeid}\}$. 

To prove \eqref{eq:ddsyslocunits}, we need to check that isomorphisms \eqref{eq:finlevlocunits} are compatible with the transition maps of the double systems. Since transition maps are defined component-wise, we can perform this check on the components corresponding to a fixed prime $\primeid$ above $\overline{\primeid}$. 

The sequences~\eqref{seq:exact_of_val} induce a short exact sequence of $\DD(\primeid)$-normic systems
\begin{equation}\label{seq:locunitsl}
0\longrightarrow \unit{\pallino[\Primeid]}\longrightarrow L_{\pallino[\Primeid]}^\times\longrightarrow \bigl(\ZZ_p,\nsu_{n,m}=\ramdeg{m}{n},\nsd_{m,n}=\indeg{m}{n}\bigr)\longrightarrow 0.
\end{equation}
The first two normic systems in~\eqref{seq:locunitsl} satisfy both~\ref{cond:inj} and~\ref{cond:surj}, and the third satisfies~\ref{cond:inj}. By applying Propositions~\ref{prop:long_cohom_up} and~\ref{prop:long_cohom_down} to~\eqref{seq:locunitsl} we obtain commutative diagrams
\[
\begin{array}{lcr}
\xymatrix@C=3.5pc{
	\HHN{0}{G_m(\primeid)}{\ZZ_p}\ar@{->}[0,1]^{\cch{0}}
	&\HH{1}{G_m(\primeid)}{\unit{\Primeid{m}}}\\
	\HHN{0}{G_n(\primeid)}{\ZZ_p}\ar@{->}[0,1]^{\cch{0}}\ar@{->}[-1,0]_{\HHm{0}{\nsu_{n,m}}}
	&\HH{2}{G_n(\primeid)}{\unit{\Primeid{n}}}\ar@{->}[-1,0]_{\HHm{2}{\arext{\Primeid{m}}{\Primeid{n}}}}
}
&\parbox[t][3.75em][c]{3.5em}{\begin{center}\text{and}\end{center}}&
\xymatrix@C=3.5pc{
	\HHN{0}{G_m(\primeid)}{\ZZ_p}\ar@{->}[0,1]^{\cch{0}}
	&\HH{2}{G_m(\primeid)}{\unit{\Primeid{m}}}\\
	\HHN{0}{G_n(\primeid)}{\ZZ_p}\ar@{->}[0,1]^{\cch{0}}\ar@{<-}[-1,0]_{\HHm{0}{\nsd_{n,m}}}
	&\HH{2}{G_n(\primeid)}{\unit{\Primeid{n}}}\ar@{<-}[-1,0]_{\HHm{2}{\arnm{\Primeid{m}}{\Primeid{n}}}}
}
\end{array}\]
Note that the cohomological connecting homomorphisms are precisely the maps inducing the isomorphisms \eqref{eq:h1loc}, so we have shown the required compatibility for $i=1$.

Concerning the case $i=2$, the sequence \eqref{seq:locunitsl} induces an injection $\sysHH[\primeid]{2}{\unit{\pallino[\Primeid]}}\hookrightarrow \sysHH[\primeid]{2}{L^\times_{\Primeid{\pallino}}}$ of double systems. The isomorphisms~\eqref{eq:h2loc} are the restriction to $\HH{2}{G_n(\primeid)}{\unit{\Primeid{n}}}$ of the reciprocity isomorphisms
\begin{equation}\label{eq:locrecisol}
\theta_n\colon\HH{2}{G_n(\primeid)}{L_{\Primeid{n}}^\times}\cong\HH{0}{G_n(\primeid)}{\ZZ_p}.
\end{equation}
Hence, it is enough to show that the isomorphisms in~\eqref{eq:locrecisol}, commute with the transition morphisms for arbitrary $m\geq n\geq 0$, which is precisely Lemma \ref{lemma:fundclasscommute}.

The final statement concerning $\limemb$ follows from the first combined with Example~\ref{example:Q/Z}.
\end{proof}

The following theorem is the main result of this section.

\begin{theorem}\label{thm:i_cinque_dell'apocalisse} 
In the quotient category $\quotsys[\cofg]{\DD}$ the following sequence is exact
\begin{equation}\label{eq:long_cohom_quot}
0\longdashrightarrow\sysHH{1}{\unit{\pallino}}\longdashrightarrow\sysHH{1}{\idunit{\pallino}}\longdashrightarrow\sysHH{1}{\q{\pallino}}\longdashrightarrow\sysHH{2}{\unit{\pallino}}\longdashrightarrow\sysHH{2}{\idunit{\pallino}}.
\end{equation}
By applying the functor $\limemb=\stepfunctor\circ\varinjlim$, it gives rise to the exact sequence
\begin{equation}\label{eq:long_cohom_limemb}
0\longdashrightarrow\limemb\bigl(\sysHH{1}{\unit{\pallino}}\bigr)\longdashrightarrow\limemb\bigl(\sysHH{1}{\idunit{\pallino}}\bigr)\longdashrightarrow\limemb\bigl(\sysHH{1}{\q{\pallino}}\bigr)\longdashrightarrow\limemb\bigl(\sysHH{2}{\unit{\pallino}}\bigr)\longdashrightarrow\limemb\bigl(\sysHH{2}{\idunit{\pallino}}\bigr)
\end{equation}
and there is a commutative diagram of exact sequences
\begin{equation}\label{diag:comm_limemb}
\begin{split}\xymatrix@C1.75em{
0\ar@{-->}[0,1]&\sysHH{1}{\unit{\pallino}}\ar@{-->}[0,1]\ar@{-->}[1,0]&\sysHH{1}{\idunit{\pallino}}\ar@{-->}[0,1]\ar@{-->}[0,1]\ar@{-->}[1,0]^{\cong}& \sysHH{1}{\q{\pallino}}\ar@{-->}[0,1]\ar@{^(-->}[1,0]&\sysHH{2}{\unit{\pallino}}\ar@{-->}[0,1]\ar@{-->}[1,0]&\sysHH{2}{\idunit{\pallino}}\ar@{-->}[1,0]^{\cong}\\
0\ar@{-->}[0,1]&\limemb\bigl(\sysHH{1}{\unit{L_\infty}}\bigr)\ar@{-->}[0,1]& \limemb\bigl(\sysHH{1}{\idunit{L_\infty}}\bigr)\ar@{-->}[0,1]& \limemb\bigl(\sysHH{1}{\q{L_\infty}}\bigr)\ar@{-->}[0,1]&\limemb\bigl(\sysHH{2}{\unit{L_\infty}}\bigr)\ar@{-->}[0,1]&\limemb\bigl(\sysHH{2}{\idunit{L_\infty}}\bigr)
}
\end{split}
\end{equation}
\end{theorem}
\begin{proof} Since the sequence~\eqref{eq:long_cohom_DS} is already exact in $\dblsys[\cofg]{\DD}$, exactness of~\eqref{eq:long_cohom_quot} amounts to
\begin{equation}\label{eq:ker_H1_bdd}
\kernel\Bigl(\sysHH{1}{\unit{\pallino}}\longrightarrow\sysHH{1}{\idunit{\pallino}}\Bigr)\in\bddsys{\DD}.
\end{equation}
To show this, observe that the ascending transition morphisms of this system are the restriction of 
\[
\HHm{1}{\nsu}\colon \HH{1}{G_n}{\unit{L_n}}\longrightarrow \HH{1}{G_m}{\unit{L_m}}.
\]
Lemma~\ref{lemma:HHm(j)_and_inf}, together with injectivity of inflation in degree $1$, implies that $\HHm{1}{\nsu}$ is injective, because $\nsu^\ast$ is the identity in this case. Moreover, for every $n\geq 0$, the kernel $\kernel\bigl(\HH{1}{G_n}{\unit{L_n}}\to\HH{1}{G_n}{\idunit{L_n}}\bigr)$ is isomorphic to the capitulation kernel $\kernel\bigl(\clsyl{F}\to\clsyl{L_n}\bigr)$ (see, for instance,~\cite[Proposition~2.2]{Nuc10}), and therefore its order is constant for $n$ large enough. It follows that the maps $\HHm{1}{\nsu}$ are isomorphisms for $n\gg 0$, and Proposition~\ref{prop:trivialbdd} yields~\eqref{eq:ker_H1_bdd}. Since $\limemb$ is an exact functor when restricted to $\quotsys[\cofg]{\DD}$ by Proposition~\ref{prop:Sexact}, exactness of~\eqref{eq:long_cohom_limemb} follows. By applying Corollary~\ref{cor:natural_transformation} to~\eqref{eq:long_cohom_limemb}, we find the commutative diagram~\eqref{diag:comm_limemb}, where the second and fifth vertical maps are isomorphisms by Corollary~\ref{cor:indunits_12}. What remains to be proved is that the third vertical morphism of~\eqref{diag:comm_limemb} is injective. To check this, recall that by Lemma~\ref{lemma:HHm(j)_and_inf} the following triangle commutes
\[\xymatrix{
	\HH{1}{G_n}{\q{L_n}}\ar@{->}[0,2]^{\HHm{1}{\nsu_{n,m}}}\ar@{->}[1,1]_{\nsu_{n,m}^{\ast,1}}&&\HH{1}{G_m}{\q{L_m}}\\
	&\HH{1}{G_n}{\q{L_m}^{G_{m,n}}}\ar@{->}[-1,1]_{\infl}
}\]
for all $m\geq n\geq 0$. In particular, $\kernel\bigl(\HHm{1}{\nsu_{n,m}}\bigr)=\kernel(\nsu_{n,m}^{\ast,1})$, because inflation is injective in degree $1$. Taking $G_n$-cohomology of the exact sequence
\begin{equation*}
0\longrightarrow \q{L_n}\longrightarrow \q{L_m}^{G_{m,n}}\longrightarrow \q{L_m}^{G_{m,n}}/\q{L_n}\longrightarrow 0
\end{equation*}
shows that $\kernel(\nsu_{n,m}^{\ast,1})$ is a quotient of a submodule of $\q{L_m}^{G_{m,n}}/\q{L_n}$. Taking direct limits, the naturality of the connecting homomorphism guarantees that in fact 
\[
\kernel\bigl(\sysHH{1}{\q{\pallino}} \longrightarrow \limemb\bigl(\sysHH{1}{\q{\pallino}}\bigr)\bigr)\quad \text{ is a quotient of a subsystem of}\quad \bigl(\q{L_\infty}^{\Gamma_n}/\q{L_n},\arext{L_m}{L_n},\arnm{L_m}{L_n}\bigr).
\]
Now, $\q{L_m}^{G_{m,n}}/\q{L_n}$ is isomorphic to $\kernel\bigl(\clsyl{L_n}\to\clsyl{L_m}\bigr)$, as shown in~\cite[proof~of Proposition~2.2]{Nuc10}, the isomorphism~\ibid being the connecting homomorphism of the snake lemma. Again these isomorphisms give rise to an isomorphism of double systems 
\[
\bigl(\q{L_\infty}^{\Gamma_n}/\q{L_n},\arext{L_m}{L_n},\arnm{L_m}{L_n}\bigr)\cong \bigl(\kernel(\clsyl{L_n}\to\clsyl{L_\infty}),\arext{L_m}{L_n},\arnm{L_m}{L_n}\bigr)
\]
thanks to the naturality of the snake lemma. Hence, in order to prove that $\nsu^{\ast,1}=(\nsu^{\ast,1}_{n,m})$ is a $\bddsys{\DD}$-monomorphism, it is enough to show that the double system $\bigl(\kernel(\clsyl{L_n}\to\clsyl{L_\infty}),\arext{L_m}{L_n},\arnm{L_m}{L_n}\bigr)$ belongs to $\bddsys{\DD}$. This follows from Proposition~\ref{prop:trivialbdd}, because the norm is eventually an isomorphism on capitulation kernels as shown, for example, in~\cite[Proposition]{Oza95}.
\end{proof}

\begin{remark}\label{rmk:quel_vecchio_volpone_del_Nuccio} In the Introduction we mentioned that two morphisms in~\eqref{diag:comm_limemb} become injective in the quotient category $\quotsys[\cofg]{\DD}$ but their kernels do not belong to the corresponding thick subcategories when regarded simply as direct (\rsp inverse) systems. The morphisms in questions are $\sysHH{1}{\q{\pallino}} \to \limemb\bigl(\sysHH{1}{\q{\pallino}}\bigr)$ and $\sysHH{1}{\unit{\pallino}}\to\sysHH{1}{\idunit{\pallino}}$, respectively. This is the main reason that led us to introduce the notion of double systems.
\end{remark}

We now move forward to the study the global counterpart of Corollary~\ref{cor:indunits_12}. Before, we need the following technical lemma:

\begin{lemma}\label{lemma:mimi_omologico_ferito_nell'onore} 
Let $A_1,\dots,A_5,B_1,\dots,B_5$ be objects of an abelian category and let $f_1,\dots,f_5$ be morphisms in the category making the following diagram with exact rows commute:
\[\xymatrix@C=3em{
	&A_1\ar@{->}[1,0]^{f_1}\ar@{->}[0,1]&A_2\ar@{->}[1,0]^{f_2}\ar@{->}[0,1]&A_3\ar@{->}[1,0]^{f_3}\ar@{->}[0,1]&
	A_4\ar@{->}[1,0]^{f_4}\ar@{->}[0,1]&A_5\ar@{->}[1,0]^{f_5}\\
	0\ar@{->}[0,1]&B_1\ar@{->}[0,1]&B_2\ar@{->}[0,1]&B_3\ar@{->}[0,1]&
	B_4\ar@{->}[0,1]&B_5
}\]
Suppose that $f_2$ is an isomorphism and $f_5$ is injective. Then there is an exact sequence
\[
0\longrightarrow\cokernel{f_1} \longrightarrow \kernel{f_3} \longrightarrow \kernel{f_4}\longrightarrow 0.
\]
In particular, if $f_3$ is injective, then $f_4$ is injective as well and $f_1$ is surjective.
\end{lemma}
\begin{proof}
Completing with $0$ the rows of the diagram of the statement, we obtain cochain complexes, which we denote by $\mathscr{A}$ and $\mathscr{B}$ respectively, satisfying $\mathscr{A}_i=A_i$ and $\mathscr{B}_i=B_i$  for $1\leq i\leq 5$. Thus we have a map $f\colon\mathscr{A}\to \mathscr{B}$ of cochain complexes and we consider the cochain complexes $\mathscr{K}=\kernel{f}$, $\mathscr{I}=\im{f}$ and $\mathscr{C}=\cokernel{f}$. We obtain two exact sequences of complexes 
\begin{equation}\label{eq:exact_cochain_cpx}
0 \longrightarrow \mathscr{K} \longrightarrow \mathscr{A} \longrightarrow \mathscr{I} \to 0\qquad\text{ and }\qquad 0 \longrightarrow\mathscr{I} \longrightarrow\mathscr{B} \longrightarrow\mathscr{C} \longrightarrow 0.
\end{equation}
The exactness of the rows of the diagram in the statement can be recast as $H^i(\mathscr{A})=0$ for all $i\neq 1,5$ and $H^i(\mathscr{B})=0$ for all $i\neq 5$. In particular, the zig-zag lemma applied to the exact sequences in~\eqref{eq:exact_cochain_cpx} yields isomorphisms
\begin{equation}\label{eq:zig-zag}
H^1(\mathscr{C})\cong H^2(\mathscr{I})\cong  H^3(\mathscr{K}) \qquad\text{ and }\qquad  H^2(\mathscr{C})\cong H^3(\mathscr{I})\cong H^4(\mathscr{K}).
\end{equation}
The exact sequence of the statement can be obtained from the exact sequence
\begin{equation}\label{eq:first_from_zz_to_lemma}
0\longrightarrow \kernel\bigl(\mathscr{K}_3\longrightarrow\mathscr{K}_4\bigr)\longrightarrow \mathscr{K}_3\longrightarrow \mathscr{K}_4\longrightarrow \mathscr{K}_4/\im\bigl(\mathscr{K}_3\bigr)\longrightarrow 0
\end{equation}
as follows. The assumptions on $f_2$ and $f_5$ mean that $\mathscr{K}_2=\mathscr{K}_5=0$ and $\mathscr{C}_2=0$: in particular, $H^2(\mathscr{C})=0$ and~\eqref{eq:zig-zag} implies $H^4(\mathscr{K})=0$. Therefore the fourth term $\mathscr{K}_4/\im\bigl(\mathscr{K}_3\bigr)$ in~\eqref{eq:first_from_zz_to_lemma} is also trivial since it is isomorphic to $H^4(\mathscr{K})$, since $\mathscr{K}_5=0$. Similarly, $\mathscr{K}_2=0$ implies $H^3(\mathscr{K})=\kernel\bigl(\mathscr{K}_3\to\mathscr{K}_4\bigr)$, and invoking again~\eqref{eq:zig-zag} we can rewrite~\eqref{eq:first_from_zz_to_lemma} as
\begin{equation*}
0\longrightarrow H^1(\mathscr{C})\longrightarrow \mathscr{K}_3=\kernel{f_3}\longrightarrow \mathscr{K}_4=\kernel{f_4}\longrightarrow 0.
\end{equation*}
The fact that $H^1(\mathscr{C})$ coincides with $\mathscr{C}_1=\cokernel{f_1}$ follows again from $\mathscr{C}_2=0$, finishing the proof.
\end{proof}

\begin{corollary}\label{cor:growth_cohom_units}
The $\DD$-double systems $\sysHH{i}{\unit{\pallino}}$ and $\limemb\bigl(\sysHH{i}{\unit{\pallino}}\bigr)=\stepfunctor(\HHN{i}{\Gamma}{\unit{L_\infty}})$ are $\bddsys{\DD}$-isomorphic for $i=1,2$. 
\end{corollary}
\begin{proof}
Apply Lemma~\ref{lemma:mimi_omologico_ferito_nell'onore} to the commutative diagram~\eqref{diag:comm_limemb} of Theorem~\ref{thm:i_cinque_dell'apocalisse}: we obtain a surjection
\begin{equation}\label{eq:surj_uno_globali}
\sysHH{1}{\unit{\pallino}}\longtwoheaddashedarrow\limemb\bigl(\sysHH{1}{\unit{\pallino}}\bigr)
\end{equation}
as well as an injection
\begin{equation}\label{eq:inj_due_globali}
\sysHH{2}{\unit{\pallino}}\longhookdashedarrow\limemb\bigl(\sysHH{2}{\unit{\pallino}}\bigr).
\end{equation}
The morphism in~\eqref{eq:surj_uno_globali} is clearly also injective (and therefore an isomorphism), because so is the composition
\[
\sysHH{1}{\unit{\pallino}}\longdashrightarrow\limemb\bigl(\sysHH{1}{\unit{\pallino}}\bigl)\longdashrightarrow\limemb\bigl(\sysHH{1}{\idunit{\pallino}}\bigr),
\]
as shown in~Theorem~\ref{thm:i_cinque_dell'apocalisse}.

Passing to cohomology in degree $2$, consider the morphisms, for $m\geq n\geq 0$,
\begin{align*}
\HHm{0}{\arext{L_m}{L_n}}\colon\unit{F}/\arnm{L_n}{F}\unit{L_n} &\longrightarrow \unit{F}/\arnm{L_m}{F}\unit{L_m}\\
x& \longmapsto x^{p^{m-n}}
\end{align*}
Set
\[
\tau_n=\varinjlim_m \HHm{0}{\arext{L_m}{L_n}}\colon \unit{F}/\arnm{L_n}{F}\unit{L_n} \longrightarrow \tors[p^n]{\bigl(\varinjlim_m\unit{F}/\arnm{L_m}{F}\unit{L_m}\bigr)}=\bigl(\limemb\bigl(\sysHH{0}{\unit{\pallino}}\bigr)\bigr)_n.
\]
By definition of $\HHm{2}{\arext{L_m}{L_n}}$ as $\HHm{2}{\arext{L_m}{L_n}}=\tateiso{m}{}\circ\HHm{0}{\arext{L_m}{L_n}}\circ\tateiso{n}{-1}$ (see~\S \ref{subsubsec:arbitrarydegree}), $\bddsys{\DD}$-surjectivity in~\eqref{eq:inj_due_globali} is equivalent to the claim that the double system $\calC$ with components $(\calC)_n=\cokernel(\tau_n)$ is in $\bddsys{\DD}$, or also that the image of $\tau=(\tau_n)$ in the quotient category $\quotsys{\DD}$ is a surjective morphism. Consider the $\DD$-double system
\begin{equation}\label{def:double_sys_OF/p}
\unit{F}/(\unit{F})^{\pallino[\expp]}=\left(\unit{F}/(\unit{F})^{p^n},\dsu_{n,m}\colon x\longmapsto x^{p^{m-n}} \mod{(\unit{F})^{p^m}},\dsd_{m,n}\colon x\longmapsto x \mod{(\unit{F})^{p^n}}\right)_{m\geq n\geq 0}.
\end{equation}
Let $\sysnorm{\unit{\pallino}}$ be the double subsystem of $\unit{F}/(\unit{F})^{\pallino[\expp]}$ whose components are $\arnm{L_n}{F}\unit{L_n}/(\unit{F})^{p^n}$. By the structure theorem of units in number fields, one immediately sees that $\unit{F}/(\unit{F})^{\pallino[\expp]}$ is $\bddsys{\DD}$-isomorphic to $\calZ^{\rank_{\ZZ_p}\unit{F}}$ and in particular belongs to $\dblsys[\cofg]{\DD}$ (and the same holds for $\sysnorm{\unit{\pallino}}$). We can therefore apply Proposition~\ref{prop:Sexact} and Corollary~\ref{cor:natural_transformation} to the tautological exact sequence $0\to\sysnorm{\unit{\pallino}}\to\unit{F}/(\unit{F})^{\pallino[\expp]}\to\sysHH{0}{\unit{\pallino}}\to 0$ and obtain a commutative diagram in $\quotsys[\cofg]{\DD}$
\begin{equation*}\begin{split}
\xymatrix{
0\ar@{-->}[0,1]&\sysnorm{\unit{\pallino}}\ar@{-->}[0,1]\ar@{-->}[1,0]&
	\unit{F}/(\unit{F})^{\pallino[\expp]}\ar@{-->}[0,1]\ar@{-->}[1,0]&
	\sysHH{0}{\unit{\pallino}}\ar@{-->}[1,0]^{\tau}\ar@{-->}[0,1]&0\\
0\ar@{-->}[0,1]&\limemb\bigl(\sysnorm{\unit{\pallino}}\bigr)\ar@{-->}[0,1]&
	\limemb\bigl(\unit{F}/(\unit{F})^{\pallino[\expp]}\bigr)\ar@{-->}[0,1]&
	\limemb\bigl(\sysHH{0}{\unit{\pallino}}\bigr)\ar@{-->}[0,1]&0
}\end{split}
\end{equation*}
Example~\ref{example:Q/Z} shows that the middle vertical map is an isomorphism. The snake lemma implies that the right vertical morphism is surjective, concluding the proof that~\eqref{eq:inj_due_globali} is an isomorphism.
\end{proof}

\begin{remark}\label{rmk:Iwasawa83}
In \cite{Iwa83} Iwasawa obtains a general description of the cohomology groups $\HHN{i}{\Gamma}{\unit{L_\infty}}$ ``up to finite groups'', which becomes explicit under the assumption that $F$ is totally real and that it satisfies Leopoldt conjecture (see~\cite[Proposition~5 and Section~4]{Iwa83}). His analysis is performed only at the limit, and conveys no information on the full ``double system'' (a notion which does not appear \ibid). One of our motivations for this work was to unravel the purely group-theoretical content of his results so that it could be generalized to other settings and other normic systems. More precisely, to obtain an Iwasawa-like formula for the growth of the orders of cohomology groups of units from the description of the cohomology groups $\HHN{i}{\Gamma}{\unit{L_\infty}}$, one is naturally led to consider the morphisms $\lambda_n$ and $\mu_n$ of~\cite[p.~199]{Iwa83}, which are precisely our maps $\sysHH{i}{\unit{\pallino}}\to\limemb\bigl(\sysHH{i}{\unit{\pallino}}\bigr)$. The groups $\HHN{i}{\Gamma}{\unit{L_\infty}}[p^n]=\bigl(\limemb\bigl(\sysHH{i}{\unit{\pallino}}\bigr)_n$ follow an Iwasawa-like formula, since
\[
\HHN{i}{\Gamma}{\unit{L_\infty}} \cong (\QQ_p/\ZZ_p)^{r_i} \oplus \text{(finite group)}
\] 
by Iwasawa's description. One can deduce from this that the cohomology groups of units at finite levels also follows an Iwasawa-like formula if both kernels and cokernels of $\lambda_n$ and $\mu_n$ are eventually constant: Corollary~\ref{cor:growth_cohom_units} makes this precise. Formul\ae{} for the cohomology groups of units at finite levels will be made explicit in the next section.
\end{remark}

\section{\texorpdfstring{Fake $\ZZ_p$-extensions of dihedral type}{Fake extensions of dihedral type}}\label{sec:prodihedral}
\subsection{Definitions}\label{subsec:fake_def} In the rest of the paper assume that $p\geq 3$ is odd, and that $\DD$ is pro-dihedral: by this, we mean that
$\splitG$ has order $2$ and the split exact sequence
\begin{equation}\tag{\ref{eq:seq_gal_groups}}
1 \longrightarrow \Gamma \longrightarrow \DD \longrightarrow  \splitG\longrightarrow 1
\end{equation}
induces a non-trivial action of $\splitG$ on $\Gamma$. 
Many results we obtain are likely to easily generalise to the case where~$\splitG$ is a finite group of order prime to $p$ whose characters are $\ZZ_p$-valued.

The image of a fixed splitting $\splitG\to\DD$ is a subgroup $\subsplit\subseteq\DD$, which intersects each $\Gamma_n$ trivially. It is hence isomorphic to its image $\subsplit\subseteq \splitG[n]$, still denoted by the same symbol by a slight abuse of notation. For $n\leq \infty$, let $K_n$ be the subfield of $L_n$ fixed by $\subsplit$ and write $k=K_0$. The corresponding diagram of fields is the following:
\begin{equation}\label{diag:fake}\begin{split}
\xymatrix@R=1.5pc{
L_\infty\ar@{-}[1,1]_{\subsplit}\ar@{-}[2,0]_{\Gamma_n}\ar@/_2.25pc/@{-}[4,0]_{\Gamma}\ar@/^5.25pc/@{-}[5,1]^{\DD\cong\splitG\ltimes\Gamma}
\\
&K_\infty\ar@{--}[2,0]\\
L_n\ar@{-}[1,1]\ar@{-}[2,0]_{G_n}\ar@{-}[3,1]_{\splitG[n]}\\
&K_n\ar@{--}[2,0]\\
F\ar@{-}[1,1]_{\splitG}\\
&k
}\end{split}\end{equation}
The collection of fields $\{K_n\}_{n\geq 0}$ share some similarities with a $\ZZ_p$-extension: indeed, each extension $K_{n+1}/K_n$ has degree $p$ and the only possible ramified primes lie above $p$, because so is for the extension $L_\infty/F$. Yet, a crucial difference is that none of the extension $K_{n+1}/K_n$, let alone the whole $K_\infty/k$, is Galois. This motivates the following
\begin{definition}\label{def:fake}
Assume that $p\geq 3$ is odd and consider the setting of diagram~\eqref{diag:fake}. The extension $K_\infty/k$ is a fake $\ZZ_p$-extension of dihedral type, the field $L_\infty$ is the Galois closure of the fake $\ZZ_p$-extension and $F/k$ is called the normalizing quadratic extension.
\end{definition}

If $B$ is a $\ZZ_p[\DD]$-module, then it is uniquely $2$-divisible. By restriction, it admits an action of $\subsplit$, which we can see as an action of $\splitG$ via the splitting $\splitG\to\DD$. Accordingly, it admits a functorial decomposition $B=B^+\oplus B^-$ where we denote by $B^+$ (\rsp $B^-$) the maximal submodule on which $\splitG$ acts trivially (\rsp as $-1$). Such decomposition is obtained by writing $b=\frac{1+\sigma}{2}b+\frac{1-\sigma}{2}b$ for $b\in B$, where $\sigma$ generates $\splitG$. Then 
\begin{equation}\label{2divid}
B^{\splitG} =B^+\cong B/B^-=B_{\splitG}.
\end{equation}
This applies in particular to $\ZZ_p[\splitG[n]]$-modules and their $G_n$-cohomology. In fact we can say a bit more in this case, as explained in the next proposition which is \cite[Proposition 2.1]{CapNuc20} (the final statement about $\Gamma$-cohomology follows from the result at finite levels by taking direct limits):

\begin{proposition}\label{prop:fixcohom} 
Let $B$ be a uniquely $2$-divisible $\splitG[n]$-module. Then, for every $i\in \ZZ$, restriction induces functorial isomorphisms 
\[
\HH{i}{\splitG[n]}{B}\cong \HH{i}{G_n}{B}^+.
\]
Moreover, the Tate isomorphism $\HH{i}{G_n}{B}\cong \HH{i+2}{G_n}{B}$ is $\splitG$-antiequivariant. In particular, for every $i\in\ZZ$, there are isomorphisms
\begin{align*}
\HH{i}{G_n}{B}^-&\cong\HH{i+2}{G_n}{B}^+\cong \HH{i+2}{\splitG[n]}{B}\\ 
\HH{i}{G_n}{B}^+&\cong\HH{i+2}{G_n}{B}^-\cong \HH{i}{\splitG[n]}{B}\end{align*}
of abelian groups. In particular, given any $\DD$-normic system $\calB=(B_{L_n})$, there are functorial isomorphisms
\[
\HHN{i}{\Gamma}{B_{L_\infty}}^+\cong\HHN{i}{\DD}{B_{L_\infty}}\qquad\text {for }i\geq 1.
\]
\end{proposition}
Proposition~\ref{prop:fixcohom} induces a decomposition of double systems $\sysHH{i}{\calB}=\sysHH{i}{\calB}^+\oplus\sysHH{i}{\calB}^-$, where, for every $n$, $\bigr(\sysHH{i}{\calB}^\pm\bigl)_n=\HH{i}{G_n}{B_{L_n}}^\pm$; an entirely analogous argument applies to the double system of the cohomology of a local normic system, inducing a decomposition $\sysHH[\primeid]{i}{\calB}=\sysHH[\primeid]{i}{\calB}^+\oplus\sysHH[\primeid]{i}{\calB}^-$ whenever $\DD(\primeid)$ is pro-dihedral.

\subsection{An Iwasawa-like formula}\label{subsec:formula}
In this section we state and prove Theorem~\ref{thm:iwasawa_formula}, which is one of our main results. Fix a fake $\ZZ_p$-extension of dihedral type $K_\infty/k$ and adopt the conventions introduced in~\S\ref{subsec:fake_def}. For future reference, we add the following piece of notation to Notation~\ref{notation:ram_and_ind}:

\begin{notation}\label{notation:primesets} 
Denote by $\Ram{L_\infty}{k}_s\subseteq\Ram{L_\infty}{k}$ the set of primes of $k$ that ramify in $L_\infty/k$ and split in $F/k$. Also, denote by $\oram{L_\infty}{k}$ and by $\oram{L_\infty}{k}_s$, respectively, the number of primes in $\Ram{L_\infty}{k}$ and in $\Ram{L_\infty}{k}_s$.
\end{notation}

The first step is Corollary~\ref{cor:growth_cohom_units_dihedral} below, which is a signed version of Corollary~\ref{cor:indunits_12} and of Corollary~\ref{cor:growth_cohom_units}. Here, $\calZ=\calZ(\Gamma)$ again denotes the double system introduced in Example~\ref{example:shift}.

For a group $H$ and a $H$-module $B$, we denote by $\hh{i}{H}{B}$ the order of $\HH{i}{H}{B}$, whenever this group is finite.

\begin{corollary}\label{cor:growth_cohom_units_dihedral}
\leavevmode \begin{enumerate}[label=\textit{\roman*}\textup{)}]
\item \label{cor:growth_cohom_units_dihedral:pt_idelic} For $i=1,2$, there are $\bddsys{\Gamma}$-isomorphisms  $\sysHH{i}{\idunit{\pallino}}^+\cong_{\bddsys{\Gamma}}\calZ^{\Ram{L_\infty}{k}}$ and  $\sysHH{i}{\idunit{\pallino}}^-\cong_{\bddsys{\Gamma}}\calZ^{\Ram{L_\infty}{k}_s}$ that induce ${\catMod[]{fin}}$-isomorphisms
\begin{equation}\label{eq:iso_idunits_QZ_signed}
\HH{i}{\DD}{\idunit{L_\infty}}\bddcong[{\catMod[]{fin}}] \bigr(\QQ_p/\ZZ_p\bigl)^{\Ram{L_\infty}{k}}\qquad\text{and}\qquad
\HH{i}{\Gamma}{\idunit{L_\infty}}^-\bddcong[{\catMod[]{fin}}] \bigr(\QQ_p/\ZZ_p\bigl)^{\Ram{L_\infty}{k}_s}\qquad i=1,2.
\end{equation}
\label{cor:growth_cohom_id_units_dihedral_h0}
\item For each choice of a symbol $\ast\in\{\emptyset,\pm \}$ and every $i\in\ZZ$, there exists integers $\lambdacohom{i}^\ast\in\ZZ_{\geq 0}$ such that
\[
\hh{i}{G_n}{\unit{L_n}}^\ast\bddeq p^{n\cdot \lambdacohom{i}^\ast}.
\]\label{cor:growth_cohom_units_dihedral_existence_global}
\item \label{cor:growth_cohom_units_dihedral:various_formulae} For every $i\in\ZZ$ these invariants satisfy
\begin{align}
\lambdacohom{i}^\pm&=\lambdacohom{i+4}^\pm\label{eq:lambdaDperiodic}\\
\intertext{as well as}
\lambdacohom{i}&=\lambdacohom{i}^{+}+\lambdacohom{i}^{-}\label{eq:lambdaGD}\\
\lambdacohom{i}^\pm&=\lambdacohom{i+2}^\mp\label{eq:lambdaGperiodic_signed}\\
\hh{i}{\splitG[n]}{\unit{L_n}}&\bddeq p^{n\cdot \lambdacohom{i}^+}\label{eq:lambdaD+}\\
\lambdacohom{2i} &=\lambdacohom{2i+1}  - 1. \label{eq:lambdaherbrand}\\
\intertext{In particular,}
\label{eq:cohom_unit_many_calcs}
\lambdacohom{2i-1}^+ - \lambdacohom{2i}^+ &=  \lambdacohom{2i-2}^+- \lambdacohom{2i+1}^+ + 1.
\end{align}
\end{enumerate}
\end{corollary}
\begin{proof}
\begin{itemize}[wide, labelindent=0pt]
\item[\ref{cor:growth_cohom_id_units_dihedral_h0}]
Recall that Corollary~\ref{cor:indunits_12} yields $\bddsys{\DD}$-isomorphisms
\[
\sysHH{i}{\idunit{\pallino}}\bddcong\calZ^{\Ram{L_\infty}{F}}\qquad \text{ for }i=1,2
\]
where $\splitG$ acts by permutation on the set $\Ram{L_\infty}{F}$. We are thus reduced to analyse the objects $\bigr(\calZ^{\Ram{L_\infty}{F}}\bigl)^\pm$. Starting with the plus component, let
\begin{equation}\label{eq:x_generic_H}
x=\bigl(x_{\Primeid{}}[\Primeid{}]\bigr)_{\Primeid{}\in\Ram{L_\infty}{F}}\qquad \text{ with }x_{\Primeid{}}\in\calZ
\end{equation}
be an element of $\bigl(\calZ^{\Ram{L_\infty}{F}}\bigr)^+$. For primes $\Primeid{}$ above a non-split prime in $F/k$ we have $\sigma\Primeid{}= \Primeid{}$, while for primes in $\Ram{L_\infty}{F}$ dividing a prime that splits in $F/k$ we have $\sigma\Primeid{}\neq \Primeid{}$. The condition $\sigma x=x$ translates into
\begin{equation}\label{eq:x_plus_H}
x_{\Primeid{}}=x_{\sigma\Primeid{}}\qquad\text{ for every }\Primeid{}\in \Ram{L_\infty}{F}\text{ dividing a split prime in }F/k
\end{equation}
and no constraint on the coefficients $x_{\Primeid{}}$ when $\Primeid{}$ lies above a non-split prime. Let $\primeid{}\in\Ram{L_\infty}{k}$ denote the prime below $\Primeid{}\in\Ram{L_\infty}{F}$: the surjection sending $x$ to the element
\[
\bigl(x_{\Primeid{}}[\primeid{}]\bigr)_{\primeid{}\in\Ram{L_\infty}{k}}
\]
is injective thanks to~\eqref{eq:x_plus_H}, so defines an isomorphism $\bigr(\calZ^{\Ram{L_\infty}{F}}\bigl)^+\cong\calZ^{\Ram{L_\infty}{k}}$. Concerning the minus component, the same analysis as above shows that an element $x$ as in~\eqref{eq:x_generic_H} satisfies $\sigma x=-x$ if and only if
\[
x_{\Primeid{}}=
\begin{cases}
0&\text{ if }\Primeid{}\text{ lies above a non-split prime in }F/k;\\
-x_{\sigma\Primeid{}}&\text{ if }\Primeid{}\text{ divides a prime that splits in }F/k.
\end{cases}
\]
Therefore, the morphism sending $x\in\bigl(\calZ^{\Ram{L_\infty}{F}}\bigr)^-$ to the element
\[
\bigl(x_{\Primeid{}}[\primeid{}]\bigr)_{\primeid{}\in\Ram{L_\infty}{k}_s}
\]
defines an isomorphism $\bigr(\calZ^{\Ram{L_\infty}{F}}\bigl)^-\cong\calZ^{{\Ram{L_\infty}{k}}_s}$. This concludes the proof of the first part of~\ref{cor:growth_cohom_id_units_dihedral_h0}, and~\eqref{eq:iso_idunits_QZ_signed} follows by taking direct limits.
\item[\ref{cor:growth_cohom_units_dihedral_existence_global}] When $\ast=\emptyset$, Corollaries~\ref{cor:ha_vinto_caputo} and~\ref{cor:growth_cohom_units} imply the statement for $i=1,2$: the case for general $i\in\ZZ$ follows by periodicity of Tate cohomology. For the signed invariants, let $i=1,2$ and take $\pm$-parts for the action of $\splitG$ on the isomorphism $\sysHH{i}{\unit{\pallino}}\bddcong\limemb\bigl(\sysHH{i}{\unit{\pallino}}\bigr)$ of Corollary~\ref{cor:growth_cohom_units}: we obtain
\begin{align*}
\sysHH{i}{\unit{\pallino}}^\pm&\bddcong[\bddsys{\Gamma}]\limemb\bigl(\sysHH{i}{\unit{\pallino}}\bigr)^\pm\\
&=\stepfunctor\bigl(\varinjlim\HH{i}{G_n}{\unit{L_n}}\bigr)^\pm\\
&=\stepfunctor\bigl(\bigl(\varinjlim\HH{i}{G_n}{\unit{L_n}}\bigr)^\pm\bigr)\\
&=\stepfunctor\bigl(\varinjlim\bigl(\HH{i}{G_n}{\unit{L_n}}^\pm\bigr)\bigr)\\
&=\limemb\bigl(\sysHH{i}{\unit{\pallino}}^\pm\bigr).
\end{align*}
Corollary~\ref{cor:ha_vinto_caputo} implies that there exist $\lambdacohom{i}^\pm\in\ZZ_{\geq 0}$ such that
\[
\gorder{\HH{i}{G_n}{\unit{L_n}}^\pm}\bddeq p^{n\cdot \lambdacohom{i}^\pm}\qquad\text{ for }i=1,2.
\]
Combining Proposition~\ref{prop:fixcohom} with the periodicity of Tate cohomology permits to extend the definition to all $i$ by setting
\begin{equation}\label{eq:def_lambda_sign}
\lambdacohom{i}^\pm=\begin{cases}
\lambdacohom{2}^\mp&\text{ if }i\equiv 0 \pmod{4}\\
\lambdacohom{1}^\pm&\text{ if }i\equiv 1 \pmod{4}\\
\lambdacohom{2}^\pm&\text{ if }i\equiv 2 \pmod{4}\\
\lambdacohom{1}^\mp&\text{ if }i\equiv 3 \pmod{4}
\end{cases}
\end{equation}
\item[\ref{cor:growth_cohom_units_dihedral:various_formulae}] The equalities in~\eqref{eq:lambdaDperiodic} and~\eqref{eq:lambdaGperiodic_signed} follow by \eqref{eq:def_lambda_sign}. That in~\eqref{eq:lambdaGD} is a direct consequence of the decomposition $\HH{i}{G_n}{\unit{L_n}}=\HH{i}{G_n}{\unit{L_n}}^+\oplus \HH{i}{G_n}{\unit{L_n}}^-$, while~\eqref{eq:lambdaD+} follows from Proposition~\ref{prop:fixcohom}. As for \eqref{eq:lambdaherbrand}, it comes from the fact that the Herbrand quotient of global units satisfies
\[
\frac{\hh{2i}{G_n}{\unit{L_n}}}{\hh{2i+1}{G_n}{\unit{L_n}}} = \frac{\hh{2}{G_n}{\unit{L_n}}}{\hh{1}{G_n}{\unit{L_n}}} = \frac{1}{p^n}
\]
together with the relation
\[
\frac{\hh{2i}{G_n}{\unit{L_n}}}{\hh{2i+1}{G_n}{\unit{L_n}}} \bddeq p^{n(\lambdacohom{2i} - \lambdacohom{2i+1})}
\]
which has been shown in the first part of the proof. The last relation~\eqref{eq:cohom_unit_many_calcs} follows by combining~\eqref{eq:lambdaherbrand} with~\eqref{eq:lambdaGD} and~\eqref{eq:lambdaGperiodic_signed}.\qedhere
\end{itemize}
\end{proof}
\begin{remark}
Observe that $\HHN{2j+3}{\Gamma}{\unit{L_\infty}}=0$ for all $j\geq 0$ because $\Gamma$ has strict cohomological dimension $2$. Nevertheless, $\lambdacohom{2j+3}=\lambdacohom{1}$ is non-zero, in general, showing that the $\bddsys{\DD}$-isomorphism between $\sysHH{i}{\unit{\pallino}}$ and $\limemb\bigl(\sysHH{i}{\unit{\pallino}}\bigr)$ of Corollary~\ref{cor:growth_cohom_units} does not hold for $i\geq 3$.
\end{remark}
Recall now~\cite[Theorem~3.14]{CapNuc20} which, applied to the dihedral extension $L_n/k$, says that 
\begin{equation}\label{eq:p-formula_dihedral}
\gorder{\clsyl{L_n}}=\frac{\gorder{\clsyl{K_n}}^2\gorder{\clsyl{F}}}{\gorder{\clsyl{k}}^2}\frac{\hh{0}{\splitG[n]}{\unit{L_n}}}{\hh{-1}{\splitG[n]}{\unit{L_n}}}.
\end{equation}

For every number field $M$, let $\hsyl{M}$ be the $p$-power of the order of the $p$-Sylow of its class group, so that
\[
\gorder{\clsyl{M}} = p^{\hsyl{M}}.
\]

Building upon Corollary~\ref{cor:growth_cohom_units_dihedral}, we find the following

\begin{theorem}\label{thm:iwasawa_formula} 
Let $K_\infty/k$ be a fake $\ZZ_p$-extension of dihedral type with Galois closure $L_\infty/k$ and normalizing quadratic extension $F/k$. There exist constants $\iwmu{\fk},\iwnu{\fk}\in\ZZ[\frac{1}{2}]$ and $\iwlambda{\fk}\in\ZZ$ such that
\begin{equation}\label{eq:yap_exp_fake}
\hsyl{K_n}=\iwmu{\fk}p^n+\iwlambda{\fk}n+\iwnu{\fk}\qquad\text{ for all }n\gg 0.
\end{equation}
Moreover, if we denote by $\iwmu{\iw},\iwlambda{\iw}$ the Iwasawa invariants of the $\ZZ_p$-extension $L_\infty/F$, then
\[
\iwmu{\fk}=\frac{\iwmu{\iw}}{2}\qquad\text{ and }\qquad\iwlambda{\fk}=\frac{\iwlambda{\iw}+\lambdacohom{-1}^+-\lambdacohom{0}^+}{2}\in\ZZ_{\geq 0}
\]
where $\lambdacohom{-1}^+$ and $\lambdacohom{0}^+$ are as in Corollary~\ref{cor:growth_cohom_units_dihedral}-\ref{cor:growth_cohom_units_dihedral_existence_global}.
\end{theorem}
\begin{proof} By taking $p$-adic valuation in~\eqref{eq:p-formula_dihedral} and rearranging terms we obtain
\begin{equation*}
2\hsyl{K_n}=\hsyl{L_n}+2\hsyl{k}-\hsyl{F}+v_p\bigl(\hh{-1}{\splitG[n]}{\unit{L_n}}\bigr)-v_p\bigr(\hh{0}{\splitG[n]}{\unit{L_n}}\bigl).
\end{equation*}
By Iwasawa's theorem quoted in the Introduction, the growth of $\hsyl{L_n}$ is controlled by the three invariants $\iwmu{\iw},\iwlambda{\iw},\iwnu{\iw}$. Formula~\eqref{eq:yap_exp_fake} follows from Corollary~\ref{cor:growth_cohom_units_dihedral}-\eqref{eq:lambdaD+}.

The only thing which is left to be checked is that $\iwlambda{\fk}$ is an integer: this follows by computing the difference $\hsyl{K_{n+1}}-\hsyl{K_n}$ for any $n$ big enough so that the~\eqref{eq:yap_exp_fake} holds. We find
\[
\hsyl{K_{n+1}}-\hsyl{K_n}=\iwmu{\fk}p^{n+1}+\iwlambda{\fk}(n+1)+\iwnu{\fk}-\iwmu{\fk}p^n-\iwlambda{\fk}n-\iwnu{\fk}=\iwmu{\fk}p^n(p-1)+\iwlambda{\fk}\in \ZZ
\]
and since $(p-1)$ is even we have $\iwmu{\fk}p^n(p-1)\in\ZZ$, whence $\iwlambda{\fk}\in\ZZ$.
\end{proof}

\begin{remark}\label{rmk:jaulent}
Under the hypothesis that $F/\QQ$ is an abelian extension such that $\Gal{F}{\QQ}$ is killed by $(p-1)$, and that $L_\infty/\QQ$ is Galois, the above result is a particular case of \cite[Théorème~3]{Jau81}.
\end{remark}

In \cite{CapNuc20}, explicit bounds for the ratio of class numbers are determined. Such bounds can be translated into bounds for $\iwlambda{\fk}$, as shown in the next corollary. For a number field $M$, set
\[
\defect{M}{p}=
\begin{cases}
0&\text{if $M$ does not contain any primitive $p$th root of unity}\\
1&\text{otherwise}
\end{cases}
\]   	

\begin{corollary}\label{cor:bounds} 
The following inequalities hold
\[
\frac{\iwlambda{\iw} + a}{2}\geq  \iwlambda{\fk} \geq \frac{\iwlambda{\iw} - b}{2}
\]
where $a=\rank_\ZZ(\unit{F}) + \defect{F}{p}+1$ and $b=\rank_\ZZ(\unit{k})+\defect{k}{p}$.
\end{corollary}
\begin{proof}
By~\cite[Corollary~3.15]{CapNuc20}, we have
\begin{equation}\label{eq:generalbounds}
-an\leq v_p\left(\frac{\gorder{\clsyl{L_n}}\cdot\gorder{\clsyl{k}}^2}{\gorder{\clsyl{F}}\cdot\gorder{\clsyl{K_n}}^2}\right)\leq bn
\end{equation}
where $v_p$ denotes the $p$-adic valuation. When $n\gg0$, using Theorem \ref{thm:iwasawa_formula} we can write the central term of the above chain of inequalities as $(\iwlambda{\iw} - 2\iwlambda{\fk})n + \kappa$ where $\kappa$ is a constant independent of $n$. We deduce that 
\[
(2\iwlambda{\fk} - \iwlambda{\iw} - a)n \leq \kappa \leq (b +2\iwlambda{\fk} - \iwlambda{\iw})n 
\]
and, since these inequalities hold for every $n\gg 0$, we must have
\[
(2\iwlambda{\fk} - \iwlambda{\iw} - a) \leq 0 \qquad\text{ and }\qquad (b +2\iwlambda{\fk} - \iwlambda{\iw}) \geq 0
\]
which is equivalent to the statement.
\end{proof}

We end this section by showing that the constant $\lambdacohom{-1}^+$ appearing in Corollary~\ref{cor:growth_cohom_units_dihedral}-\ref{cor:growth_cohom_units_dihedral_existence_global} and in Theorem~\ref{thm:iwasawa_formula}  has a more direct description in terms of arithmetic invariants of the extension $L_\infty/k$. In what follows, adopt the notation introduced in \S\ref{subsection:arithmetic_set-up} as well as those of Notation~\ref{notation:primesets}.

\begin{definition}\label{def:aboveT} 
For every $n$ such that every ramified prime in $L_\infty/L_n$ is totally ramified, denote by $\aboveT{L_n}\subseteq \clsyl{L_n}$ the subgroup  generated by the classes of the form
\begin{equation*}
\bigl[\algnm{G_n/G_n(\primeid)}\Primeid{n}\bigr]\in \clsyl{L_n}
\end{equation*}
for $\primeid$ running through the set $\Ram{L_n}{F}$. One sees immediately that $\aboveT{L_n}$ is in fact a $\splitG[n]$-submodule of $\clsyl{L_n}$, so that the plus and minus components $\aboveT{L_n}^\pm$ are well-defined.
\end{definition}

According to the above definition, the subgroups $\aboveT{L_n}$ are defined only for $n\gg 0$. This shall cause no harm to our arguments, because we will be interested in the behaviour of their order for $n\to +\infty$: in particular, they do not directly give rise to a normic system, although the relations of Definition~\ref{def:normic_system} are satisfied for $m\geq n\gg 0$.
\begin{proposition}\label{prop:arithmetic_interpretation_lambda1}
There exists $\iwlambda{\aboveT{}},\iwlambda[\pm]{\aboveT{}}\in \ZZ_{\geq 0}$ such that $\gorder{\aboveT{L_n}}\bddeq p^{\iwlambda{\aboveT{}}n}$ and $\gorder{\aboveT{L_n}^\pm}\bddeq p^{\iwlambda[\pm]{\aboveT{}}n}$. They are related to the invariants obtained in Corollary~\ref{cor:growth_cohom_units_dihedral}-\ref{cor:growth_cohom_units_dihedral_existence_global} as follows:
\begin{equation}\label{eq:relation_lambdaT_lambdaunit}
\lambdacohom{1}=\oram{L_\infty}{F}-\iwlambda{\aboveT{}},\quad \lambdacohom{1}^+=\oram{L_\infty}{k}-\iwlambda[+]{\aboveT{}}\quad\text{ and }\quad \lambdacohom{-1}^+=\oram{L_\infty}{k}_s-\iwlambda[-]{\aboveT{}}.
\end{equation}
\end{proposition}
\begin{proof}
Let $\calY=(Y_n,\dsu_{n,m},\dsd_{m,n})$ be the cokernel of the first map in~\eqref{eq:long_cohom_DS}, so that, by Theorem \ref{thm:i_cinque_dell'apocalisse}, there is a short exact sequence 
\begin{equation}\label{eq:short_cohom_quot}
0\longdashrightarrow\sysHH{1}{\unit{\pallino}}\longdashrightarrow\sysHH{1}{\idunit{\pallino}}\longdashrightarrow\calY\longdashrightarrow 0
\end{equation}
in $\quotsys[\cofg]{\DD}$. The first two terms satisfy $\sysHH{1}{\unit{\pallino}}\bddcong\limemb\bigl(\sysHH{1}{\unit{\pallino}}\bigr)$ and $\sysHH{1}{\idunit{\pallino}}\bddcong\limemb\bigl(\sysHH{1}{\idunit{\pallino}}\bigr)$, by Corollary~\ref{cor:growth_cohom_units} and by Corollary~\ref{cor:indunits_12}, respectively. Since the functor $\limemb$ is exact, we obtain $\calY\bddcong\limemb\bigl(\calY\bigr)$. Thanks to Corollary~\ref{cor:ha_vinto_caputo}, there exists $\iwlambda{\calY}$ such that $\gorder{Y_n}\bddeq p^{\iwlambda{\calY}n}$ and Corollary~\ref{cor:indunits_12} implies 
\begin{equation}\label{eq:lambda_units_X}
\lambdacohom{1}=\oram{L_\infty}{F}-\iwlambda{\calY}.
\end{equation}
Moreover, given any $\DD$-double system $\calX$ satisfying $\calX\bddcong\limemb(\calX)$, we deduce $\calX^\pm\bddcong[\bddsys{\Gamma}]\limemb(\calX^\pm)$ because the functor $\limemb$ commutes with the action of $\splitG$: hence, $\eqref{eq:lambda_units_X}$ admits the signed versions
\begin{equation}\label{eq:lambda_units_X_signed}
\lambdacohom{1}^+=\oram{L_\infty}{k}-\iwlambda[+]{\calY}\qquad\text{ and }\qquad\lambdacohom{1}^-=\oram{L_\infty}{k}_s-\iwlambda[-]{\calY}.
\end{equation}
To conclude the proof we show that
\begin{equation}\label{eq:piy}
\gorder{\aboveT{L_m}^\ast}\bddeq \gorder{Y_n^\ast}
\end{equation}
for all choices of a sign $\ast\in\{\emptyset,\pm\}$: one obtains the sought-for constant $\iwlambda{\aboveT{}}$ by setting $\iwlambda[\ast]{\aboveT{}}=\iwlambda[\ast]{\calY}$. We give the proof of~\eqref{eq:piy} for $\ast = \emptyset$; the other cases follow since all 
morphisms considered in the proof are $\splitG$-equivariant. 

For $m\geq n\geq 0$, let $\cch{0}_{m,n}\colon\HHN{0}{G_{m,n}}{\clsyl{L_m}}\to\HH{1}{G_{m,n}}{\q{L_n}}$ and $\cchT{0}_{m,n}\colon\HH{0}{G_{m,n}}{\clsyl{L_m}}\to \HH{1}{G_{m,n}}{\q{L_m}}$, respectively, be the connecting homomorphisms induced by the bottom row in~\eqref{diag:magic} for regular (\rsp Tate) cohomology. For brevity, write $\cch{0}_{n}$ (\rsp $\cchT{0}_n$) to denote~$\cch{0}_{n,0}$ (\rsp $\cchT{0}_{n,0}$). Observe that, for all $m\geq n\geq 0$, the equality
\[
\cchT{0}_{m,n}(\overline{\aboveT{L_m}})=\cch{0}_{m,n}(\aboveT{L_m})
\]
holds, where $\overline{\aboveT{L_m}}$ denotes the image of $\aboveT{L_m}\subseteq\HHN{0}{G_{m,n}}{\clsyl{L_m}}$ in $\HH{0}{G_{m,n}}{\clsyl{L_m}}$. In particular,~\cite[Proposition~2.1]{Nuc10}, yields an isomorphism
\begin{equation}\label{eq:iso_from_cambridge}
Y_n\cong\cchT{0}_n(\overline{\aboveT{L_n}})=\cch{0}_n(\aboveT{L_n})\qquad\text{ for all }n\geq 0.
\end{equation}
Therefore, to prove \eqref{eq:piy} we are left to show that $\kernel\bigl(\cch{0}_n\vert_{\aboveT{L_n}}\bigr)$ has constant order for $n$ large enough.
Let $\indextr[\empty]$ be the first index such that all ramified primes in $L_\infty/L_{\indextr[\empty]}$ are totally ramified and consider, for $n\geq \indextr[\empty]$, the commutative diagram
\begin{equation}\label{diag:delta_asymptotic}
\begin{split}\xymatrix@C=2.89em{
	\HH{0}{G_{n,\indextr[\empty]}}{\idunit{L_n}}\ar@{->}[1,0]\ar@{->}[1,1]^{\alpha}\\
	\HH{0}{G_{n,\indextr[\empty]}}{\q{L_n}}\ar@{->}[0,1]&
\HH{0}{G_{n,\indextr[\empty]}}{\idclg{L_n}}\ar@{->}[0,1]&\HH{0}{G_{n,\indextr[\empty]}}{\clsyl{L_n}}\ar@{->}[0,1]^{\cchT{0}_{n,\indextr[\empty]}}&\HH{1}{G_{n,\indextr[\empty]}}{\q{L_n}}\ar@{->}[0,1]&0
}\end{split}
\end{equation}
The compatibility of local and global class field theory shows that $\alpha$ is surjective: indeed, its image is the product of all inertia subgroups in $G_{n,\indextr[\empty]}$, and this is the whole $G_{n,\indextr[\empty]}$ since $L_n/L_{\indextr[\empty]}$ is totally ramified at all primes where it ramifies. It follows that $\cchT{0}_{n,\indextr[\empty]}$ is injective and therefore, for all $n\geq \indextr[\empty]$,
\begin{equation}\label{eq:kernel_algnorm}
\kernel(\cch{0}_{n,\indextr[\empty]})=\algnm{G_{n,\indextr[\empty]}}(\clsyl{L_n})\qquad\text{ and }\qquad\kernel(\cch{0}_{n,\indextr[\empty]})\cap\aboveT{L_n}=\algnm{G_{n,\indextr[\empty]}}(\clsyl{L_n})\cap\aboveT{L_n}.
\end{equation}
Since $L_n/L_{\indextr[\empty]}$ is totally ramified, the arithmetic norm is surjective on ideal class groups and we obtain
\begin{equation}\label{eq:norm_eventually_ext}
\algnm{G_{n,\indextr[\empty]}}(\clsyl{L_n})=\arext{L_n}{L_{\indextr[\empty]}}\circ \arnm{L_n}{L_{\indextr[\empty]}}(\clsyl{L_n})
=\arext{L_n}{L_{\indextr[\empty]}}(\clsyl{L_{\indextr[\empty]}}).
\end{equation}
This already shows that the order $\gorder{\algnm{G_{n,\indextr[\empty]}}(\clsyl{L_n})}$ is bounded independently of $n$; the same holds, \emph{a fortiori} for the order $\gorder{\algnm{G_{n,\indextr[\empty]}}(\clsyl{L_n})\cap\aboveT{L_n}}$. Moreover, for $m\geq n\geq \indextr[\empty]$, \eqref{eq:norm_eventually_ext} yields identifications
\[
\arext{L_m}{L_n}\bigl(\algnm{G_{n,\indextr[\empty]}}(\clsyl{L_n})\bigr)=\arext{L_m}{L_n}\bigl(\arext{L_n}{L_{\indextr[\empty]}}(\clsyl{L_{\indextr[\empty]}})\bigr)=\arext{L_m}{L_{\indextr[\empty]}}(\clsyl{L_{\indextr[\empty]}})=\algnm{G_{m,\indextr[\empty]}}(\clsyl{L_m}).
\]
It follows that, for all $m\geq n\gg\indextr[\empty]$, the maps $\arext{L_m}{L_n}$ are surjections between groups of bounded orders: they are therefore isomorphisms. In particular, through~\eqref{eq:kernel_algnorm}, they induce injections $\arext{L_m}{L_n}\colon\kernel(\cch{0}_{n,\indextr[\empty]})\cap\aboveT{L_n}\hookrightarrow\kernel(\cch{0}_{m,\indextr[\empty]})\cap\aboveT{L_m}$;
again, since the orders $\gorder{\kernel(\cch{0}_{n,\indextr[\empty]})\cap\aboveT{L_n}}=\gorder{\algnm{G_{n,\indextr[\empty]}}(\clsyl{L_n})\cap\aboveT{L_n}}$ are bounded, this shows that the extension maps are actually isomorphisms
\begin{equation}\label{eq:iso_kernel_mn}
\arext{L_m}{L_n}\colon\kernel(\cch{0}_{n,\indextr[\empty]})\cap\aboveT{L_n}=\kernel\bigl(\cch{0}_{n,\indextr[\empty]}\vert_{\aboveT{L_n}}\bigr)\overset{\cong}{\longrightarrow}\kernel(\cch{0}_{m,\indextr[\empty]})\cap\aboveT{L_m}=\kernel\bigl(\cch{0}_{m,\indextr[\empty]}\vert_{\aboveT{L_m}}\bigr)
\end{equation}
 for all $m\geq n\gg\indextr[\empty]$. Observe now that, for all $n\geq \indextr[\empty]$, there is a commutative triangle
\begin{equation*}
\xymatrix@C=7em{
\aboveT{L_n}\ar@{->}[1,0]_{\cch{0}_n}\ar@{->}[1,1]^{\cch{0}_{n,\indextr[\empty]}}\\
\cch{0}_n(\aboveT{L_n})\ar@{->}[0,1]^{\res_{G_{n,\indextr[\empty]}}}&\cch{0}_{n,\indextr[\empty]}(\aboveT{L_n})
}
\end{equation*}
It induces an inclusion $\kernel\bigl(\cch{0}_n\vert_{\aboveT{L_n}}\bigr)\subseteq\kernel\bigl(\cch{0}_{n,\indextr[\empty]}\vert_{\aboveT{L_n}}\bigr)$. This implies that $\gorder{\kernel\bigl(\cch{0}_n\vert_{\aboveT{L_n}}\bigr)}$ is bounded independently of $n$, a bound being the order $\gorder{\kernel\bigl(\cch{0}_{n,\indextr[\empty]}\vert_{\aboveT{L_n}}\bigr)}$, which is independent of $n\gg \indextr[\empty]$ by~\eqref{eq:iso_kernel_mn}. Moreover, restricting the isomorphism in~\eqref{eq:iso_kernel_mn} to $\kernel\bigl(\cch{0}_n\vert_{\aboveT{L_n}}\bigr)$ induces inclusions, for $m\geq n\gg \indextr[\empty]$,
\[
\arext{L_m}{L_n}\colon\kernel\bigl(\cch{0}_n\vert_{\aboveT{L_n}}\bigr)\longhookrightarrow\kernel\bigl(\cch{0}_m\vert_{\aboveT{L_m}}\bigr).
\]
Once more, these are inclusions among groups of bounded order, and are therefore isomorphisms, showing that $\gorder{\kernel\bigl(\cch{0}_n\vert_{\aboveT{L_n}}\bigr)}$ is constant for $n$ large enough.
\end{proof}

\begin{remark}\label{rmk:ker_delta_bdd}
For later reference, we extract from the above proof the fact that $\kernel(\cchT{0}_n)$ is bounded independently of $n$ (and actually eventually trivial, although we will not need this). Indeed, we have shown that $\cchT{0}_{n,\indextr[\empty]}$ is injective for all $n\geq \indextr[\empty]$, and the commutative triangle
\begin{equation*}
\xymatrix@C=7em{
\HH{0}{G_n}{\clsyl{L_n}}\ar@{->}[1,0]_{\cchT{0}_n}\ar@{->}[1,1]^{\cchT{0}_{n,\indextr[\empty]}}\\
\HH{1}{G_n}{\q{L_n}}\ar@{->}[0,1]^{\res_{G_{n,\indextr[\empty]}}}&\HH{1}{G_{n,\indextr[\empty]}}{\q{L_n}}
}
\end{equation*}
yields an injection $\kernel(\cchT{0}_n)\subseteq \kernel(\cchT{0}_{n,\indextr[\empty]})$. In particular, $\kernel(\cchT{0}_{n,\indextr[\empty]})=0$ for all $n\geq \indextr[\empty]$.
\end{remark}

\subsection{\texorpdfstring{Formul\ae\ for $\lambdacohom{0}^+$ in special cases}{Formul\ae\ for the zero-th invariant in special cases}}\label{subsec:special_cases}
In this section, we fix a fake $\ZZ_p$-extension of dihedral type $K_\infty/k$ and we place ourselves in the setting introduced in~\S\ref{subsec:fake_def}. The invariant $\iwlambda{\fk}$ of Theorem~\ref{thm:iwasawa_formula} is explicitly related to the constants $\lambdacohom{-1}^+$  and $\lambdacohom{0}^+$, and we have provided an arithmetic description of the former in Proposition~\ref{prop:arithmetic_interpretation_lambda1}. In this section, we show that in some cases also the latter can be interpreted in more explicit arithmetic terms. In the next proposition, write $\realemb{k}$ and $\cpxemb{k}$, respectively, for the number of real embeddings and of pairs of conjugate complex embeddings of $k$, respectively.
\begin{proposition}\label{prop:formula_con_hp}
Suppose that one of the following conditions holds:
\begin{enumerate}[label=\textup{(}\textit{\alph*}\textup{)}]
\item $F/k$ is totally imaginary and $k$ is totally real, so that $F$ is a CM field with $k$ as totally real subfield;\label{list:thm_formula:CM}
\item $k$ is either $\QQ$ or an imaginary quadratic field;\label{list:thm_formula:quad}
\item $p$ is totally split in $F/\QQ$, all primes of $F$ above $p$ ramify in $L_\infty/F$ and the Leopoldt conjecture holds for $F$ \textup{(}and $p$\textup{)}.\label{list:thm_formula:Ass}
\end{enumerate}
With the same notation introduced in Notation~\ref{notation:primesets} and in Proposition~\ref{prop:arithmetic_interpretation_lambda1}, we have
\begin{equation*}
\lambdacohom{0}^+=
\begin{cases}
	\oram{L_\infty}{F}-\iwlambda{\aboveT{}}-1&\text{ in case \ref{list:thm_formula:CM}}\\
	0&\text{ in case \ref{list:thm_formula:quad}}\\
	\realemb{k}+\cpxemb{k}-1&\text{ in case \ref{list:thm_formula:Ass}}\\
\end{cases}
\end{equation*}
\end{proposition}
\begin{proof}
\begin{itemize}[wide, labelindent=0pt]
\item[\ref{list:thm_formula:CM}] By \cite[Theorem~4.12]{Was97} we know (recall that all abelian groups under consideration are $\ZZ_p$-modules) that $\unit{F}=\unit{k}=(\unit{F})^+$: it follows that $\HH{0}{G_n}{\unit{L_n}}^-=0$, so that $\lambdacohom{0}^-=0$. Corollary~\ref{cor:growth_cohom_units_dihedral}-\eqref{eq:lambdaGD} yields $\lambdacohom{0}^+=\lambdacohom{0}$ and~Corollary~\ref{cor:growth_cohom_units_dihedral}-\eqref{eq:lambdaherbrand} finally gives $\lambdacohom{0}^+=\lambdacohom{1}-1$. The result follows from Proposition~\ref{prop:arithmetic_interpretation_lambda1}.
\item[\ref{list:thm_formula:quad}] When $k$ is $\QQ$ or imaginary quadratic the unit group $\unit{k}$ is finite, so $\lambdacohom{0}^+=0$.
\item[\ref{list:thm_formula:Ass}] To prove the statement, we first define an ancillary normic system. For every $n\geq 0$, let
\[
\semiloc{L_n}=\prod_{\begin{smallmatrix}\otherprime\subseteq\rint{L_n}\\\otherprime\mid p\end{smallmatrix}} \unit{\otherprime}
\]
be the $\ZZ_p$-module of semi-local $p$-adic units of $L_n$. Clearly, this gives rise to a $\DD$-normic system $\semiloc{\empty}$ whose $n$th component is $\semiloc{L_n}$. The diagonal embedding $\unit{L_n}\hookrightarrow \semiloc{L_n}$ induces a morphism of normic systems $\unit{\pallino}\hookrightarrow \semiloc{\pallino}$ which gives rise, through~\eqref{diag:f*_functorial}, to a morphism $\sysHH{0}{\unit{\pallino}}\overset{u}{\longdashrightarrow} \sysHH{0}{\semiloc{\pallino}}$ in $\quotsys[\cofg]{\DD}$. Consider now the $\DD$-double system $\unit{F}/(\unit{F})^{\pallino[\expp]}$ introduced in~\eqref{def:double_sys_OF/p} together with the semi-local counterpart
\[
\semiloc{F}/(\semiloc{F})^{\pallino[\expp]}=\left(\semiloc{F}/(\semiloc{F})^{p^n},\dsu_{n,m}\colon x\mapsto x^{p^{m-n}}\mod{(\semiloc{F})^{p^m}},\dsd_{m,n}\colon x\mapsto x\mod{(\semiloc{F})^{p^n}}\right)_{m\geq n\geq 0};
\]
the diagonal embedding induces a morphism $\unit{F}/(\unit{F})^{\pallino[\expp]}\overset{v}{\longdashrightarrow} \semiloc{F}/(\semiloc{F})^{\pallino[\expp]}$ in $\quotsys[\cofg]{\DD}$.
The inclusions $(\unit{F})^{p^n}\subseteq \arnm{L_n}{F}\unit{L_n}$ (\rsp $(\semiloc{F})^{p^n}\subseteq\arnm{L_n}{F}\semiloc{L_n}$) give rise to a commutative diagram in $\quotsys[\cofg]{\DD}$
\begin{equation}\begin{aligned}\label{diag:miki}
\xymatrix{
&1\ar@{-->}[1,0]&1\ar@{-->}[1,0]&1\ar@{-->}[1,0]&\\
1\ar@{-->}[0,1]&\kernel w\ar@{-->}[1,0]\ar@{-->}[1,0]\ar@{-->}[0,1]&\kernel v\ar@{-->}[1,0]\ar@{-->}[0,1]&\kernel u\ar@{-->}[1,0]\\
1\ar@{-->}[0,1]&\sysnorm{\unit{\pallino}}\ar@{-->}[0,1]\ar@{-->}[1,0]^w&\unit{F}/(\unit{F})^{\pallino[\expp]}\ar@{-->}[0,1]\ar@{-->}[1,0]^v&\sysHH{0}{\unit{\pallino}}\ar@{-->}[0,1]\ar@{-->}[1,0]^u&1\\
1\ar@{-->}[0,1]&\sysnorm{\semiloc{\pallino}}\ar@{-->}[0,1]
&\semiloc{F}/(\semiloc{F})^{\pallino[\expp]}\ar@{-->}[0,1]
&\sysHH{0}{\semiloc{\pallino}}\ar@{-->}[0,1]
&1
}\end{aligned}\end{equation}
where $\sysnorm{\unit{\pallino}}$ and $\sysnorm{\semiloc{\pallino}}$ denote, respectively, the double system whose components are $\arnm{L_n}{F}\unit{L_n}/(\unit{F})^{p^n}$ and $\arnm{L_n}{F}\semiloc{L_n}/(\semiloc{F})^{p^n}$. We now make the following claims:
\begin{enumerate}[before=\medskip, label=\textup{Claim \Alph*)}, align=left, leftmargin=3.5\parindent, labelwidth=35pt, itemsep=.5\bigskipamount,   after=\medskip]
\item When $p$ is totally split in $F/\QQ$ and all primes above $p$ ramify in $L_\infty/F$, the double system $\sysnorm{\semiloc{\pallino}}$ belongs to $\bddsys{\DD}$;\label{prop:formula_con_hp:claim:split}
\item Assuming the Leopoldt conjecture for $F$ and $p$, the morphism $v$ is a $\bddsys{\DD}$-monomorphism.\label{prop:formula_con_hp:claim:Leopoldt} 
\end{enumerate}
Assuming the claims, let us derive point~\ref{list:thm_formula:Ass}. \ref{prop:formula_con_hp:claim:Leopoldt} implies that $\kernel w = 0$, so \ref{prop:formula_con_hp:claim:split}  yields $\sysnorm{\unit{\pallino}}\bddcong 0$: we get an isomorphism
\begin{equation}\label{seq:miki_and_split}\
\unit{F}/(\unit{F})^{\pallino[\expp]}\bddcong\sysHH{0}{\unit{\pallino}}
\end{equation}
and, taking plus parts,
\begin{equation}\label{seq:miki_global}
\unit{k}/(\unit{k})^{\pallino[\expp]}\bddcong[\bddsys{\Gamma}]\sysHH{0}{\unit{\pallino}}^+.
\end{equation}
The structure theorem of global units gives $\unit{k}/(\unit{k})^{\pallino[\expp]}\bddcong[\bddsys{\Gamma}]\calZ^{\realemb{k}+\cpxemb{k}-1}$ and~\ref{list:thm_formula:Ass} follows, by combining Example~\ref{example:Q/Z} with Corollary~\ref{cor:ha_vinto_caputo}. We are therefore left with the proof of the claims. 
\end{itemize}
Concerning~\ref{prop:formula_con_hp:claim:split}, fix $n\gg 0$ such that all primes above $p$ are totally ramified in $L_\infty/L_n$: for every $\primeid\mid p$ in $F$ and every $\otherprime_n\mid\primeid$ in $L_n$, the extension $L_{\otherprime_n}/F_\primeid$ is a cyclic extension of $\QQ_p$, with ramification degree $p^{n-\logindeg[\primeid]}$ (where $\logindeg[\primeid]$ is as in Notation~\ref{notation:ram_and_ind}), and local class field theory implies that $\arnm{\otherprime_n}{\primeid}\unit{\otherprime_n}=(\ZZ_p^\times)^{p^{n-\logindeg[\primeid]}}$. It follows that
\[
\arnm{L_n}{F}\semiloc{L_n}\cong\prod_{\begin{smallmatrix}\otherprime_n\subseteq\rint{L_n}\\\otherprime\mid p\end{smallmatrix}}(\ZZ_p^\times)^{p^{n-\logindeg[\primeid]}}
\]
and thus the double system $\sysnorm{\semiloc{\pallino}}$ coincides with the ``constant'' double system
\[
\prod_{\begin{smallmatrix}\otherprime_n\subseteq\rint{L_n}\\\otherprime\mid p\end{smallmatrix}}\bigl(\ZZ/p^{\logindeg[\primeid]},\id,\phantom{\cdot}p^{m-n})\bigr)_{m\geq n\geq 0}
\]
which lies in $\bddsys{\DD}$ by Proposition~\ref{prop:trivialbdd}.

Concerning~\ref{prop:formula_con_hp:claim:Leopoldt},~\cite[Lemma~2]{Mik87} shows that the Leopoldt conjecture is equivalent to the existence of $\nu\geq 0$ such that
\begin{equation*}
\unit{F}\cap (\semiloc{F})^{p^{n}}=\bigl(\unit{F}\cap (\semiloc{F})^{p^{n-1}}\bigr)^p\qquad\text{ for all }n> \nu
\end{equation*}
or, equivalently, such that
\begin{equation}\label{eq:miki_lemma2}
\unit{F}\cap (\semiloc{F})^{p^{n}}=\bigl(\unit{F}\cap (\semiloc{F})^{p^{\nu}}\bigr)^{p^{n-\nu}}\qquad\text{ for all }n> \nu.
\end{equation}
On the other hand, the $n$th component of $\kernel{v}$ is $\bigl(\unit{F}\cap(\semiloc{F})^{p^n}\bigr)/(\unit{F})^{p^n}$ which can be rewritten, using \eqref{eq:miki_lemma2}, as
\[
\bigl(\kernel{v}\bigr)_n=\bigl(\unit{F}\cap (\semiloc{F})^{p^{\nu}}\bigr)^{p^{n-\nu}}/\bigl((\unit{F})^{p^{\nu}}\bigr)^{p^{n-\nu}}\qquad\text{ for all }n> \nu.
\]
Set $M=\unit{F}\cap (\semiloc{F})^{p^{\nu}}$ and $N=(\unit{F})^{p^{\nu}}$: then $N\subseteq M$ are two finitely generated $\ZZ_p$-modules such that $M/N=\bigl(\kernel{v}\bigr)_0$ is finite. It follows that they have the same rank, and this shows that
\[
\bigl(\kernel{v}\bigr)_n=M^{p^{n-\nu}}/N^{p^{n-\nu}}\qquad\text{ has order bounded independently of }n\geq 0.
\]
To conclude the proof, we show that the $\DD$-double system $\kernel{v}$ satisfies condition~\ref{point:lemma_eqB:cap=nuc}, namely that $\nuczero{\kernel{v}_n}=\capzero{\kernel{v}_n}$ for all $n\gg 0$. We need to compute the direct and inverse limits of the $\DD$-double system
\[
\bigl(M^{p^n}/N^{p^n},\dsu_{n,m}\colon x\longmapsto x^{p^{m-n}}\mod N^{p^m},\dsd_{m,n}\colon x\longmapsto x\mod N^{p^n}\bigl).
\]
Since both $M$ and $N$ are finitely generated $\ZZ_p$-modules, the inverse limit functor is exact and
\[
\varprojlim \bigl(\kernel{v}\bigr)_n=\varprojlim \bigl(M^{p^n}/N^{p^n}\bigr)=\varprojlim M^{p^n}/\varprojlim N^{p^n}=\bigcap_{n\geq 0}M^{p^n}/\bigcap_{n\geq 0}N^{p^n}=0.
\]
In particular,
\begin{equation}\label{eq:nuczero_kerv}
\nuczero{\kernel{v}_n}=0\qquad\text {for all }n\geq 0.
\end{equation}
On the other hand, the torsion submodule of $M$ is finite, say $\tors[p^\infty]{M}=\tors[p^t]{M}$. It follows that the morphism $\dsu_{n}\colon M^{p^n}/N^{p^n}\to \varinjlim M^{p^n}/N^{p^n}$ is injective for $n\geq t$: to see this, suppose $\bar{x}\in\kernel (\dsu_n)$ is represented by some $x\in M^{p^n}$ and let $x_0\in M$ be such that $x=x_0^{p^n}$. To say that $\dsu_n(\bar{x})=0$ means that there exists $m\geq n$ such that $\dsu_{n,m}(\bar{x})=0$ or, equivalently, $x^{p^{m-n}}\in N^{p^m}$. Hence, there exists $y\in N$ such that $x_0^{p^m}=x^{p^{m-n}}=y^{p^m}$: therefore $x_0=y\cdot\zeta$ for some $\zeta\in\tors[p^\infty]{M}=\tors[p^t]{M}$. Since $n\geq t$, this yields $x=x_0^{p^n}=y^{p^n}\in N^{p^n}$, so that $\bar{x}=0$. In particular,
\begin{equation}\label{eq:capzero_kerv}
\capzero{\kernel{v}_n}=0\qquad\text {for all }n\gg 0
\end{equation}
and combining~\eqref{eq:nuczero_kerv} with~\eqref{eq:capzero_kerv} establishes condition~\ref{point:lemma_eqB:cap=nuc}, concluding the proof of~\ref{prop:formula_con_hp:claim:Leopoldt}.
\end{proof}
Before proceeding, we need a small digression concerning the Gross conjecture. Retain the same notation as in the rest of the paper but suppose that $F$ is a CM number field\footnote{The Gross conjecture can actually be stated in greater generality, but we will not need this.} and write $\Lcyc$ for its cyclotomic extension. Denote by $M_0$ the maximal $p$-extension of $\Lcyc$ which is unramified everywhere and abelian over~$F$, and by $M_0'\subseteq M_0$ the maximal $p$-extension of $\Lcyc$ which is totally split at every prime above $p$ in $\Lcyc$ and abelian over~$F$. Denote by $k=F^+$ the maximal totally real subfield of $F$ and by $c$ the non-trivial element of $\Gal{F}{k}$. It acts upon $\Gal{M'_0}{F}$, decomposing it as a direct sum
\[
\Gal{M'_0}{F}=\Gal{M'_0}{F}^{c=1}\oplus\Gal{M'_0}{F}^{c=-1}.
\]
The Gross conjecture is

\begin{conjecture}[Gross]  
The group $\Gal{M'_0}{F}^{c=-1}$ is finite.
\end{conjecture}

The original formulation of Gross conjecture as in \cite[Conjecture~1.15 and~(1.21)]{Gro81} is proven to be equivalent to the above formulation\footnote{We have made a small change of notation with respect to \cite{FedGro81}: in particular, the field that we call $M_0$ is called $L_0$ in \loccit and the superscript $(-)^{c=-1}$ is simply denoted $(-)^-$ in \loccit These changes were forced by the fact that $L_0=F$ already has a meaning in our work, and the $(-)^-$-notation is normally used for $\splitG[n]$-modules rather than for $\splitG\times\Gal{(\Lcyc)_n}{L_0}$-modules.} in \cite[(6.6) and p.~457]{FedGro81}.

\begin{lemma}\label{lemma:Gross}
Let $F$ be a CM field which is Galois over $\QQ$ and which satisfies both the Leopoldt and the Gross conjectures. Suppose that $L_\infty/F$ is a $\ZZ_p$-extension which is Galois over $\QQ$ and is different from the cyclotomic $\ZZ_p$-extension $\Lcyc/F$: write $\Ldue$ for the compositum $\Ldue=L_\infty\Lcyc$. Assume that a prime of $L_\infty$ above $p$ is unramified in $\Ldue/L_\infty$. Then every prime $\Otherprime'$ of $L_\infty$ above $p$ is unramified in $\Ldue/L_\infty$ and $\Frob{\Otherprime'}{\Ldue/L_\infty}$ has infinite order in the group $\Gal{\Ldue}{L_\infty}$.
\end{lemma}
\begin{proof}
Let $\Otherprime$ be a prime of $L_\infty$ above $p$ that is unramified in $\Ldue/L_\infty$, and let $\widetilde{\Otherprime}$ be a prime of $\Ldue$ above $\Otherprime$; set $\Otherprime[\cyc]=\widetilde{\Otherprime}\cap \Lcyc$ for its restriction to $\Lcyc$ and $\otherprime=\widetilde{\Otherprime}\cap F$ for its restriction to $F$. Let $\Otherprime[\cyc]=\Otherprime[\cyc]^{(1)},\hdots,\Otherprime[\cyc]^{(s)}$ be all the primes in $\Lcyc$ above $p$. The diagram of fields and primes is as follows:
\[\xymatrix{
\overline{\Otherprime}\ar@{-}[0,2]&&\Ldue\ar@{-}[1,1]\ar@{-}[1,-1]\\
\Otherprime\ar@{-}[0,1]&L_\infty\ar@{-}[1,1]&&\Lcyc\ar@{-}[1,-1]&\{\Otherprime[\cyc]=\Otherprime[\cyc]^{(1)},\hdots,\Otherprime[\cyc]^{(s)}\}\ar@{-}[0,-1]\\
\otherprime\ar@{-}[0,2]&&F\ar@{-}[1,0]\\
p\ar@{-}[0,2]&&\QQ
}\]
The assumption that $L_\infty$ is Galois over $\QQ$ implies that $\Otherprime$ is infinitely ramified because $L_\infty/F$ must eventually ramify at one, and hence at all primes above $p$ in $F$. It also implies that $\Ldue/\QQ$ is Galois and therefore each prime $\Otherprime'$ of $L_\infty$ dividing $p$ is unramified in $\Ldue/L_\infty$, establishing the first part of the lemma. For the same reason, if we can prove that $\Frob{\Otherprime,\Ldue/L_\infty}$ has infinite order, then the same is true for all other $\Frob{\Otherprime'}{\Ldue/L_\infty}$. A similar argument shows that the splitting behaviour of $\Otherprime[\cyc]^{(i)}$ in $\Ldue/\Lcyc$ is the same as that of $\Otherprime[\cyc]$, for all $1\leq i\leq s$.

Suppose by contradiction that $\Frob{\Otherprime}{\Ldue/L_\infty}$ has finite order or, equivalently, that it is trivial in the group $\Gal{\Ldue}{L_\infty}\cong\ZZ_p$: in particular, $(\Lcyc)_{\Otherprime[\cyc]}(L_\infty)_{\Otherprime}=(L_\infty)_{\Otherprime}$ or equivalently $(\Lcyc)_{\Otherprime[\cyc]}\subseteq (L_\infty)_{\Otherprime}$. This means in particular that $(\Lcyc)_{\Otherprime[\cyc]}=(L_\infty)_{\Otherprime}$ since  $\Gal{(L_\infty)_{\Otherprime}}{F_\otherprime}\cong \ZZ_p$ does not have any non-trivial quotient isomorphic to $\ZZ_p$. Then the extension $\Ldue/\Lcyc$ is totally split at $\Otherprime[\cyc]$ and hence, by the above considerations, it is totally split at each of the $\Otherprime[\cyc]^{(i)}$. This implies that $\Ldue\subseteq M'_0$ (we follow the same notation as in the above formulation of the Gross conjecture). By hypothesis, the totally real subfield $k$ of $F$ satisfies the Leopoldt conjecture and hence admits a unique $\ZZ_p$-extension, so that $\Gal{\Ldue}{\Lcyc}=\Gal{\Ldue}{\Lcyc}^{c=-1}$. The inclusion $\Ldue\subseteq M'_0$ then implies that ${\Gal{M'_0}{\Lcyc}}^{c=-1}$ is infinite, because $\Gal{\Ldue}{\Lcyc}\cong\ZZ_p$, violating the Gross conjecture. This provides the required contradiction and establishes the lemma.
\end{proof}

In the remainder of this section, we focus on some special cases to relate the invariant $\iwlambda{\fk}$ appearing in Theorem~\ref{thm:iwasawa_formula} to the classical invariants $\iwlambda{\iw}$.

\begin{corollary}\label{cor:formula_quad} 
Suppose that $k=\QQ$, so that $F$ is an imaginary quadratic field and $L_\infty/F$ is the anticyclotomic $\ZZ_p$-extension of $F$. Then, with notation as in Proposition~\ref{prop:arithmetic_interpretation_lambda1},
\[
\lambdacohom{-1}^+= \iwlambda{\aboveT{}}=\iwlambda[+]{\aboveT{}} = \oram{L_\infty}{\QQ}_s\in\{0,1\}\qquad \text{ and }\qquad \lambdacohom{0}^+ =\iwlambda[-]{\aboveT{}} = 0.
\]
In particular,
\[
\iwlambda{\fk}=
\begin{cases}
\frac{\iwlambda{\iw}+1}{2}&\text{ if $p$ splits in $F/\QQ$}\\&\\
\frac{\iwlambda{\iw}}{2}&\text{ if $p$ does not split in $F/\QQ$}
\end{cases}
\]
and therefore $\iwlambda{\iw}$ is odd if $p$ splits in $F$ and it is even if $p$ does not split.
\end{corollary}
\begin{proof} Once the formul\ae\ for $\lambdacohom{0}$ and $\lambdacohom{-1}$ are established, the formula for $\iwlambda{\fk}$ and the claim about its parity follow readily from Theorem~\ref{thm:iwasawa_formula}.

In the setting where $k=\QQ$ and $F$ is imaginary quadratic, both conditions~\ref{list:thm_formula:CM} and~\ref{list:thm_formula:quad} of Proposition~\ref{prop:formula_con_hp} are satisfied. The equality $\lambdacohom{0}^+ = 0$ is proved Proposition~\ref{prop:formula_con_hp}-\ref{list:thm_formula:quad} and Proposition~\ref{prop:formula_con_hp}-\ref{list:thm_formula:CM} implies that
\begin{equation}\label{eq:oram_quad}
\oram{L_\infty}{F}=\iwlambda{\aboveT{}}+1.
\end{equation}
Assume first that $p$ does not split in $F/\QQ$, so $\oram{L_\infty}{F}=1$. Then~\eqref{eq:oram_quad} reads $1=\iwlambda{\aboveT{}}+1$, so $\iwlambda{\aboveT{}}=\iwlambda[+]{\aboveT{}}=\iwlambda[-]{\aboveT{}}=0=\oram{L_\infty}{\QQ}_s$ and Proposition~\ref{prop:arithmetic_interpretation_lambda1} gives  $\lambdacohom{-1}^+=\oram{L_\infty}{\QQ}_s-\iwlambda[-]{\aboveT{}}=0$, establishing the corollary in this case.

If $p$ splits in $F/\QQ$, write $p\rint{F}=\primeid\conjprimeid$ and let $\Primeid{}=\Primeid{}^{(1)},\ldots,\Primeid{\empty}^{(s)}$ (\rsp $\Conjprimeid{\empty}=\Conjprimeid{\empty}^{(1)},\ldots\Conjprimeid{\empty}^{(s)}$) be the prime ideals of $L_\infty$ that lie above $\primeid$ (\rsp above $\conjprimeid$). Let $\Lcyc$ be the cyclotomic $\ZZ_p$-extension of $F=L_0$ and set $\Ldue=L_{cyc}L_\infty$. Since $F/\QQ$ is abelian and totally imaginary, it satisfies the assumptions of Lemma~\ref{lemma:Gross} by~\cite[Corollary~2.14]{Gro81} and~\cite[Corollary~5.32]{Was97}. Now, observe that $\Primeid{}$ is unramified in $\Ldue/L_\infty$ because $F_{\primeid_i}\cong \QQ_p$ and $\QQ_p$ admits no totally ramified $\ZZ_p^2$-extension. Lemma~\ref{lemma:Gross} implies that the global extension $\Ldue/L_\infty$ is unramified everywhere, because it is clearly unramified outside of $p$. Denoting by $\Minfty$ the maximal abelian $p$-extension of $L_\infty$ which is unramified everywhere, we find the following diagram of fields and of Galois groups:
\[\xymatrix{
&\Minfty\ar@{-}[1,0]^H\ar@{-}@/_1pc/[2,-1]_{X}\\
&\Ldue\ar@{-}[1,1]\ar@{-}[1,-1]^{Y=X/H}\\
L_\infty\ar@{-}[1,1]\ar@{-}@/_1pc/[2,1]_{\DD}&&\Lcyc\ar@{-}[1,-1]\\
&F\ar@{-}[1,0]\\
&\QQ
}\]
For all $1\leq j\leq s$, consider the Frobenius elements $\Frob{\Primeid{\empty}^{(j)}}{\Minfty/L_\infty}$ and $\Frob{\Conjprimeid{\empty}^{(j)}}{\Minfty/L_\infty}$ in $X=\Gal{\Minfty}{L_\infty}$. The group $\DD$ acts on the set $\bigl\{\Frob{\Primeid{\empty}^{(j)}}{\Minfty/L_\infty},\Frob{\Conjprimeid{\empty}^{(j)}}{\Minfty/L_\infty}\bigr\}_{\;j=1,\ldots,s}$ by conjugation: given $1\leq j\leq s$, for each $\tau\in \DD$ and a lift $\tilde{\tau}\in \Gal{\Minfty}{\QQ}$ we have
\[
\tilde{\tau}\Frob{\Primeid{\empty}^{(j)}}{\Minfty/L_\infty}\tilde{\tau}^{-1}= \Frob{\tau\Primeid{\empty}^{(j)}}{\Minfty/L_\infty}
\]
and similarly when replacing $\Primeid{\empty}^{(j)}$ by $\Conjprimeid{\empty}^{(j)}$. On the one hand, this action becomes trivial on the classes
\[
\Frob{\Primeid{\empty}^{(j)}}{\Minfty/L_\infty}\pmod H\qquad\text{ and }\qquad\Frob{\Conjprimeid{\empty}^{(j)}}{\Minfty/L_\infty}\pmod H
\]
in $Y$, because $\Gal{\Ldue}{\QQ}\cong\DD\times \Gal{\Lcyc}{F}$ with $\DD$ acting trivially on $Y$. On the other, the action factors through $\splitG[{\indextr[\empty]}]=\Gal{L_{\indextr[\empty]}}{\QQ}$, where $L_{\indextr[\empty]}$ is such that $L_\infty/L_{\indextr[\empty]}$ is totally ramified at one (and hence all) primes above $p$, because $\Gamma_{\indextr[\empty]}$ acts trivially on all $\Primeid{}^{(j)},\Conjprimeid{}^{(j)}$. Since $\Frob{\Primeid{}^{(j)}}{\Minfty/L_\infty}\pmod{H}=\Frob{\Primeid{}^{(j)}}{\Ldue/L_\infty}$ and similarly when replacing $\Primeid{\empty}^{(j)}$ by $\Conjprimeid{\empty}^{(j)}$, it follows that, in $Y$, the equality
\begin{equation}\label{eq:frob_modH}
\prod_j\Frob{\Primeid{}^{(j)}}{\Ldue/L_\infty}\Frob{\Conjprimeid{}^{(j)}}{\Ldue/L_\infty}\pmod{H}=\prod_{\tau\in\splitG[{\indextr[\empty]}]}\Frob{\tau\Primeid{}}{\Ldue/L_\infty}=\Frob{\Primeid{}}{\Ldue/L_\infty}^{2s}
\end{equation}
holds. Having set this up, consider the element
\[
x=\prod_j\Frob{\Primeid{}^{(j)}}{\Minfty/L_\infty}\Frob{\Conjprimeid{}^{(j)}}{\Minfty/L_\infty}\in X.
\]
Lemma~\ref{lemma:Gross} and \eqref{eq:frob_modH} imply that $x\pmod{H}$ is the $(2s)$th power of an element of infinite order, so it again has infinite order: in particular, $x$ has infinite order as element of $X$. Via the Artin map of global class field theory, $x$ corresponds to the inverse limit of the generators of $\aboveT{n}^+$: in particular, the orders $\gorder{\aboveT{n}^+}$ cannot be bounded, showing that $\iwlambda[+]{\aboveT{}}\geq 1$. Since $\oram{L_\infty}{F}=2$, \eqref{eq:oram_quad} yields $\iwlambda{\aboveT{}}=\iwlambda[+]{\aboveT{}}+\iwlambda[-]{\aboveT{}}=1$, and therefore we must have $\iwlambda[+]{\aboveT{}}=1$, finishing the proof of the corollary.
\end{proof}
\begin{remark}\label{rmk:Finis_Hida_Tilouine}
It has been shown in~\cite[Corollary~26]{HubWas18} that, in some special cases, $\iwmu{\iw}=0$ for the anticyclotomic extension of an imaginary quadratic field: in particular, in these cases, $\iwmu{\fk}=0$. Finis and Hida (see~\cite{Fin06, Hid10} and the references therein) prove that, under mild assumptions, the \emph{analytic} anticyclotomic invariant $\iwmu[\mathit{an}]{\iw}$ vanishes. As Hida and Tilouine observe in \cite[\S 1, Comment 4]{HidTil94} the fact that their proof of the anticyclotomic main conjecture requires tensoring with $\QQ$ prevents a direct comparison between $\iwmu[\mathit{an}]{\iw}$ and $\iwmu{\iw}$. 
\end{remark}

As it is clear from the proof of Corollary~\ref{cor:formula_quad}, the situation where $p$ does not split in $F/\QQ$ is significantly easier to consider. In particular, the full assumptions that $F/\QQ$ is imaginary quadratic is not needed, and we single out this case in the next corollary:

\begin{corollary}\label{cor:special_case_inert}
Suppose that $F$ is a CM field and that $k$ is its totally real subfield. Assume that there is a unique prime above $p$ in $F/\QQ$. Then
\[
\iwlambda{\fk}=\frac{\iwlambda{\iw}}{2}
\]
so that $\iwlambda{\iw}$ is even.
\end{corollary}
\begin{proof} By Proposition~\ref{prop:formula_con_hp}-\ref{list:thm_formula:CM}, we have $\lambdacohom{0}^+=\oram{L_\infty}{F}-\iwlambda{\aboveT{}}-1$. On the other hand, since there is a unique prime above $p$ in $F$, the constant $\oram{L_\infty}{F}$ equals $1$ and from the trivial bound $\lambdacohom{0}^+\geq 0$ we obtain
\[
\lambdacohom{0}^+=\iwlambda{\aboveT{}}=0.
\]
In particular, $\iwlambda[-]{\aboveT{}}=0$. Moreover, since the unique prime in $k$ above $p$ does not split in $F/k$, we also have $\oram{L_\infty}{k}_s=0$, so that Proposition~\ref{prop:arithmetic_interpretation_lambda1} yields $\lambdacohom{-1}^+=0$. The result follows from Theorem~\ref{thm:iwasawa_formula}.
\end{proof}
\begin{remark}\label{rmk:Gillard_Carrol_Kisilevsky} The parity result concerning the Iwasawa invariant $\iwlambda{\iw}$ of the anticyclotomic $\ZZ_p$-extension of an imaginary quadratic field $F$ obtained in Corollary~\ref{cor:formula_quad} is not new: it has already been found by Gillard in \cite[Théorème~2 and Théorème~1--Corollaire]{Gil76}, and by Carroll and Kisilevsky in \cite[Theorem~5]{CarKis81}. Moreover, Carroll and Kisilevsky obtain \ibid a result similar to Corollary~\ref{cor:special_case_inert}, which holds for more general semi-direct products beyond the pro-dihedral case, although under the additional assumption that $F/\QQ$ is abelian with Galois group of exponent divisible by $(p-1)$.
\end{remark}

The next corollary improves Corollary~\ref{cor:bounds} under the assumption that $F$ is a CM field, generalizing Corollary~\ref{cor:formula_quad} to more general number fields beyond the imaginary quadratic case. As observed in Remark~\ref{rmk:hps_cor_CM}, its assumptions are satisfied when $F/\QQ$ is abelian.

\begin{corollary}\label{cor:Ass+CM} 
Suppose that $F/\QQ$ is a CM field of degree $2d$, and that $k$ is its totally real subfield. Assume that $p$ splits completely in $F/\QQ$ and that $L_\infty/\QQ$ is Galois. Then, with notation as in Proposition~\ref{prop:arithmetic_interpretation_lambda1},
\begin{equation}\label{ineq:cor_Ass+CM}
\iwlambda{\fk}\geq \frac{\iwlambda{\iw}+1-d}{2}.
\end{equation}
If, moreover, $F/\QQ$ is Galois and satisfies both the Gross and the Leopoldt conjectures, then $\iwlambda{\aboveT{}}=d$,  $\iwlambda[-]{\aboveT{}}\leq d-1$, $\iwlambda[+]{\aboveT{}}\geq 1$ and
\begin{equation}\label{ineq:cor_Ass+CM+abelian}
\iwlambda{\fk}\geq \frac{\iwlambda{\iw}+2-d}{2}.
\end{equation}
\end{corollary}
\begin{proof}
The hypotheses imply that $\oram{L_\infty}{F}=2d$ and $\oram{L_\infty}{k}_s=[k:\QQ]=d$. Moreover, condition~\ref{list:thm_formula:CM} of Proposition~\ref{prop:formula_con_hp} is satisfied: we find
\begin{equation}\label{eq:lambda0D_Ass+CM}
\lambdacohom{0}^+=\oram{L_\infty}{F}-\iwlambda{\aboveT{}}-1=2d-\iwlambda{\aboveT{}}-1.
\end{equation}
The formula of Theorem~\ref{thm:iwasawa_formula} becomes, through~\eqref{eq:lambda0D_Ass+CM} and Proposition~\ref{prop:arithmetic_interpretation_lambda1},
\[
2\iwlambda{\fk}=\iwlambda{\iw}+\bigl(d-\iwlambda[-]{\aboveT{}}\bigr)-(2d-\iwlambda{\aboveT{}}-1)=\iwlambda{\iw}-d+\iwlambda[+]{\aboveT{}}+1.
\]
Since $\iwlambda[+]{\aboveT{}}\geq 0$, the estimate~\eqref{ineq:cor_Ass+CM} follows.

Assume the Gross and the Leopoldt conjectures for $F$ and that $F/\QQ$ is Galois: the estimate in case~\ref{list:thm_formula:Ass} of Proposition~\ref{prop:formula_con_hp} reads $\lambdacohom{0}^+=d-1$, which yields $\iwlambda{\aboveT{}}=d$ through \eqref{eq:lambda0D_Ass+CM}. We now proceed to prove that $\iwlambda[+]{\aboveT{}}\geq 1$:  this will imply both estimates $\iwlambda[-]{\aboveT{}}\leq d-1$ and~\eqref{ineq:cor_Ass+CM+abelian}. The proof goes along the same lines as those of Corollary~\ref{cor:formula_quad}. Denote by $\GG$ the Galois group $\Gal{L_\infty}{\QQ}$ and set $\Ldue=L_\infty\Lcyc$, where $\Lcyc$ is the cyclotomic $\ZZ_p$-extension of $F$. Since $F/\QQ$ is Galois, so is $\Lcyc/\QQ$ and the Galois group $\Gal{\Ldue}{\QQ}$ is the direct product $\GG\times\Gal{\Lcyc}{F}$. Write $p\rint{k}=\overline{\conjprimeid_1}\cdots\overline{\conjprimeid_d}$ and $\overline{\conjprimeid_i}\rint{F}=\conjprimeid_{i}\conjprimeid_{d+i}$: each prime $\conjprimeid_j$ is divided by the same number of primes in $L_\infty/F$, so let $\conjprimeid_j\rint{L_\infty}=\Conjprimeid{j}^{(1)}\cdots\Conjprimeid{j}^{(s)}$ be its factorisation in $L_\infty$, for $1\leq j\leq 2d$. The localisations~$F_{\conjprimeid_j}$ are all isomorphic to $\QQ_p$, and hence admit no totally ramified $\ZZ_p^2$-extension, showing that $\Conjprimeid{\empty}=\Conjprimeid{1}^{(1)}$ is unramified in $\Ldue/L_\infty$ and so Lemma~\ref{lemma:Gross} guarantees that $\Ldue\subseteq \Minfty$, where $\Minfty$ is the maximal $p$-extension of~$L_\infty$ which is unramified everywhere. The diagram of fields is as follows:
\[\xymatrix{
&\Minfty\ar@{-}[1,0]^H\ar@{-}@/_1pc/[2,-1]_{X}\\
&\Ldue\ar@{-}[1,1]\ar@{-}[1,-1]^{Y=X/H}\\
L_\infty\ar@{-}[1,1]\ar@{-}@/_1pc/[2,1]_{\GG}&&\Lcyc\ar@{-}[1,-1]\\
&F\ar@{-}[1,0]\\
&\QQ
}\]
The group $\GG$ permutes the ideals $\Conjprimeid{j}^{(i)}$ and acts trivially on $\Gal{\Ldue}{L_\infty}$, so the equivalent of~\eqref{eq:frob_modH} reads
\begin{equation}\label{eq:frob_modH_Ass+CM+abelian}
\prod_{i,j}\Frob{\Conjprimeid{j}^{(i)}}{M_\infty/L_\infty}\pmod{H}=\prod_{\tau\in\Gal{L_{\indextr[\empty]}}{\QQ}}
\Frob{\tau\Conjprimeid{}}{\Ldue/L_\infty}=\Frob{\Conjprimeid{}}{\Ldue/L_\infty}^{2ds}\in Y
\end{equation}
where, as \ibid, $L_{\indextr[\empty]}$ is such that $L_\infty/L_{\indextr[\empty]}$ is totally ramified at all primes above $p$.
The element
\[
x=\prod_{i,j}\Frob{\Conjprimeid{j}^{(i)}}{\Minfty/L_\infty}\in X
\]
has infinite order by Lemma~\ref{lemma:Gross} and corresponds, via the Artin map of global class field theory, to an element in $\varprojlim \aboveT{n}^+$. This shows $\iwlambda{\aboveT{}^+}\geq 1$, finishing the proof.
\end{proof}

\begin{remark}\label{rmk:hps_cor_CM}
In the proof of~\eqref{ineq:cor_Ass+CM}, the  hypothesis that $L_\infty/\QQ$ is Galois is only needed to ensure that all primes in $F$ above $p$ ramify completely in $L_\infty/\QQ$, in order to apply Proposition~\ref{prop:formula_con_hp}-\ref{list:thm_formula:Ass}. In particular, the inequality holds also under this weaker assumption.

One condition that guarantees that all assumptions of Corollary~\ref{cor:Ass+CM} are satisfied is that the group $\Gal{F}{\QQ}$ is abelian with exponent divisible by $(p-1)$: this is shown in \cite[Theorem~1]{CarKis81}, together with the result that both the Gross and the Leopoldt conjectures hold for abelian number fields.
\end{remark}

\subsection{The structure of the limits of class groups}\label{subsec:structure_limit}
As in the previous section, fix a fake $\ZZ_p$-extension of dihedral type $K_\infty/k$ and adopt the conventions introduced in~\S\ref{subsec:fake_def}. Let $\calB=(B_{L_n})\in\{\unit{\pallino},\idunit{\pallino},\q{\pallino},\idclg{\pallino},\clsyl{\pallino}\}$ be any of the five normic systems occurring either in the left column or in the bottom row of~\eqref{diag:magic}. As in~\S\ref{subsec:cohom_units}, set $B_{L_\infty}=\varinjlim B_{L_n}$ (\rsp $B_{K_\infty}=\varinjlim B_{K_n}$, letting $B_{K_n}=B_{L_n}^+$) where the limit is taken with respect to $\arext{L_m}{L_n}$ (\rsp with respect to $\arext{K_m}{K_n}$).

In this section we make no use of double systems: only direct and inverse systems will appear. Indeed, the focus will be on the structure of direct and inverse \emph{limits}, which are coarser invariants. Moreover, one easily realises that the direct system of class groups with ideals extension as morphisms satisfies neither~\ref{cond:inj} nor~\ref{cond:surj}, and the assumptions of Proposition~\ref{prop:long_cohom_down} are not fulfilled, so its $G_n$-cohomology does not give rise to a double system fitting in a cohomological long exact sequence. In fact, the analysis of the direct and inverse limits of class groups relies on our previous work~\cite{CapNuc20}, and we need to consider morphisms between $0$th Tate cohomology groups which are different from the map $\HHm{0}{\arext{L_m}{L_n}}$ defined in \S\ref{subsubsub:HH0}, because of the factor $p^{m-n}$ appearing there. To define these morphisms, observe that, for $m\geq n \gg 0$, $\arnm{L_m}{L_n}\colon\clsyl{L_m}\to\clsyl{L_n}$ is surjective and so, by the arguments of Remark~\ref{rmk:h0surjnorm}, we have maps
\[
\hatnsu[\arextsymb]{L_m/L_n}\colon\HH{0}{G_n}{\clsyl{L_n}}\longrightarrow \HH{0}{G_m}{\clsyl{L_m}}\qquad\text{ and }\qquad \hatnsu[\arextsymb]{L_m/L_n}^+\colon\HH{0}{\splitG[n]}{\clsyl{L_n}}\longrightarrow \HH{0}{\splitG[m]}{\clsyl{L_m}}.
\]
We can therefore take direct limits with respect to these maps, and define the groups\footnote{In \cite[Definition~I.1.9.3]{NeuSchWin08}, the notation $\HH{0}{\mathcal{G}}{M}$ denotes a different construction, which we do not need, hence we allow ourselves this inconsistency.}
\[
\HH{0}{\Gamma}{\clsyl{L_\infty}} := \varinjlim_{\hatnsu[\arextsymb]{L_m/L_n}}\HH{0}{G_n}{\clsyl{L_n}} \qquad \text{ and }\qquad\HH{0}{\DD}{\clsyl{L_\infty}} := \varinjlim_{\hatnsu[\arextsymb]{L_m/L_n}^+} \HH{0}{\splitG[n]}{\clsyl{L_n}}.
\]
Being direct limits of finite $\splitG$-modules, these are $\ZZ_p$-torsion $\ZZ_p[\splitG]$-modules; moreover, via Proposition~\ref{prop:fixcohom}, they satisfy $\HH{0}{\Gamma}{\clsyl{L_\infty}}^+\cong \HH{0}{\DD}{\clsyl{L_\infty}}$, extending the isomorphisms \ibid for $i\geq 1$. More generally, in this section, we tacitly identify plus parts of all $G_n$-cohomology groups (\rsp of all $\Gamma$-cohomology groups) with $\splitG[n]$-cohomology (\rsp $\DD$-cohomology).

\subsubsection{An Euler--Poincaré formula for fake $\ZZ_p$-extensions of dihedral type}
The main result of this section is Theorem~\ref{thm:herbrand_coranks}, which can be seen as a relationship between dihedral versions of Herbrand quotients of class groups and of units, at least in the limit. This theorem will also be pivotal in the next section to interpret the invariant $\iwlambda{\fk}$ as the dimension of a $\QQ_p$-vector space. 

We begin by relating $\HHN{1}{\DD}{\clsyl{L_\infty}}$ to $\HHN{2}{\DD}{\q{L_\infty}}$ on the one hand, and $\HH{0}{\DD}{\clsyl{L_\infty}}$ to $\HHN{1}{\DD}{\q{L_\infty}}$ on the other. 

\begin{proposition}\label{prop:qacinfty}
There is an isomorphism of $\ZZ_p$-modules
\[
\HH{0}{\DD}{\clsyl{L_\infty}}\cong\HHN{1}{\DD}{\q{L_\infty}}
\] 
and an exact sequence
\[
0\longdashrightarrow \HHN{1}{\DD}{\clsyl{L_\infty}} \longdashrightarrow \HHN{2}{\DD}{\q{L_\infty}}\longdashrightarrow \HHN{2}{\DD}{\idclg{L_\infty}}\longdashrightarrow 0
\]
in the category $\catMod[]{\cofg}/\catMod[]{fin}$. Moreover, $\HHN{2}{\DD}{\idclg{L_\infty}}\cong \QQ_p/\ZZ_p$.
\end{proposition}
\begin{proof}
By Proposition \ref{prop:long_cohom_up}, we have a commutative diagram
\begin{equation}\begin{split}\label{diag:2.10+}
\xymatrix@C=4em{
	\HHN{0}{\splitG[m]}{\clsyl{L_m}}\ar@{->}[0,1]^{\cch{0}}&\HH{1}{\splitG[m]}{\q{L_m}}\\
	\HHN{0}{\splitG[n]}{\clsyl{L_n}}\ar@{->}[0,1]^{\cch{0}}\ar@{->}[-1,0]_{\arext{L_m}{L_n}}&\HH{1}{\splitG[n]}{\q{L_n}}\ar@{->}[-1,0]_{\HHm{1}{\arext{L_m}{L_n}}}
}
\end{split}\end{equation}
The definition of $\HHm{1}{\arext{L_m}{L_n}}$ given in Remark~\ref{rmk:h0surjnorm} makes it clear that~\eqref{diag:2.10+} induces the commutative diagram
\[
\xymatrix@C=4em{
	\HH{0}{\splitG[m]}{\clsyl{L_m}}\ar@{->}[0,1]^{\cchT{0}}&\HH{1}{\splitG[m]}{\q{L_m}}\\
	\HH{0}{\splitG[n]}{\clsyl{L_n}}\ar@{->}[0,1]^{\cchT{0}}\ar@{->}[-1,0]_{\hatnsu[\arextsymb]{L_m/L_n}}&\HH{1}{\splitG[n]}{\q{L_n}}\ar@{->}[-1,0]_{\HHm{1}{\arext{L_m}{L_n}}}
}
\]
From \cite[Proposition 3.4]{CapNuc20}, the rows of the above diagram are isomorphisms, and so by taking limits we obtain the first statement of the proposition.

Arguing as in the proof of \cite[Proposition 3.4]{CapNuc20}, we get an exact sequence
\[
0\longrightarrow\HH{1}{G_n}{\clsyl{L_n}} \longrightarrow \HH{2}{G_n}{\q{L_n}}\longrightarrow \HH{2}{G_n}{\idclg{L_n}} \longrightarrow \HH{2}{G_n}{\clsyl{L_n}} \longrightarrow \HH{3}{G_n}{\q{L_n}} \longrightarrow 0
\] 
with its plus counterpart
\[
0\longrightarrow \HH{1}{\splitG[n]}{\clsyl{L_n}} \longrightarrow \HH{2}{\splitG[n]}{\q{L_n}}\longrightarrow \HH{2}{\splitG[n]}{\idclg{L_n}} \longrightarrow \HH{2}{\splitG[n]}{\clsyl{L_n}} \longrightarrow \HH{3}{\splitG[n]}{\q{L_n}} \longrightarrow 0.
\]
By Proposition~\ref{prop:fixcohom}, this exact sequence can be shortened as
\begin{equation}\label{eq:qcafinite}
0\longrightarrow \HH{1}{\splitG[n]}{\clsyl{L_n}} \longrightarrow \HH{2}{\splitG[n]}{\q{L_n}}\longrightarrow \HH{2}{\splitG[n]}{\idclg{L_n}} \longrightarrow \kernel{\bigl(\cchT{0}_n\bigr)^-}\longrightarrow 0
\end{equation}
where $\cchT{0}_n\colon\HH{0}{G_n}{\clsyl{L_n}} \to \HH{1}{G_n}{\q{L_n}}$ is the morphism analysed in Remark~\ref{rmk:ker_delta_bdd}. As shown \ibid, $\kernel(\cchT{0}_n)$ is bounded independently of $n$, and taking direct limits of \eqref{eq:qcafinite} we obtain the exact sequence of the statement.

Concerning the final isomorphism, an argument similar to the one used to prove the commutativity of the left-hand diagram of Lemma~\ref{lemma:fundclasscommute} shows that the direct system $(\HH{2}{\splitG[n]}{\idclg{L_n}},\HHm{2}{\arext{L_m}{L_n}})$ is isomorphic to $(\HH{0}{\splitG[n]}{\ZZ_p},\HHm{0}{\nsu_{n,m}})\cong (\frac{1}{p^n}\ZZ_p/\ZZ_p,\operatorname{incl})$ (where $\operatorname{incl}$ denotes the inclusion $\frac{1}{p^n}\ZZ\hookrightarrow\frac{1}{p^m}\ZZ$), finishing the proof.
\end{proof}

\begin{remark}
The proof of the above proposition actually shows that there is an exact sequence
\[
0\longrightarrow \HHN{1}{\DD}{\clsyl{L_\infty}} \longrightarrow \HHN{2}{\DD}{\q{L_\infty}}\longrightarrow \HHN{2}{\DD}{\idclg{L_\infty}}\longrightarrow V \longrightarrow0
\]
in the category $\catMod[]{\cofg}$ with $V$ in $\catMod[]{fin}$. Since we are only interested in comparing coranks of divisible $\ZZ_p$-modules, here and in what follows we content ourselves to work in $\catMod[]{\cofg}/\catMod[]{fin}$.
\end{remark}

We can now proceed to the proof of the main result of this section, that we regard as an Euler--Poincaré formula:

\begin{theorem}\label{thm:herbrand_coranks}
Given a fake $\ZZ_p$-extension of dihedral type $K_\infty/k$ with normal closure $L_\infty/k$, we have
\[
\corank \HH{0}{\DD}{\clsyl{L_\infty}} - \corank\HHN{1}{\DD}{\clsyl{L_\infty}} = \corank\HHN{2}{\DD}{\unit{L_\infty}} - \corank\HHN{1}{\DD}{\unit{L_\infty}} + 1.
\]
\end{theorem}
\begin{proof}
Taking direct limits in~\eqref{eq:long_cohom_quot} we obtain the exact sequence in $\catMod[]{\cofg}/\catMod[]{fin}$
\begin{equation*}
0\longdashrightarrow\HHN{1}{\Gamma}{\unit{L_\infty}} \longdashrightarrow \HHN{1}{\Gamma}{\idunit{L_\infty}}\longdashrightarrow \HHN{1}{\Gamma}{\q{L_\infty}} \longdashrightarrow \HHN{2}{\Gamma}{\unit{L_\infty}}\overset{\overline{v}_\infty}{\longdashrightarrow} \HHN{2}{\Gamma}{\idunit{L_\infty}}.
\end{equation*}
By construction, the rightmost morphism $\overline{v}_\infty\colon\HHN{2}{\Gamma}{\unit{L_\infty}}\dashrightarrow \HHN{2}{\Gamma}{\idunit{L_\infty}}$ is the image, in the quotient category, of the direct limit of the morphisms $\HHN{2}{G_n}{\unit{L_n}}\to \HHN{2}{G_n}{\idunit{L_n}}$. It follows that we can complete the above sequence as
\[
0\dashrightarrow\HHN{1}{\Gamma}{\unit{L_\infty}} \dashrightarrow \HHN{1}{\Gamma}{\idunit{L_\infty}}\dashrightarrow \HHN{1}{\Gamma}{\q{L_\infty}} \dashrightarrow \HHN{2}{\Gamma}{\unit{L_\infty}}\dashrightarrow \HHN{2}{\Gamma}{\idunit{L_\infty}}\dashrightarrow \HHN{2}{\Gamma}{\q{L_\infty}} \dashrightarrow 0
\]
where the final $0$ is a consequence of $\Gamma$ having strict cohomological dimension equal to $2$. By taking plus parts we get, through Corollary~\ref{cor:growth_cohom_units_dihedral}-\eqref{eq:iso_idunits_QZ_signed},
\begin{align*}
\corank \HHN{1}{\DD}{\q{L_\infty}} - \corank\HHN{2}{\DD}{\q{L_\infty}}&	=
	\begin{aligned}[t] \corank\HHN{1}{\DD}{\idunit{L_\infty}} - \corank \HHN{2}{\DD}{\idunit{L_\infty}} + \corank&\HHN{2}{\DD}{\unit{L_\infty}} \\
	&- \corank\HHN{1}{\DD}{\unit{L_\infty}}
\end{aligned}\\
&= \corank\HHN{2}{\DD}{\unit{L_\infty}} - \corank\HHN{1}{\DD}{\unit{L_\infty}}.
\end{align*}
Moreover, by Proposition~\ref{prop:qacinfty},
\begin{align*}
\corank\HHN{2}{\DD}{\q{L_\infty}} &= \corank\HHN{1}{\DD}{\clsyl{L_\infty}} + 1\\
\intertext{and}
\corank\HHN{1}{\DD}{\q{L_\infty}} &= \corank\HH{0}{\DD}{\clsyl{L_\infty}},
\end{align*}
finishing the proof.
\end{proof}

\subsubsection{Interpretation of $\iwlambda{\fk}$} 
Let $\iwLambda=\ZZ_p[\![\Gamma]\!]$ denote the Iwasawa algebra of the group $\Gamma=\Gal{L_\infty}{F}$. In this section we show that the invariant $\iwlambda{\fk}$ obtained in Theorem~\ref{thm:iwasawa_formula} is also the dimension of the $\QQ_p$-vector space $\iwXK\otimes_{\ZZ_p}\QQ_p$, where 
\[
\iwXK = \varprojlim A_{K_n},
\]
the inverse limit being taken with respect to norms. The analogous result for $\iwXL$ (with $\iwlambda{\iw}$ replacing $\iwlambda{\fk}$) is well-known and comes from the structure theorem for finitely generated torsion $\Lambda$-modules together with a classical descent argument (see~\cite[\S1.2 and~Theorem~7]{Iwa73a}). Since $\iwXK$ is not stabilized by $\Gamma$, and hence it is not a $\Lambda$-submodule of $\iwXL$, we argue in a different way. Our approach is to work with direct, rather than inverse, limits because they are easier to analyse in our setting. We then recover the result about inverse limits via a purely algebraic argument due to Yamashita.

Tensoring the exact sequence~\cite[(3.1)]{CapNuc20} with $\ZZ_p$, we obtain
\begin{equation}\label{eq:nuccioskernel_finite}
1\longrightarrow \kernel{\nuceta{n}}\longrightarrow \clsyl{K_n} \oplus \clsyl{K_n'} \overset{\nuceta{n}}{\longrightarrow} \clsyl{L_n} \longrightarrow \cokernel{\nuceta{n}}\longrightarrow 1
\end{equation}
where the map $\nuceta{n}$ is defined via the formula
\begin{equation*}
\nuceta{n}\left(c,c'\right)=\arext{L_n}{K_n}(c)\arext{L_n}{K_n'}(c') \qquad \text{ for }c\in \clsyl{K_n}, c'\in \clsyl{K_n'}.
\end{equation*}
The explicit definition of $\nuceta{n}$ makes it evident that, for $m\geq n$, the relations
\begin{equation}\label{eq:nuceta_ext}
\nuceta{m}\circ\bigl(\arext{K_m}{K_n}\oplus\arext{K_m'}{K_n'}\bigr)=\arext{L_m}{L_n}\circ\nuceta{n}
\end{equation}
hold. We can therefore take direct limits with respect to the extension maps to obtain an exact sequence
\begin{equation}\label{eq:nuccioskernel_infinite}
1\longrightarrow \kernel{\nuceta{\infty}}\longrightarrow \clsyl{K_\infty} \oplus \clsyl{K'_\infty} \overset{\nuceta{\infty}}{\longrightarrow} \clsyl{L_\infty} \longrightarrow \cokernel{\nuceta{\infty}}\longrightarrow 1
\end{equation}
where $\nuceta{\infty}=\varinjlim\nuceta{n}$.

The first step of our strategy is to interpret, in Lemmas~\ref{lemma:nkerinfty} and~\ref{lemma:ncokerinfty}, $\kernel{\nuceta{\infty}}$ and $\cokernel{\nuceta{\infty}}$  in terms of the cohomology groups of $\clsyl{L_\infty}$. We then invoke Theorem~\ref{thm:herbrand_coranks} to compare their coranks.

\begin{lemma}\label{lemma:nkerinfty}
The modules $\kernel{\nuceta{\infty}}$ and $\HH{0}{\DD}{\clsyl{L_\infty}}$ are $\catMod[]{fin}$-isomorphic. 
\end{lemma}
\begin{proof}
Recall from \cite[Proposition~3.2]{CapNuc20} that
\[
\kernel{\nuceta{n}} \cong \HHN{0}{\splitG[n]}{\clsyl{L_n}}
\]
and the isomorphism, called $\theta$ \ibid, is defined by sending $(x, x')\in\kernel{\nuceta{n}}\subseteq\clsyl{K_n}\oplus\clsyl{K'_n}$ to $\arext{L_n}{K_n}(x)$. Hence, the diagram for $m\geq n\geq 0$
\[
\xymatrix@C=5em{
	\kernel{\nuceta{m}}\ar@{->}[0,1]^(.40){\sim}& \HHN{0}{\splitG[m]}{\clsyl{L_m}}\\
	\kernel{\nuceta{n}}\ar@{->}[-1,0]^{\arext{K_m}{K_n}\oplus\arext{K_m'}{K_n'}}\ar@{->}[0,1]^(.40){\sim}& \HHN{0}{\splitG[n]}{\clsyl{L_n}}\ar@{->}[-1,0]_{\arext{L_m}{L_n}}
}
\]
commutes, showing that 
\begin{equation}\label{eq:nuckerinf}
\kernel{\nuceta{\infty}} \cong \HHN{0}{\DD}{\clsyl{L_\infty}}.
\end{equation}
The definition of Tate cohomology groups gives another commutative diagram
\begin{align*}
\xymatrix{
	0\ar@{->}[0,1]&\algnm{\splitG[m]}\clsyl{L_m} \ar@{->}[0,1]&\HHN{0}{\splitG[m]}{\clsyl{L_m}}\ar@{->}[0,1]& \HH{0}{\splitG[m]}{\clsyl{L_m}}\ar@{->}[0,1]&0\\
	0\ar@{->}[0,1]&\algnm{\splitG[n]}\clsyl{L_n}\ar@{->}[-1,0]^{\arext{L_m}{L_n}}\ar@{->}[0,1]&\HHN{0}{\splitG[n]}{\clsyl{L_n}}\ar@{->}[-1,0]^{\arext{L_m}{L_n}}\ar@{->}[0,1]& \HH{0}{\splitG[n]}{\clsyl{L_n}}\ar@{->}[-1,0]^{\hatnsu[\arextsymb]{L_m/L_n}}\ar@{->}[0,1]&0
}
\end{align*}
yielding the exact sequence
\[
0\longrightarrow \varinjlim\algnm{\splitG[n]}\clsyl{L_n}\longrightarrow \HHN{0}{\DD}{\clsyl{L_\infty}} \longrightarrow \HH{0}{\DD}{\clsyl{L_\infty}}\longrightarrow 0.
\] 
The groups $\algnm{\splitG[n]}\clsyl{L_n}$ have bounded order since 
\[
\gorder{\algnm{\splitG[n]}\clsyl{L_n}} = \gorder{\arext{L_n}{k}\arnm{L_n}{k}\clsyl{L_n}}\leq \gorder{\arext{L_n}{k}\clsyl{k}}\leq \gorder{\clsyl{k}},
\] 
so their direct limit is certainly a finite group. We conclude using \eqref{eq:nuckerinf}.
\end{proof}

\begin{lemma}\label{lemma:ncokerinfty}
The modules $\cokernel{\nuceta{\infty}}$ and $\HHN{1}{\DD}{\clsyl{L_\infty}}$ are $\catMod[]{fin}$-isomorphic.
\end{lemma}
\begin{proof}
As in the proof of~\cite[Proposition~3.2]{CapNuc20}, Remark~3.1 \ibid yields the exact sequence
\[
0 \longrightarrow \clsyl{L_m}[\algnm{G_m}]/\bigl(\clsyl{L_m}[\algnm{G_m}]\cap\im{\nuceta{m}}\bigr) \longrightarrow \cokernel{\nuceta{n}}\overset{\algnm{G_n}}{\longrightarrow} \algnm{G_n}\clsyl{L_n}/\algnm{\splitG[n]}\clsyl{L_n}\longrightarrow 0.
\]
We actually have
commutative diagrams, for $m\geq n\geq 0$,
\[
\xymatrix{
	0\ar@{->}[0,1]&\clsyl{L_m}[\algnm{G_m}]/\bigl(\clsyl{L_m}[\algnm{G_m}]\cap\im{\nuceta{m}}\bigr)\ar@{->}[0,1]&\cokernel{\nuceta{m}}\ar@{->}[0,1]^(.29){\algnm{G_m}}&\algnm{G_m}(\clsyl{L_m})/\algnm{G_m}(\im{\nuceta{L_m}})\ar@{->}[0,1]&0\\
	0\ar@{->}[0,1]&\clsyl{L_n}[\algnm{G_n}]/\bigl(\clsyl{L_n}[\algnm{G_n}]\cap\im{\nuceta{n}}\bigr)\ar@{->}[-1,0]^{\arext{L_m}{L_n}}\ar@{->}[0,1]&\cokernel{\nuceta{n}}\ar@{->}[-1,0]^{\arext{L_m}{L_n}}\ar@{->}[0,1]^(.29){\algnm{G_n}}&\algnm{G_n}(\clsyl{L_n})/\algnm{G_n}(\im{\nuceta{n}})\ar@{->}[-1,0]_{p^{m-n}\arext{L_m}{L_n}}\ar@{->}[0,1]&0	
}
\]
By~\cite[(2.7)]{CapNuc20} and by definition of $\HHm{-1}{\arext{L_m}{L_n}}$ (see~\S\ref{subsubsec:minus1}), this can be rewritten as
\[
\xymatrix@C=4em{
	0\ar@{->}[0,1]&\HH{-1}{G_m}{\clsyl{L_m}}^-\ar@{->}[0,1]&\cokernel{\nuceta{m}}\ar@{->}[0,1]^(.29){\algnm{G_m}}&\algnm{G_m}(\clsyl{L_m})/\algnm{G_m}(\im{\nuceta{L_m}})\ar@{->}[0,1]&0\\
	0\ar@{->}[0,1]&\HH{-1}{G_n}{\clsyl{L_n}}^-\ar@{->}[-1,0]^{\HHm{-1}{\arext{L_m}{L_n}}}\ar@{->}[0,1]&\cokernel{\nuceta{n}}\ar@{->}[-1,0]^{\arext{L_m}{L_n}}\ar@{->}[0,1]^(.29){\algnm{G_n}}&\algnm{G_n}(\clsyl{L_n})/\algnm{G_n}(\im{\nuceta{n}})\ar@{->}[-1,0]_{p^{m-n}\arext{L_m}{L_n}}\ar@{->}[0,1]&0	
}
\]
Moreover, the cup product with $\varkappa_n$ induces an isomorphism $\HH{-1}{G_n}{\clsyl{L_n}}^-\cong\HH{1}{\splitG[n]}{\clsyl{L_n}}$, and by definition of $\HHm{1}{\arext{L_m}{L_n}}$, together with Lemma~\ref{lemma:HHm(j)_and_inf}, the diagram becomes
\[
\xymatrix@C=4.35em{
	0\ar@{->}[0,1]&\HH{1}{\splitG[m]}{\clsyl{L_m}}\ar@{->}[0,1]&\cokernel{\nuceta{m}}\ar@{->}[0,1]^(.29){\algnm{G_m}}&\algnm{G_m}(\clsyl{L_m})/\algnm{G_m}(\im{\nuceta{L_m}})\ar@{->}[0,1]&0\\
	0\ar@{->}[0,1]&\HH{1}{\splitG[n]}{\clsyl{L_n}}\ar@{->}[-1,0]^{\infl\circ \arext{L_m}{L_n}^*}\ar@{->}[0,1]&\cokernel{\nuceta{n}}\ar@{->}[-1,0]^{\arext{L_m}{L_n}}\ar@{->}[0,1]^(.29){\algnm{G_n}}&\algnm{G_n}(\clsyl{L_n})/\algnm{G_n}(\im{\nuceta{n}})\ar@{->}[-1,0]_{p^{m-n}\arext{L_m}{L_n}}\ar@{->}[0,1]&0	
}
\]
The groups on the right have bounded order since $\gorder{\algnm{G_n}\clsyl{L_n}}\leq \gorder{\clsyl{F}}$, by the same argument as at the end of the proof of Lemma~\ref{lemma:nkerinfty}. Taking direct limits we obtain the statement.
\end{proof}

Combined with Theorem~\ref{thm:herbrand_coranks}, the above lemmas yield the following

\begin{corollary}\label{cor:corankslambda}
We have
\[
2 \corank \clsyl{K_\infty} - \corank \clsyl{L_\infty}= \corank \kernel{\nuceta{\infty}} -\corank \cokernel{\nuceta{\infty}} =  \lambdacohom{2}^+ - \lambdacohom{1}^+ + 1
\]
where $\lambdacohom{1}^+$ and $\lambdacohom{2}^+$ are the invariants defined in Corollary~\ref{cor:growth_cohom_units_dihedral}-\ref{cor:growth_cohom_units_dihedral_existence_global}.
\end{corollary}
\begin{proof}
From the exact sequence \eqref{eq:nuccioskernel_infinite}, we obtain
\begin{align*}
2 \corank \clsyl{K_\infty} - \corank \clsyl{L_\infty}&= \corank \kernel{\nuceta{\infty}} -\corank \cokernel{\nuceta{\infty}}\\
\text{(Lemmas \ref{lemma:nkerinfty} and \ref{lemma:ncokerinfty})} &=  \corank \HH{0}{\DD}{\clsyl{L_\infty}} -\corank\HHN{1}{\DD}{\clsyl{L_\infty}}\\
\text{(Theorem \ref{thm:herbrand_coranks})} &=  \corank\HHN{2}{\DD}{\unit{L_\infty}} - \corank\HHN{1}{\DD}{\unit{L_\infty}} + 1 \\
\text{(Corollary~\ref{cor:ha_vinto_caputo} and Corollary~\ref{cor:growth_cohom_units_dihedral}-\ref{cor:growth_cohom_units_dihedral_existence_global}		)}&=  \lambdacohom{2}^+ - \lambdacohom{1}^+ + 1.\qedhere
\end{align*}
\end{proof}

We are now in shape to relate the invariant $\iwlambda{\fk}$ to the $\ZZ_p$-module structure of $\iwXK$, but we first need a final lemma:

\begin{lemma}\label{lemma:rkcorkclgp}
The equalities
\[
\corank \clsyl{L_\infty} = \dim_{\QQ_p} \iwXL\otimes_{\ZZ_p} \QQ_p \qquad \text{ and } \qquad \corank \clsyl{K_\infty} = \dim_{\QQ_p} \iwXK\otimes_{\ZZ_p} \QQ_p
\]	
hold.
\end{lemma}
\begin{proof}
The first equality can be deduced from~\cite[Theorem 11]{Iwa73a}. Here we give an alternative proof, following that of \cite[Lemma~2]{Yam84}, which also yields the second equality.

Since $\iwXL$ is a finitely generated torsion $\iwLambda$-module, we know, by the structure theorem of such modules, that there exists $\beta\in\mathbb{N}$ such that $p^\beta\iwXL$ is a finitely generated torsion $\iwLambda$-module with trivial $\iwmu{}$-invariant. Since $p^\beta\iwXL=\varprojlim p^\beta\clsyl{L_n}$, we deduce that the $\ZZ_p$-ranks of $p^\beta\clsyl{L_n}$ are bounded and thus so are those of $p^\beta\clsyl{K_n}$. As arithmetic norms are eventually surjective on $\clsyl{L_n}$, their restriction to $p^\beta\clsyl{L_n}$ is also eventually surjective, and likewise when restricting to $p^\beta\clsyl{K_n}$. Similarly, the kernel of the extension map restricted to $p^\beta\clsyl{L_n}$ has bounded order and, again, the same holds when restricting to $p^\beta\clsyl{K_n}$. We can therefore apply \cite[Lemma 1]{Yam84} to obtain that 
\[
\corank (p^\beta\clsyl{L_\infty}) = \dim_{\QQ_p} (p^\beta\iwXL\otimes_{\ZZ_p} \QQ_p) \qquad\text{ and }\qquad \corank (p^\beta\clsyl{K_\infty}) = \dim_{\QQ_p} (p^\beta\iwXK\otimes_{\ZZ_p} \QQ_p).
\]
This concludes the proof because, for every $\ZZ_p$-module $M$, $\dim_{\QQ_p} (p^\beta M\otimes_{\ZZ_p} \QQ_p) = \dim_{\QQ_p} (M\otimes_{\ZZ_p} \QQ_p)$ and $\corank (p^\beta M) = \corank M$.
\end{proof}

Finally, we obtain the main result of this section.

\begin{theorem}\label{thm:lambda_equal_dim}
Given a fake $\ZZ_p$-extension of dihedral type $K_\infty/k$, we have $\iwlambda{\fk}=\dim_{\QQ_p} \iwXK\otimes_{\ZZ_p} \QQ_p$.
\end{theorem}
\begin{proof}
By the structure theorem of finitely generated $\iwLambda$-modules, we know that
\begin{equation}\label{eq:lambdaforLambdamodules}
\dim_{\QQ_p} \iwXL\otimes_{\ZZ_p} \QQ_p = \iwlambda{\iw}.
\end{equation}
Now, using Theorem~\ref{thm:iwasawa_formula}, we deduce 	
\begin{align*}
2\iwlambda{\fk} &= \iwlambda{\iw} + \lambdacohom{-1}^+-\lambdacohom{0}^+\\
\text{(Corollary~\ref{cor:growth_cohom_units_dihedral}-\eqref{eq:cohom_unit_many_calcs} and~-\eqref{eq:lambdaDperiodic})}& = \iwlambda{\iw} + \lambdacohom{2}^+-\lambdacohom{1}^+ + 1 \\
\text{(by~\eqref{eq:lambdaforLambdamodules})}& = \dim_{\QQ_p} \iwXL\otimes_{\ZZ_p} \QQ_p + \lambdacohom{2}^+-\lambdacohom{1}^+ + 1 \\
\text{(Lemma \ref{lemma:rkcorkclgp})}& = \corank \clsyl{L_\infty} + \lambdacohom{2}^+-\lambdacohom{1}^+ + 1\\
\text{(Corollary \ref{cor:corankslambda})}& = 2\corank \clsyl{K_\infty}\\
\text{(Lemma \ref{lemma:rkcorkclgp})}& = 2\dim_{\QQ_p}\left(\iwXK\otimes_{\ZZ_p}\QQ_p\right).\qedhere
\end{align*}
\end{proof}


\begin{thebibliography}{{Nuc}10}
\providecommand{\url}[1]{\texttt{#1}}
\providecommand{\urlprefix}{URL }
\expandafter\ifx\csname urlstyle\endcsname\relax
  \providecommand{\doi}[1]{doi:\discretionary{}{}{}#1}\else
  \providecommand{\doi}{doi:\discretionary{}{}{}\begingroup
  \urlstyle{rm}\Url}\fi
\providecommand{\eprint}[2][]{\url{#2}}

\bibitem[AW67]{AtiWal67}
\textsc{Atiyah, M.F.}; \textsc{Wall, C.T.C.}
\newblock Cohomology of groups.
\newblock Algebraic {N}umber {T}heory ({P}roc. {I}nstructional {C}onf.,
  {B}righton, 1965),  (Thompson, Washington, D.C.1967) 94--115. \mrev{0219512} \zbl{1494.20079}

\bibitem[Bou61]{Bou61}
\textsc{Bourbaki, N.}
\newblock \'{E}l\'{e}ments de math\'{e}matique. {F}ascicule {XXVII}.
  {A}lg\`ebre commutative. {C}hapitre 1: {M}odules plats. {C}hapitre 2:
  {L}ocalisation.
\newblock Actualit\'{e}s Scientifiques et Industrielles, No. 1290,  (Herman,
  Paris 1961). \mrev{0217051} \zbl{0108.04002}

\bibitem[CK82]{CarKis81}
\textsc{Carroll, J.}; \textsc{Kisilevsky, H.}
\newblock On {I}wasawa's {$\lambda $}-invariant for certain {${\bf
  Z}_{l}$}\thinspace -extensions.
\newblock \textit{Acta Arith.} \textbf{40}(1) (1981/82) 1--8.
\mrev{0654019} \zbl{0517.12003}

\bibitem[CN20]{CapNuc20}
\textsc{Caputo, L.}; \textsc{N{uccio~Mortarino~Majno~di~Capriglio}, F.A.E.}
\newblock Class number formula for dihedral extensions.
\newblock \textit{Glasgow Math. J.} \textbf{62}(2) (2020) 323–353.
\mrev{4085046} \zbl{1450.11117}

\bibitem[FG81]{FedGro81}
\textsc{Federer, L.J.}; \textsc{Gross, B.H.}
\newblock Regulators and {I}wasawa modules.
\newblock \textit{Invent. Math.} \textbf{62}(3) (1981) 443--457.
\mrev{0604838} \zbl{0468.12005}

\bibitem[Fin06]{Fin06}
\textsc{Finis, T.}
\newblock The {$\mu$}-invariant of anticyclotomic {$L$}-functions of imaginary
  quadratic fields.
\newblock \textit{J. Reine Angew. Math.} \textbf{596} (2006) 131--152.
\mrev{2254809} \zbl{1111.11048}

\bibitem[Gab62]{Gab62}
\textsc{Gabriel, P.}
\newblock Des cat{\'e}gories ab{\'e}liennes.
\newblock \textit{Bull. Soc. Math. France} \textbf{90} (1962) 323--448.
\mrev{0232821} \zbl{0201.35602}

\bibitem[Gil76]{Gil76}
\textsc{Gillard, R.}
\newblock Remarques sur certaines extensions prodi\'{e}drales de corps de
  nombres.
\newblock \textit{C. R. Acad. Sci. Paris S\'{e}r. A-B} \textbf{282}(1) (1976)
  A13--A15.
\mrev{0392910} \zbl{0335.12017}

\bibitem[Gre73]{Gre73}
\textsc{Greenberg, R.}
\newblock On a certain {$\ell$}-adic representation.
\newblock \textit{Invent. Math.} \textbf{21} (1973) 117--124.
\mrev{0335468} \zbl{0268.12004}

\bibitem[Gro57]{Gro57}
\textsc{Grothendieck, A.}
\newblock Sur quelques points d'alg{\`e}bre homologique.
\newblock \textit{T{\^o}hoku Math. J. (2)} \textbf{9} (1957) 119--221.
\mrev{0102537} \zbl{0118.26104}

\bibitem[Gro81]{Gro81}
\textsc{Gross, B.H.}
\newblock {$p$}-adic {$L$}-series at {$s=0$}.
\newblock \textit{J. Fac. Sci. Univ. Tokyo Sect. IA Math.} \textbf{28}(3)
  (1981) 979--994 (1982).
\mrev{0656068} \zbl{0507.12010}


\bibitem[Hid10]{Hid10}
\textsc{Hida, H.}
\newblock The {I}wasawa {$\mu$}-invariant of {$p$}-adic {H}ecke
  {$L$}-functions.
\newblock \textit{Ann. of Math. (2)} \textbf{172}(1) (2010) 41--137.
\mrev{2680417} \zbl{1223.11131}

\bibitem[HT94]{HidTil94}
\textsc{Hida, H.}; \textsc{Tilouine, J.}
\newblock On the anticyclotomic main conjecture for {CM} fields.
\newblock \textit{Invent. Math.} \textbf{117}(1) (1994) 89--147.
\mrev{1269427} \zbl{0819.11047}

\bibitem[HW18]{HubWas18}
\textsc{Hubbard, D.}; \textsc{Washington, L.C.}
\newblock Iwasawa invariants of some non-cyclotomic
  {$\mathbb{Z}_p$}-extensions.
\newblock \textit{J. Number Theory} \textbf{188} (2018) 18--47.
\mrev{3778621} \zbl{1455.11147}

\bibitem[Iwa73]{Iwa73a}
\textsc{Iwasawa, K.}
\newblock On $\mathbf{Z}_l$-extensions of algebraic number fields.
\newblock \textit{Ann. Math. (2)} \textbf{98} (1973) 246--326.
\mrev{0349627} \zbl{0285.12008}

\bibitem[Iwa81]{Iwa81}
---{}---.
\newblock {Riemann--Hurwitz formula and $p$-adic Galois representations for
  number fields}.
\newblock \textit{{Tohoku Math. J. (2)}} \textbf{33} (1981) 263--288.
\mrev{0624610} \zbl{0468.12004}

\bibitem[Iwa83]{Iwa83}
---{}---.
\newblock On cohomology groups of units for {${\bf Z}_{p}$}-extensions.
\newblock \textit{Amer. J. Math.} \textbf{105}(1) (1983) 189--200.
\mrev{0692110} \zbl{0525.12009}

\bibitem[Jau81]{Jau81}
\textsc{Jaulent, J.F.}
\newblock Th\'{e}orie d'{I}wasawa des tours m\'{e}tabeliennes.
\newblock \textit{Seminar on {N}umber {T}heory}, 1980--1981, 21, 1--16,  (Univ. Bordeaux
  I, Talence 1981)
\mrev{0644654} \zbl{0482.12007}


\bibitem[Mik87]{Mik87}
\textsc{Miki, H.}
\newblock On the {L}eopoldt conjecture on the {$p$}-adic regulators.
\newblock \textit{J. Number Theory} \textbf{26}(2) (1987) 117--128.
\mrev{0889379} \zbl{0621.12009}


\bibitem[NSW08]{NeuSchWin08}
\textsc{Neukirch, J.}; \textsc{Schmidt, A.}; \textsc{Wingberg, K.}
\newblock Cohomology of number fields, \textit{Grundlehren der Mathematischen
  Wissenschaften}, volume 323.
\newblock Second edition,  (Springer-Verlag, Berlin 2008).
\mrev{2392026} \zbl{1136.11001}


\bibitem[{Nuc}10]{Nuc10}
\textsc{{Nuccio~Mortarino~Majno~di~Capriglio}, F.A.E.}
\newblock Cyclotomic units and class groups in {$\mathbf{Z}_p$}-extensions of
  real abelian fields.
\newblock \textit{Math. Proc. Cambridge Philos. Soc.} \textbf{148}(1) (2010)
  93--106.
\mrev{2575375} \zbl{1279.11108}

\bibitem[Oza95]{Oza95}
\textsc{Ozaki, M.}
\newblock A note on the capitulation in {${\bf Z}_p$}-extensions.
\newblock \textit{Proc. Japan Acad. Ser. A Math. Sci.} \textbf{71}(9) (1995)
  218--219.
\mrev{1373385} \zbl{0857.11064}

\bibitem[Rot09]{Rot09}
\textsc{Rotman, J.}
\newblock An introduction to homological algebra.
\newblock Universitext, second edition,  (Springer, New York 2009).
\mrev{2455920} \zbl{1157.18001}

\bibitem[Ser62]{Ser62}
\textsc{Serre, J.P.}
\newblock {Corps locaux}, \textit{{Publications de l'Institut de Math\'ematique
  de l'Universit\'e de Nancago 8. Actualit\'es Scientifiques et
  Industrielles}}, volume 1296,  ({Hermann, Paris} 1962).
\mrev{0150130} \zbl{0137.02601}


\bibitem[Ser67]{Ser67}
---{}---.
\newblock Local class field theory.
\newblock Algebraic {N}umber {T}heory ({P}roc. {I}nstructional {C}onf.,
  {B}righton, 1965),  (Thompson, Washington, D.C. 1967) 128--161.
\mrev{0220701} \zbl{1492.11158}

\bibitem[Tat67]{Tat67}
\textsc{Tate, J.}
\newblock Global class field theory.
\newblock Algebraic {N}umber {T}heory ({P}roc. {I}nstructional {C}onf.,
  {B}righton, 1965),  (Thompson, Washington, D.C. 1967) 162--203.
\mrev{0220697} \zbl{1179.11041}


\bibitem[Vau09]{Vau09}
\textsc{Vauclair, D.}
\newblock Sur la dualit\'e et la descente d'{I}wasawa.
\newblock \textit{Ann. Inst. Fourier} \textbf{59}(2) (2009) 691--767.
\mrev{2521434} \zbl{1254.11098}

\bibitem[Was97]{Was97}
\textsc{Washington, L.C.}
\newblock Introduction to cyclotomic fields, \textit{Graduate Texts in
  Mathematics}, volume~83.
\newblock Second edition,  (Springer-Verlag, New York 1997).
\mrev{1421575} \zbl{0966.11047}

\bibitem[Wei59]{Wei59}
\textsc{Weiss, E.}
\newblock A deflation map.
\newblock \textit{J. Math. Mech.} \textbf{8} (1959) 309--329.
\mrev{0132090} \zbl{0086.02303}

\bibitem[Yam84]{Yam84}
\textsc{Yamashita, H.}
\newblock The second cohomology groups of the group of units of a {${\bf
  Z}_{p}$}-extension.
\newblock \textit{Tôhoku Math. J. (2)} \textbf{36}(1) (1984) 75--80.
\mrev{0733620} \zbl{0539.12006}

\end{thebibliography}
\end{document}